\documentclass[11pt]{article}

\usepackage{amsmath,amssymb,amsfonts}
\usepackage{pstricks,pst-node}
\usepackage{xr}

\numberwithin{equation}{section}
\newtheorem{thm}{Theorem}[section] 
\newtheorem{prp}[thm]{Proposition}
\newtheorem{lmm}[thm]{Lemma}   
\newtheorem{crl}[thm]{Corollary} 
\newtheorem{dfn}[thm]{Definition}

\def\A{\mathcal A}
\def\C{\mathbb C}
\def\cC{\mathcal C}
\def\D{\mathfrak D}
\def\cD{\mathcal D}
\def\E{\mathbb E}
\def\F{\mathfrak F}
\def\cF{\mathcal F}
\def\G{\mathfrak G}
\def\H{\mathcal H}
\def\I{\mathfrak i}
\def\L{\mathfrak L}
\def\M{\mathfrak M}
\def\cM{\mathcal M}
\def\R{\mathbb R}
\def\Q{\mathbb Q}
\def\T{\mathcal T}
\def\U{\mathcal U}
\def\fU{\mathfrak U}
\def\V{\mathcal V}
\def\W{\mathcal W}
\def\X{\mathfrak X}
\def\Z{\mathbb Z}

\def\e_ref#1{(\ref{#1})}
\def\under#1{\underline{#1}}
\def\ov#1{\overline{#1}}
\def\wt#1{\widetilde{#1}}
\def\ti#1{\tilde{#1}}

\def\lra{\longrightarrow}
\def\Lra{\Longrightarrow}

\def\Llra{\Longleftrightarrow}

\def\lan{\langle}
\def\ran{\rangle}
\def\lr#1{\lan{#1}\ran}
\def\blr#1{\big\lan{#1}\big\ran}
\def\llrr#1{\lan\!\lan{#1}\ran\!\ran}

\def\al{\alpha}
\def\be{\beta}
\def\de{\delta}
\def\ep{\epsilon}
\def\ga{\gamma}
\def\io{\iota}

\def\na{\nabla}
\def\om{\omega}
\def\si{\sigma}

\def\ups{\upsilon}
\def\ve{\varepsilon}
\def\vph{\varphi}

\def\ze{\zeta}

\def\De{\Delta}
\def\Ga{\Gamma}
\def\La{\Lambda}
\def\Om{\Omega}
\def\Si{\Sigma}

\def\Aut{\textnormal{Aut}}

\def\const{\textnormal{const}}

\def\es{\textnormal{es}}
\def\Gr{\textnormal{Gr}}
\def\GW{\textnormal{GW}}

\def\ev{\textnormal{ev}}
\def\Hom{\textnormal{Hom}}
\def\id{\textnormal{id}}
\def\Im{\textnormal{Im}\,}
\def\ind{\textnormal{ind}\,}
\def\PD{\textnormal{PD}}
\def\rk{\textnormal{rk}\,}
\def\st{\textnormal{s.t.}}

\def\P{{\mathbb P}^n}
\def\PP{{\mathbb P}^2}
\def\PPP{{\mathbb P}^3}
\def\Pf{{\mathbb P}^4}
\def\bP{\mathbb P}
\def\i{\infty}
\def\eset{\emptyset}
\def\bpar{\bar\partial}
\def\ale{\aleph}

\begin{document}

\title{On the Genus-One Gromov-Witten Invariants\linebreak 
of Complete Intersections}

\author{Jun Li\thanks{Partially supported by an NSF grant} 
~and 
Aleksey Zinger\thanks{Partially supported by an NSF Postdoctoral Fellowship}}

\date{\today}
\maketitle

\begin{abstract}
\noindent
As shown in a previous paper, certain naturally arising cones of
holomorphic vector bundle sections over the main component $\ov\M_{1,k}^0(\P,d)$
of the moduli space of stable genus-one holomorphic maps into~$\P$ have 
a well-defined euler class.
In this paper, we extend this result to moduli spaces of perturbed, in a restricted way,
$J$-holomorphic maps.
We show that euler classes of such cones relate the reduced genus-one 
Gromov-Witten invariants of complete intersections to the corresponding GW-invariants 
of the ambient projective space.
As a consequence, the standard genus-one GW-invariants of complete intersections can be 
expressed in terms of the genus-zero and genus-one GW-invariants of projective spaces.
We state such a relationship explicitly for complete-intersection threefolds.
A relationship for higher-genus invariants is conjectured as~well.
\end{abstract}

\thispagestyle{empty}

\tableofcontents

\section{Introduction}
\label{intro_sec}

\subsection{Gromov-Witten Invariants and Complete Intersections}
\label{hyperpstat_subs}

\noindent
The GW-invariants of symplectic manifolds have been an area of 
much research in the past decade.
A great deal of attention has been devoted in particular to Calabi-Yau manifolds.
These manifolds play a prominent role in theoretical physics,
and as a result physicists have made a number of important predictions concerning CY-manifolds.
Some of these predictions have been verified mathematically; others have~not.\\

\noindent
If $Y$ is a compact Kahler submanifold of the complex projective space $\P$,
one could try to compute the GW-invariants of~$Y$ by relating them to
the GW-invariants of~$\P$.
For example, suppose $Y$ is a  hypersurface in~$\P$ of degree~$a$. 
In other words, if $\ga\!\lra\!\P$ is the tautological  line bundle and
$\L\!=\!\ga^{*\otimes a}\!\lra\!\P$, then
$$Y=s^{-1}(0),$$
for some $s\!\in\!H^0(\P;\L)$ such that $s$ is transverse to the zero set.
If $g$, $k$, and~$d$ are nonnegative integers, let 
$\ov\M_{g,k}(\P,d)$ and $\ov\M_{g,k}(Y,d)$ denote the moduli spaces of stable $J_0$-holomorphic 
\hbox{degree-$d$} maps from genus-$g$ Riemann surfaces with $k$ marked points
to~$\P$ and~$Y$, respectively.
These moduli spaces determine the genus-$g$ degree-$d$ GW-invariants of $\P$ and~$Y$.\\

\noindent
By definition, the moduli space $\ov\M_{g,k}(Y,d)$ is a subset
of the moduli space $\ov\M_{g,k}(\P,d)$. In fact,
\begin{equation}\label{modulirel_e1}
\ov\M_{g,k}(Y,d)=\big\{[\cC,u]\!\in\!\ov\M_{g,k}(\P,d)\!:
s\!\circ\!u=0\in H^0\big(\cC;u^*\L\big)\big\}.
\end{equation}
Here $[\cC,u]$ denotes the equivalence class of the holomorphic map 
$u\!:\cC\!\lra\!\P$ from a genus-$g$ curve~$\cC$ with $k$ marked points.
The relationship~\e_ref{modulirel_e1} can be restated more globally as follows. 
Let 
$$\pi_{g,k}^d\!\!:\fU_{g,k}(\P,d)\lra\ov\M_{g,k}(\P,d)$$ 
be 
the semi-universal family and let
$$\ev_{g,k}^d\!:\fU_{g,k}(\P,d) \lra \P$$ 
be the natural evaluation~map.
In other words, the fiber of $\pi_{g,k}^d$ over $[\cC,u]$ is the curve~$\cC$ 
with $k$ marked points, while 
$$\ev_{g,k}^d\big([\cC,u;z]\big)=u(z)
\qquad\hbox{if}\quad z\!\in\!\cC.$$
We define a section~$s_{g,k}^d$ of the sheaf 
$\pi_{g,k*}^d\ev_{g,k}^{d*}\L\!\lra\!\ov\M_{g,k}(\P,d)$ by
$$s_{g,k}^d\big([\cC,u]\big)=[s\circ u].$$
By~\e_ref{modulirel_e1}, $\ov\M_{g,k}(Y,d)$ is the zero set of this section.\\

\noindent
The previous paragraph suggests that it should be possible
to relate the genus-$g$ degree-$d$ GW-invariants of the hypersurface~$Y$
to the moduli space $\ov\M_{g,k}(\P,d)$ in general and to the~sheaf 
$$\pi_{g,k*}^d\ev_{g,k}^{d*}\L\lra\ov\M_{g,k}(\P,d)$$ 
in particular. 
In~fact, it can be shown that
\begin{equation}\label{genus0_e}
\GW_{0,k}^Y(d;\psi)\!\equiv\!
\blr{\psi,\big[\ov\M_{0,k}(Y,d)\big]^{vir}}
=\blr{\psi\cdot e\big(\pi_{0,k*}^d\ev_{0,k}^{d*}\L\big),
\big[\ov\M_{0,k}(\P,d)\big]}
\end{equation}
for all $\psi\!\in\!H^*(\ov\M_{0,k}(\P,d);\Q)$; see~\cite{Bea} for example.
The moduli space $\ov\M_{0,k}(\P,d)$ is a smooth orbivariety  and
\begin{equation}\label{g0sheaf_e}
\pi_{0,k*}^d\ev_{0,k}^{d*}\L\lra\ov\M_{0,k}(\P,d)
\end{equation}
is a locally free sheaf, i.e.~a vector bundle.
The right-hand side of~\e_ref{genus0_e}
can be computed via the classical localization theorem of~\cite{ABo},
though the complexity of this computation increases rapidly with the degree~$d$.\\

\noindent
A {\tt hyperplane property}, i.e.~a relationship such as~\e_ref{genus0_e}, 
for positive-genus GW-invariants has been elusive since 
the early days of the Gromov-Witten theory.
If $g\!>\!0$, the sheaf 
$$\pi_{g,k*}^d\ev_{g,k}^{d*}\L \lra\ov\M_{g,k}(\P,d)$$ 
is not locally free and need not define an euler class.
Thus, the right-hand side of~\e_ref{genus0_e} may not even
make sense if $0$ is replaced by $g\!>\!0$.
Instead one might try to generalize~\e_ref{genus0_e} as
\begin{equation}\label{allgenus_e}\begin{split}
\GW_{g,k}^Y(d;\psi)
&\equiv\blr{\psi,\big[\ov\M_{g,k}(Y,d)\big]^{vir}}\\
&\stackrel{?}{=}\blr{ \psi\cdot
e\big(R^0\pi_{g,k*}^d\ev_{g,k}^{d*}\L-R^1\pi_{g,k*}^d\ev_{g,k}^{d*}\L\big),
\big[\ov\M_{g,k}(\P,d)\big]^{vir}},
\end{split}\end{equation}
where $R^i\pi_{g,k*}^d\ev_{g,k}^{d*}\L\!\lra\!\ov\M_{g,k}(\P,d)$
is the $i$th direct image sheaf.
The right-hand side of~\e_ref{allgenus_e} can be computed via 
the virtual localization theorem of~\cite{GrP1}.
However, 
$$N_1(d)\!\equiv\!\hbox{GW}_{1,0}^Y(d;1)\neq
\blr{e\big(R^0\pi_{1,0*}^d\ev_{1,0}^{d*}\L\!-\!R^1\pi_{1,0*}^d\ev_{1,0}^{d*}\L\big),
\big[\ov\M_1(\Pf,d)\big]^{vir}},$$
according to a low-degree check of~\cite{GrP2} and~\cite{K}
for a quintic threefold $Y\!\subset\!\Pf$.\\

\noindent
In this paper we prove a hyperplane property for genus-one GW-invariants.
We denote by
$$\ov\M_{1,k}^0(\P,d) \subset \ov\M_{1,k}(\P,d)$$
the closure in $\ov\M_{1,k}(\P,d)$, either in the stable-map or Zariski topology,
of the subspace
$$\M_{1,k}^0(\P,d)=\big\{[\cC,u]\!\in\!\ov\M_{1,k}(\P,d)\!:
\cC~\hbox{is smooth}\big\}.$$
If $Y\!\subset\!\P$ is a hypersurface as above, let
$$\ov\M_{1,k}^0(Y,d)=\ov\M_{1,k}(Y,d)\cap\ov\M_{1,k}^0(\P,d).$$
Since $\ov\M_{1,k}^0(\P,d)$ is a unidimensional orbi-variety, it carries a fundamental class. 
By Corollary~\ref{g1comp2-fundclass_crl} in~\cite{g1comp2}, 
$\ov\M_{1,k}^0(Y,d)$ carries a virtual fundamental class.
It can be used to define {\tt reduced} genus-one Gromov-Witten invariants:
$$\GW_{1,k}^{0;Y}(d;\psi) \equiv
\blr{\psi,\big[\ov\M_{1,k}^0(Y,d)\big]^{vir}}\in\Q,$$
where $\psi$ is a tautological (cohomology) class on $\ov\M_{1,k}^0(Y,d)$; see below.
We show in this paper that the reduced genus-one GW-invariants satisfy
a natural analogue of~\e_ref{genus0_e}.

\begin{thm}
\label{hyperp_thm}
Suppose $d$ and $a$ are positive integers, $k$ is a nonnegative integer,
$\L\!=\!\ga^{*\otimes a}\!\lra\!\P$, 
$$\pi_{1,k}^d\!:\fU_{1,k}(\P,d)\lra\ov\M_{1,k}^0(\P,d)
\quad\hbox{and}\quad
\ev_{1,k}^d\!:\fU_{1,k}(\P,d)\lra\P$$ 
are the semi-universal family and the natural evaluation~map, respectively.
If $Y\!\subset\!\P$ is a smooth degree-$a$ hypersurface, then
\begin{equation}\label{hyperp_thm_e}
\GW_{1,k}^{0;Y}(d;\psi)=
\blr{\psi\cdot e(\pi_{1,k*}^d\ev_{1,k}^{d*}\L),\big[\ov\M_{1,k}^0(\P,d)\big]}
\end{equation}
for every tautological class $\psi$ on $\ov\M_{1,k}(\P,d)$.
\end{thm}

\noindent
The tautological classes on $\ov\M_{1,k}(\P,d)$ are certain natural cohomology
classes. They include all geometric classes defined in Subsection~\ref{mainthmcrl_subs}.
We describe the space of all cohomology classes~$\psi$ to which 
Theorem~\ref{hyperp_thm} applies in Subsection~\ref{hyperppf_subs}.\\

\noindent
Implicit in the statement of Theorem~\ref{hyperp_thm} is that the euler class of
the~sheaf
\begin{equation}\label{g1sheaf_e}
\pi_{1,k*}^d\ev_{1,k}^{d*}\L\lra\ov\M_{1,k}^0(\P,d)
\end{equation}
is well-defined, even though it is not locally free.
This is the case by Theorem~\ref{g1cone-thm1} in~\cite{g1cone}.\\

\noindent
The right-hand side of~\e_ref{hyperp_thm_e} should in principle
be computable via localization directly.
However, since the space $\ov\M_{1,k}^0(\P,d)$ is not smooth
and the sheaf~\e_ref{g1sheaf_e} is not locally free,
the localization theorem of~\cite{ABo} is not immediately applicable.
A desingularization of the space~$\ov\M_{1,k}^0(\P,d)$, 
i.e.~a smooth orbivariety $\wt\M_{1,k}^0(\P,d)$ and a~map
$$\ti\pi\!:\wt\M_{1,k}^0(\P,d)\lra\ov\M_{1,k}^0(\P,d),$$
which is biholomorphic onto $\M_{1,k}^0(\P,d)$, is constructed in~\cite{VZ}.
This desingularization of $\ov\M_{1,k}^0(\P,d)$ comes with a desingularization
of the sheaf~\e_ref{g1sheaf_e}, i.e.~a vector bundle 
$$\ti\V_{1,k}^d\lra\wt\M_{1,k}^0(\P,d) \qquad\hbox{s.t.}\qquad
\ti\pi_*\ti\V_{1,k}^d=\pi_{1,k*}^d\ev_{1,k}^{d*}\L.$$
In particular,
\begin{equation}\label{design_e}
\blr{\psi\cdot e\big(\pi_{1,k*}^d\ev_{1,k}^{d*}\L\big),
\big[\ov\M_{1,k}^0(\P,d)\big]}=
\blr{\ti\pi^*\psi\cdot e(\ti\V_{1,k}^d),
\big[\wt\M_{1,k}^0(\P,d)\big]}.
\end{equation}
Since a group action on $\P$ induces actions on $\wt\M_{1,k}^0(\P,d)$
and on $\ti\V_{1,k}^d$, the localization theorem of~\cite{ABo} is directly applicable
to the right-hand side of~\e_ref{design_e}, for a natural cohomology class~$\psi$.\\ 

\noindent
By itself, Theorem~\ref{hyperp_thm} does not provide a way of computing the standard
genus-one GW-invariants of~$Y$.
However, the reduced genus-one GW-invariants capture the contribution of
$\ov\M_{1,k}^0(Y,d)$ to the standard genus-one GW-invariants.
Thus, the difference between the two invariants is completely determined
by the genus-zero invariants of~$Y$; see Subsection~\ref{g1comp2-gwdiff_subs}
and Proposition~\ref{g1comp2-bdcontr_prp} in~\cite{g1comp2}.
We give explicit formulas in some special cases in Subsection~\ref{mainthmcrl_subs} below.\\

\noindent
{\it Remark 1:}
Theorem~\ref{hyperp_thm} generalizes to arbitrary smooth complete intersections
in projective spaces.
More precisely, if 
$$\L=\ga^{*\otimes a_1}\!\oplus\!\ldots\!\oplus\!\ga^{*\otimes a_m}
\lra\P,$$
with $a_1,\ldots,a_m\!\in\!\Z^+$, $s\!\in\!H^0(\P;\L)$ is transverse 
to the zero set in~$\L$, and $Y\!=\!s^{-1}(0)$, then
\begin{equation}\label{rem1_e}
\GW_{1,k}^{0;Y}(d;\psi)
=\blr{\psi\cdot e(\pi_{1,k*}^d\ev_{1,k}^{d*}\L),\big[\ov\M_{1,k}^0(\P,d)\big]},
\end{equation}
for every geometric cohomology class $\psi$ on $\ov\M_{1,k}(\P,d)$.\\

\noindent
{\it Remark 2:} In turn, Remark~1 generalizes as follows.
Suppose $(X,\om,J)$ is a compact almost Kahler manifold,
$$A\in H_2(X;\Z)^*\!\equiv\!H_2(X;\Z)\!-\!\{0\},$$ 
$(\L,\na)\!\lra\!X$ is a complex vector bundle with connection,
and $s$ is a $\na$-holomorphic section of $\L$; see 
Subsections~\ref{eulercone_subs} and~\ref{hyperppf_subs} for terminology.
If $J$ is genus-one $A$-regular in the sense of Definition~\ref{g1comp-g1reg_dfn} in~\cite{g1comp},
$s$ is transverse to the zero set in $\L$, and $(\L,\na)$ splits into line bundles that 
are $(\om,A)$-positive in the sense of Definition~\ref{pos_dfn} below, then
\begin{equation}\label{rem2_e}\begin{split}
\GW_{1,k}^{0;Y}(A;\psi) &=\blr{\psi\cdot e(\V_{1,k}^A),\big[\ov\M_{1,k}^0(X,A;J)\big]^{vir}}\\
&\equiv\blr{\psi,\PD_{\ov\M_{1,k}^0(X,A;J)}e(\V_{1,k}^A)},
\end{split}\end{equation}
where $Y\!=\!s^{-1}(0)$, $\psi$ is a tautological class, and the cone 
$$\V_{1,k}^A\lra\ov\M_{1,k}^0(X,A;J)$$ 
is the geometric analogue of the sheaf $\pi_{1,k*}^d\ev_{1,k}^{d*}\L$.
It consists of $\na$-holomorphic sections of the vector bundle~$\L$,
as defined in Subsection~\ref{eulercone_subs} below.
By Corollary~\ref{fundclass_crl}, the Poincare dual of
its euler class is well defined as long as 
$(\L,\na)$ is a direct sum of $(\om,A)$-positive line bundles.\\

\noindent
Theorem~\ref{hyperp_thm} and Remarks~1 and~2 have a natural, but rather speculative,
generalization to higher-genus invariants.
Suppose that the main component
$$\ov\M_{g,k}^0(X,A;J)\subset \ov\M_{g,k}(X,A;J)$$
is well defined and carries a virtual fundamental class.
If so, it determines {\tt reduced} genus-$g$ GW-invariants $\GW_{g,k}^{0;Y}(A;\psi)$.
Suppose further that (the Poincare dual of) the euler of the cone
$$\V_{g,k}^A \lra \ov\M_{g,k}^0(X,A;J)$$
corresponding to the vector bundle $(\L,\na)\!\lra\!X$ is well defined.
If constructions of these objects are direct generalizations of the corresponding
constructions in~Subsection~\ref{eulercone_subs} and in \cite{g1comp}-\cite{g1comp2},
then the proof of Theorem~\ref{hyperp_thm} can be generalized to show that
\begin{equation}\label{gengenus_e}
\GW_{g,k}^{0;Y}(A;\psi)=
\blr{\psi\cdot e(\V_{g,k}^d),\big[\ov\M_{g,k}^0(X,A;J)\big]^{vir}},
\end{equation}
provided appropriate generalizations of the assumptions in Remark~2 hold.
Along with an equally speculative generalization of 
Theorem~\ref{g1comp2-gwdiff_thm} in~\cite{g1comp2} stated
in Subsection~\ref{g1comp2-gwdiff_subs} of~\cite{g1comp2}, 
\e_ref{gengenus_e}~would, if true, provide an algorithm
for computing arbitrary-genus GW-invariants of complete intersections.\\

\noindent
From the point of view of algebraic geometry as described in~\cite{BehFa} and in~\cite{LT2},
the genus-$g$ degree-$d$ GW-invariant $\GW_{g,k}^Y(d;\psi)$
is the evaluation of $\psi$ on the virtual fundamental class $[\ov\M_{g,k}(Y,d)]^{vir}$.
Using the more concrete point of view of symplectic topology as described in
\cite{FuO} and~\cite{LT1}, $\GW_{g,k}^Y(d;\psi)$
can be interpreted as the euler class of a vector bundle, albeit of an infinite-rank
vector bundle over a space of the ``same" dimension.
As in the finite-dimensional case, this euler class is the number of zeros, 
counted with appropriate multiplicities, of a transverse (multivalued, generic) section.
It is shown in~\cite{LT3} and in~\cite{Si} that the two approaches are equivalent.
In this paper, we take the latter point of view.
Similarly, we view the euler class of the sheaf~\e_ref{g1sheaf_e}
as the zero set of a generic section of its geometric analogue~$\V_{1,k}^d$ defined
in Subsection~\ref{eulercone_subs}.\\

\noindent
Theorem~\ref{hyperp_thm} and Remark~1 are special cases of Remark~2,
which is the same as Theorem~\ref{hyperpgen_thm}.
It is proved in Subsection~\ref{hyperppf_subs} by showing that 
the zero sets of two bundle sections whose cardinalities are the two expressions 
in~\e_ref{rem2_e} are the same set.
In fact, Theorem~\ref{hyperpgen_thm}, just like its genus-zero analogue, 
follows easily from definitions of the two sides in~\e_ref{rem2_e},
once it is established that these definition are well-posed.

\subsection{Cones of Holomorphic Bundle Sections}
\label{eulercone_subs}

\noindent
Let $(X,\om,J)$ be a compact almost Kahler manifold.
In other words, $(X,\om)$ is a symplectic manifold and 
$J$ is an almost complex structure on $X$ tamed by $\om$,~i.e.
$$\om(v,Jv)>0 \qquad\forall~ v\in TX\!-\!X.$$
If $g,k$ are nonnegative integers and $A\!\in\!H_2(X;\Z)$,
we denote by $\ov\M_{g,k}(X,A;J)$ the moduli space of (equivalence classes of) 
stable $J$-holomorphic maps from genus-$g$ Riemann surfaces with $k$ marked points 
in the homology class~$A$.
Let $\ov\M_{1,k}^0(X,A;J)$ be the main component of the moduli space $\ov\M_{1,k}(X,A;J)$ 
described by Definition~\ref{g1comp-degen_dfn} in~\cite{g1comp};
see also Definition~\ref{degen_dfn} below.
This closed subspace of $\ov\M_{1,k}(X,A;J)$ contains the subspace $\M_{1,k}^{0}(X,A;J)$
consisting of the stable maps $[\Si,u]$ such that the domain $\Si$
is a smooth Riemann surface.
If $J$ is sufficiently regular, $\ov\M_{1,k}^0(X,A;J)$ is the closure of
$\M_{1,k}^0(X,A;J)$ in $\ov\M_{1,k}(X,A;J)$.\\

\noindent
Suppose $\L\!\lra\!X$ is a complex line bundle and $\na$ is a connection in~$\L$.
If $(\Si,j)$ is a Riemann surface and $u\!:\Si\!\lra\!X$ is a smooth map, 
let
$$\na^u\!:\Ga(\Si;u^*\L)\lra \Ga\big(\Si;T^*\Si\!\otimes\!u^*\L\big)$$
be the pull-back of $\na$ by $u$.
If $b\!=\!(\Si,j;u)$, we define the corresponding $\bar{\partial}$-operator by
\begin{equation}\label{vdfn_e}
\bpar_{\na,b}\!:\Ga(\Si;u^*\L)\lra
\Ga\big(\Si;\La_{\I,j}^{0,1}T^*\Si\!\otimes\!u^*\L\big),
\quad
\bpar_{\na,b}\xi=\frac{1}{2}
\big(\na^u\xi+\I\na^u\xi\circ j\big),
\end{equation}
where $\I$ is the complex multiplication in the bundle $u^*\L$ and
$$\La_{\I,j}^{0,1}T^*\Si\!\otimes\!u^*\L=
\big\{\eta\!\in\!\hbox{Hom}(T\Si,u^*\L)\!:\eta\!\circ\!j=-\I\eta\big\}.$$
The kernel of $\bpar_{\na,b}$ is necessarily a finite-dimensional complex vector space.\\

\noindent
We denote by $\X_{1,k}(X,A)$ the space of all degree-$A$ stable smooth maps from 
genus-one Riemann surfaces with $k$ marked points into~$X$ and~by
$$\V_{1,k}^A\lra\X_{1,k}(X,A) $$
the cone, or the bundle of (orbi-)vector spaces, such that
$$\V_{1,k}^A\big|_{[b]}=\ker\bpar_{\na,b}\big/\Aut(b)
\qquad\forall~[b]\in\X_{1,k}(X,A).$$
The spaces  $\X_{1,k}(X,A)$ and $\V_{1,k}^A$ have natural topologies; 
see  Subsection~\ref{review_subs} below.
By Theorem~\ref{g1cone-thm1} in~\cite{g1cone}, if $(X,\om,J)$ is the complex projective
space $(\P,\om_0,J_0)$ with its standard Kahler structure and 
$(\L,\na)$ is a positive power of the hyperplane line bundle,
i.e.~the dual of the tautological line bundle,  $\ga^*$ 
with its standard connection, then the euler class~of
$$\V_{1,k}^A\lra \ov\M_{1,k}^0(X,A;J)$$
and its Poincare dual are well defined.
By Theorem~\ref{g1cone-main_thm} in~\cite{g1cone}, this is also the case if 
$J$ is an almost complex structure on $\P$ sufficiently close to~$J_0$.\\

\noindent
The argument in~\cite{g1cone} easily generalizes to all $(X,\om,J)$, $(\L,\na)$, and
$A$ such that $(\L,\na)$ is a split positive vector bundle with connection and
$J$ satisfies a certain regularity condition.
This regularity condition, which is described by 
Definition~\ref{g1comp-g1reg_dfn} in~\cite{g1comp},
implies that $\ov\M_{1,k}^0(X,A;J)$ has the expected topological structure of a 
unidimensional orbivariety. 
In this paper, we show that the {\it Poincare dual of the euler class} of $\V_{1,k}^A$ 
over $\ov\M_{1,k}^0(X,A;J)$ is well defined without any condition on~$J$,
as long as $(\L,\na)$ satisfies the requirement of Definition~\ref{pos_dfn};
see Corollary~\ref{fundclass_crl} below.

\begin{dfn}
\label{pos_dfn}
Suppose $(X,\om)$ is a symplectic manifold and $A\!\in\!H_2(X;\Z)$.
A complex line bundle $\L\!\lra\!X$ is {\tt $(\om,A)$-positive} if 
$$\lr{c_1(\L),B}>0 \qquad\forall~~ B\!\in\!H_2(X;\Z)^* ~~\st~~ 
B\!=\!A ~~\hbox{or}~~ \lr{\om,B}\!<\!\lr{\om,A}.$$\\
\end{dfn}

\noindent
We note that $\V_{1,k}^A\!\lra\!\ov\M_{1,k}^0(X,A;J)$ is {\it not} a vector bundle,
as the fibers of $\V_{1,k}^A$ are of two possible dimensions.
In Subsections~\ref{g1cone-appr_subs} and~\ref{g1cone-mainres_subs} of~\cite{g1cone},
the Poincare dual of the euler class of $\V_{1,k}^A$ is defined as the zero set
of a generic multisection~$\vph$ of $\V_{1,k}^A$ over $\ov\M_{1,k}^0(X,A;J)$.
This zero set determines a homology class in $\ov\M_{1,k}^0(X,A;J)$ 
if $\vph$ is sufficiently regular.
In Section~\ref{g1cone-main_sec} of~\cite{g1cone}, it is shown that $\V_{1,k}^A$
contains a vector subbundle of a sufficiently high rank over a neighborhood of
every stratum of $\ov\M_{1,k}^0(X,A;J)$.
The existence of such subbundles implies that regular sections of $\V_{1,k}^A$ 
exist; see Subsection~\ref{g1cone-g1conestr_subs} in~\cite{g1cone}.\\

\noindent
If $J$ does not satisfy the regularity condition of 
Definition~\ref{g1comp-g1reg_dfn} in~\cite{g1comp},
the moduli space $\ov\M_{1,k}^0(X,A;J)$ itself need not carry a fundamental class.
In this case, we cannot define the Poincare dual of the euler class of $\V_{1,k}^A$ 
as the zero set of a section of~$\V_{1,k}^A$ over $\ov\M_{1,k}^0(X,A;J)$.
On the other hand, in~\cite{g1comp2}, the definition of 
$$\ov\M_{1,k}^0(X,A;J)\subset \ov\M_{1,k}(X,A;J)$$
given in~\cite{g1comp} is generalized to define the main component 
$\ov\M_{1,k}^0(X,A;J,\nu)$ 
of the moduli space $\ov\M_{1,k}(X,A;J,\nu)$ of $(J,\nu)$-holomorphic maps
for  an {\it effectively supported perturbation} $\nu$ of the $\bpar_J$-operator;
see Definitions~\ref{pert_dfn} and~\ref{degen_dfn} below.
By Theorem~\ref{g1comp2-reg_thm} in~\cite{g1comp2}, if $\nu$ is sufficiently small and 
generic, $\ov\M_{1,k}^0(X,A;J,\nu)$ determines a rational homology class in 
a small neighborhood of $\ov\M_{1,k}^0(X,A;J)$ in $\X_{1,k}(X,A)$.
This rational homology class is independent of the choice of~$\nu$.
We will define the Poincare dual of the euler class of $\V_{1,k}^A$ 
as the zero set of a generic multisection of~$\V_{1,k}^A$ over $\ov\M_{1,k}^0(X,A;J,\nu)$.\\

\noindent
If $J$ is an almost complex structure on $X$ and $\under{J}\!\equiv\!(J_t)_{t\in[0,1]}$
is a family of almost complex structures on $X$, we denote~by 
$$\G_{1,k}^{\es}(X,A;J) \qquad\hbox{and}\qquad \G_{1,k}^{\es}(X,A;\under{J})$$
the spaces of effectively supported perturbations of the $\bpar_J$-operator on 
$\X_{1,k}(X,A)$ and of effectively supported families of perturbations of
the $\bpar_{J_t}$-operators on $\X_{1,k}(X,A)$; see Subsection~\ref{review_subs} for details.
If 
$$\bar\nu\!\equiv\!(\nu_t)_{t\in[0,1]}\in\G_{1,k}^{\es}(X,A;\under{J}),$$ 
we put
$$\ov\M_{1,k}^0(X,A;\under{J},\under{\nu})=
\big\{(t,b)\!\in\![0,1]\!\times\!\X_{1,k}(X,A)\!:
b\!\in\!\ov\M_{1,k}^0(X,A;J_t,\nu_t)\big\}.$$
We denote by $\bar\Z^+$ the set of nonnegative integers. Let
$$\dim_{1,k}(X,A;\L)=\dim_{1,k}(X,A)-2\lr{c_1(\L),A}
=2\big(\lr{c_1(TX)\!-\!c_1(\L),A}+k\big).$$

\begin{thm}
\label{main_thm}
Suppose $(X,\om,J)$ is a compact almost Kahler manifold,
$A\!\in\!H_2(X;\Z)^*$, $k\!\in\!\bar\Z^+$,
$(\L,\na)\!\lra\!X$ is an $(\om,A)$-positive line bundle with connection,
 $\V_{1,k}^A\!\lra\!\X_{1,k}(X,A)$ is the corresponding cone, and
$W$ is a neighborhood of $\ov\M_{1,k}^0(X,A;J)$ in $\X_{1,k}(X,A)$.
If $\nu\!\in\!\G_{1,k}^{\es}(X,A;J)$ is sufficiently small and generic and
$\vph$ is a generic multisection of $\V_{1,k}^A$ over $\ov\M_{1,k}^0(X,A;J,\nu)$,
then $\vph^{-1}(0)$ determines a rational homology class in~$W$.
Furthermore, if $\under{J}\!=\!(J_t)_{t\in[0,1]}$ is a family of $\om$-tamed
almost complex structures on $X$, such that $J_0\!=\!J$ and $J_t$ is sufficiently close
to~$J$ for all~$t$, $\nu_0$ and $\nu_1$ are sufficiently small generic 
effectively supported perturbations of $\bar{\partial}_{J_0}$ and $\bar{\partial}_{J_1}$,
and $\vph_0$ and $\vph_1$ are generic multisections of $\V_{1,k}^A$ over $\ov\M_{1,k}^0(X,A;J_0,\nu_0)$ and $\ov\M_{1,k}^0(X,A;J_1,\nu_1)$, 
then there exist homotopies 
$$\under\nu\!=\!(\nu_t)_{t\in[0,1]} \in \G^{\es}_{1,k}(X,A;\under{J})
\qquad\hbox{and}\qquad \Phi\in\Ga\big(
\ov\M_{1,k}^0(X,A;\under{J},\under{\nu});\V_{1,k}^A\big)$$
between $\nu_0$ and $\nu_1$ and between $\vph_0$ and $\vph_1$ 
such that $\Phi^{-1}(0)$ determines a chain in $W$~and
$$\partial\Phi^{-1}(0) = \vph_1^{-1}(0)-\vph_0^{-1}(0).$$
\end{thm}

\begin{crl}
\label{fundclass_crl}
If $(X,\om,J)$, $A$, $k$, and $(\L,\na)$ are as in Theorem~\ref{main_thm},
the cone $\V_{1,k}^A\!\lra\!\X_{1,k}(X,A)$ 
corresponding to~$(\L,\na)$ determines a well-defined homology class
$$\PD_{\ov\M_{1,k}^0(X,A;J)}e(\V_{1,k}^A)\in
H_{\dim_{1,k}(X,A;\L)}\big(\ov\M_{1,k}^0(X,A;J);\Q).$$
This class is an invariant of $(X,\om)$ and~$(\L,\na)$.
\end{crl}

\noindent
As in~\cite{g1cone}, we will describe the local structure of the cone~$\V_{1,k}^A$.
In contrast to~\cite{g1cone}, we will not construct a high-rank vector subbundle 
of $\V_{1,k}^A$ over a neighborhood of every stratum of $\ov\M_{1,k}^0(X,A;J,\nu)$.
Instead, we will use the Tietze Extension Theorem to construct  a sufficiently 
regular multisection of~$\V_{1,k}^A$. 
Its zero set determines a homology class in a small neighborhood of
$\ov\M_{1,k}^0(X,A;J)$ in the space $\X_{1,k}(X,A)$.\\

\noindent
For a generic $\nu$, $\ov\M_{1,k}^0(X,A;J,\nu)$ can be stratified by orbifolds $\U_{\al}$
of even dimensions; see Subsection~\ref{g1str_subs} and Remark~1 at the end of Subsection~\ref{top_subs}.
The main stratum of $\ov\M_{1,k}^0(X,A;J,\nu)$,
$$\M_{1,k}^0(X,A;J,\nu)\equiv  \ov\M_{1,k}^0(X,A;J,\nu)
\cap \X_{1,k}^0(X,A),$$
is of dimension $\dim_{1,k}(X,A)$, where 
$$\X_{1,k}^0(X,A)\subset \X_{1,k}(X,A)$$
is the subspace of stable maps with smooth domains.
In Subsection~\ref{g1conestr_subs}, we describe a subcone $\W_{1,k}^A$ of $\V_{1,k}^A$
such that $\W_{1,k}^A|_{\U_{\al}}$ is a smooth vector bundle for every stratum~$\U_{\al}$.
By analyzing the obstruction to extending holomorphic bundle sections from singular
to smooth domains in Section~\ref{gluing_sec}, we show that $\W_{1,k}^A$ is
a regular obstruction-free cone in the sense of Definition~\ref{coneprop_dfn}.
By Proposition~\ref{conetop_prp}, for a generic multisection~$\vph$ of
$\W_{1,k}^A\!\subset\!\V_{1,k}^A$ over $\ov\M_{1,k}^0(X,A;J,\nu)$,
$\vph|_{\U_{\al}}$ is then transverse to the zero set in~$\W_{1,k}^A|_{\U_{\al}}$.
By the rank statements of Proposition~\ref{conestr_prp}, 
$\vph^{-1}(0)$ is stratified by smooth orbifolds of even dimensions.
Furthermore, the main stratum of $\vph^{-1}(0)$ is of dimension $\dim_{1,k}(X,A;\L)$
and is contained in $\M_{1,k}^0(X,A;J,\nu)$.
We can then choose an arbitrarily small neighborhood $U$ of the boundary of $\vph^{-1}(0)$
such that
$$H_l(U;\Q)=\{0\} \qquad\forall~l\ge \dim_{1,k}(X,A;\L)-1.$$
Since $\vph^{-1}(0)\!-\!U$ is compact,
via the pseudocycle construction of Chapter~7 in~\cite{McSa} and
Section~1 of~\cite{RT}, $\vph^{-1}(0)$ determines a homology class
\begin{equation*}\begin{split}
\big[\vph^{-1}(0)\big]&\in 
H_{\dim_{1,k}(X,A;\L)}(W,U;\Q)\\
&\qquad\qquad \approx H_{\dim_{1,k}(X,A;\L)}(W;\Q).
\end{split}\end{equation*}
The second part of Theorem~\ref{main_thm} is a parametrized version of this construction.
Corollary~\ref{fundclass_crl} is an immediate consequence of Theorem~\ref{main_thm};
see also Remark~2 in Subsection~\ref{g1comp2-pertmaps_subs} of~\cite{g1comp2} and
the comments at the end of Subsection~\ref{g1comp2-res_subs} in~\cite{g1comp2}.\\

\noindent
The statement of Corollary~\ref{fundclass_crl} is not needed to show that 
the expressions on the right-hand sides of~\e_ref{hyperp_thm_e} and~\e_ref{rem1_e}
are well defined, as this is the case by Theorem~\ref{g1cone-thm1} in~\cite{g1cone}.
However, the detailed statement of Theorem~\ref{main_thm} is useful for proving  
Theorem~\ref{hyperp_thm} and its generalizations in Remarks~1 and~2 
whenever $Y$ is not a Fano complete intersection.
If $Y$ is Fano, Theorem~\ref{hyperp_thm} can be obtained from~\cite{g1cone}
by working just with $J$-holomorphic, instead of $(J,\nu)$-holomorphic, maps.\\

\noindent
{\it Remark:} If $\L$ is a direct sum of $(\om,A)$-positive line bundles,
the Poincare dual of the euler class of the corresponding cone is defined to be 
the cap product of the Poincare duals of the euler classes of the cones corresponding 
to the component line bundles.
The cap product can be defined by intersecting pseudocycle representatives for
the above homology classes;
see Subsection~\ref{g1cone-appr_subs} in~\cite{g1cone}.

\subsection{Some Special Cases}
\label{mainthmcrl_subs}

\noindent
Proposition~\ref{g1comp2-bdcontr_prp} in~\cite{g1comp2} implies that
the difference between the standard and reduced genus-one invariants of 
a symplectic manifold $(Y,\om)$ is a combination of the genus-zero invariants of~$Y$.
The exact form of this combination can be determined in each specific case.\\

\noindent
If $(Y,\om,J)$ is an almost Kahler manifold, for each $l\!=\!1,\ldots,k$ let
$$\ev_l\!: \ov\M_{g,k}(Y,A;J)\lra Y, \qquad 
\big[\Si,y_1,\ldots,y_k;u\big] \lra  u(y_l),$$
be the evaluation map at the $l$th marked point.
We will call a cohomology class $\psi$  on $\ov\M_{g,k}(Y,A;J)$ {\tt geometric}
if $\psi$ is a product of the classes $\ev_l^*\mu_l$ for $\mu_l\!\in\!H^*(Y;\Z)$. 
By Theorem~\ref{g1comp2-gwdiff_thm} in~\cite{g1comp2}, if $A\!\in\!H_2(Y;\Z)^*$, 
then
\begin{equation}\label{gwdiff_e}
\GW_{1,k}^Y(A;\psi)-\GW_{1,k}^{0;Y}(A;\psi)
=\begin{cases}
0,& \hbox{if}~\dim_{\R}Y\!=\!4;\\
\frac{2-\lr{c_1(TY),A}}{24}\GW_{0,k}^Y(A;\psi),& \hbox{if}~\dim_{\R}Y\!=\!6,
\end{cases}\end{equation}
for every geometric cohomology class $\psi$ on $\ov\M_{1,k}(Y,A;J)$.\\

\noindent 
In the rest of this subsection, we discuss some implications of Theorem~\ref{hyperp_thm} 
and Remarks~1 and~2, combined with~\e_ref{gwdiff_e}, 
focusing on Calabi-Yau complete-intersection threefolds.
We note that if $Y$ is a Calabi-Yau threefold, then the expected dimension of 
the moduli space $\ov\M_{g,0}(Y,A;J)$ is zero for every $g$ and~$A$.\\

\noindent
With notation as in Theorem~\e_ref{main_thm}, if  $a\!=\!5$, $Y$ is a quintic threefold.
It can easily seen that $c_1(TY)\!=\!0$.
Let
$$N_g(d)=\GW_{g,0}^Y(d;1).$$
Theorem~\ref{main_thm} and equation~\e_ref{gwdiff_e} then give

\begin{crl}
\label{cy_crl}
Suppose $d$ is a positive integer, $\L\!=\!\ga^{*\otimes 5}\!\lra\!\Pf$, and
$$\pi_1^d\!:\fU_1(\Pf,d)\lra\ov\M_1^0(\Pf,d) \quad\hbox{and}\quad
\ev_1^d\!:\fU_1(\Pf,d)\lra\Pf$$ 
are the semi-universal family and the natural evaluation~map, respectively.
If $Y\!\subset\!\Pf$ is a smooth quintic threefold,
\begin{equation}\label{cy_crl_e}
N_1(d)=\frac{1}{12}N_0(d)+
\blr{e\big(\pi_{1*}^d\ev_1^{d*}\L\big),\big[\ov\M_1^0(\Pf,d)\big]}.
\end{equation}\\
\end{crl}

\noindent
The middle number in~\e_ref{cy_crl_e} can be computed using~\e_ref{genus0_e}.
This has been done for every~$d$ in \cite{Ber}, \cite{Ga1}, \cite{Gi}, \cite{Le}, and~\cite{LLY}.
As mentioned in Subsection~\ref{hyperpstat_subs}, the last number in~\e_ref{cy_crl_e}
can be computed, for each given~$d$, via the classical localization theorem of~\cite{ABo}.
Similarly to the genus-zero case, the complexity of computing the last term 
in~\e_ref{cy_crl_e} increases rapidly with the degree~$d$.
Its low-degree values are shown in the second row of Table~\ref{main_table}.
The numbers $n_1(d)$ that appear in the last row of this table are defined~by
$$N_0(d)=\sum_{k|d}\frac{n_0(d/k)}{k^3},\qquad
N_1(d)=\frac{1}{12}\sum_{k|d}\frac{n_0(d/k)}{k}
+\sum_{k|d}\frac{n_1(d/k)}{k}.$$
The numbers $n_0(d)$ and $n_1(d)$ are of importance in theoretical physics.
Conjecturally, $n_g(d)$~is a count of $J$-holomorphic degree-$d$ genus-$g$ curves in~$Y$
for a generic almost complex structure~$J$ on~$Y$.
The low-degree numbers~$n_1(d)$ obtained from Theorem~\ref{main_thm}
agree with those predicted by physicists;
see Table~1 in~\cite{BCOV}.\\

\begin{table}
\begin{center}
\begin{tabular}{||c|c|c|c|c||}
\hline\hline
$d$&                                  1&   2& 3& 4\\
\hline
$\lan\ldots\ran$&                     0&  {\small$\frac{2,875}{32}$}&  
   {\small$\frac{49,355,000}{81}$}& {\small$\frac{952,691,384,375}{256}$}\\
\hline 
$N_1(d)$&    {\small$\frac{2,875}{12}$}& {\small$\frac{407,125}{8}$}&  
    {\small$\frac{243,388,750}{9}$}&    {\small$\frac{366,163,353,125}{16}$}\\
\hline 
$n_1(d)$&                             0&   0& 609,250&  3,721,431,625\\
\hline\hline
\end{tabular}
\caption{The Low-Degree GW-Invariants of a Quintic Threefold}
\label{main_table}
\end{center}
\end{table}

\noindent
With notation as in Remark~1 in Subsection~\ref{hyperpstat_subs},
if $a_1\!+\!\ldots\!+\!a_n\!=\!n\!+\!1$ and $Y$ is a corresponding complete intersection,
then $Y$ is a Calabi-Yau threefold.
Let
$$N_g^Y(d)=\GW_{g,0}^Y(d;1).$$
The identities in Remark~1 and in~\e_ref{gwdiff_e} then give
$$N_1^Y(d)=\frac{1}{12}N_0^Y(d)+
\blr{e\big(\pi_{1*}^d\ev_1^{d*}\L\big),\big[\ov\M_1^0(\P,d)\big]}.$$
Once again, both terms on the right-hand side are computable via 
\e_ref{genus0_e} and the classical localization theorem.\\

\noindent
In the more general case of Remark~2 in Subsection~\ref{hyperpstat_subs},
$Y$ is a Calabi-Yau threefold if
$$c_1(\L)-c_1(TX)=0 \qquad\hbox{and}\qquad \dim_{\R}X-2\rk_{\C}\L=6.$$
In such a case,
$$N_1^Y(A)=\frac{1}{12}N_0^Y(A)+
\blr{e(\V_{1,k}^A),\big[\ov\M_1^0(X,A;J)\big]},$$
where $\V_{1,k}^A$ is the cone corresponding to $(\L,\na)$
and $N_g^Y(A)=\GW_{g,0}^Y(A;1)$.\\

\noindent
Two completely different approaches to computing positive-genus GW-invariants
of complete intersections have been proposed in~\cite{Ga2} and~\cite{MaP}.
Both approaches use degenerations and relative Gromov-Witten invariants.
The first approach can be used to compute the genus-one and -two GW-invariants
of a quintic threefold.
The latter can be used to compute arbitrary-genus GW-invariants of 
a quintic threefold as well as of some other low-degree low-dimensional 
complete intersections.
In contrast, Theorem~\ref{hyperp_thm} above and Proposition~\ref{g1comp2-bdcontr_prp} 
in~\cite{g1comp2} are at the present restricted to genus-one GW-invariants only, 
but are applicable to arbitrary complete intersections.\\

\noindent
We would like to thank D.~Maulik, R.~Pandharipande, G.~Tian,
and R.~Vakil for a number of helpful discussions.

\section{Hyperplane Property for Genus-One GW-Invariants}
\label{hyperpthm_sec}

\subsection{Review of Definitions}
\label{review_subs}

\noindent
Suppose $X$ is a compact manifold,
$A\!\in\!H_2(X;\Z)$, and $g,k\!\in\!\bar\Z^+$.
Let $\X_{g,k}(X,A)$  denote the space of equivalence classes
of stable smooth maps $u\!:\Si\!\lra\!X$ from genus-$g$ Riemann surfaces
with $k$~marked points, which may have simple nodes, to~$X$ of degree~$A$, i.e. 
$$u_*[\Si]=A\in H_2(X;\Z).$$
The spaces $\X_{g,k}(X,A)$ are topologized using $L^p_1$-convergence on compact subsets 
of smooth points of the domain and certain convergence requirements near the nodes;
see Section~3 in~\cite{LT1}.
Here and throughout the rest of the paper, $p$~denotes a real number 
greater than~two.
The spaces $\X_{g,k}(X,A)$ can be stratified by the smooth infinite-dimensional orbifolds
$\X_{\T}(X)$ of stable maps from domains of the same geometric type and with
the same degree distribution between the components of the domain;
see Subsections~\ref{notation0_subs} and~\ref{notation1_subs}.
The closure of the main stratum, $\X_{g,k}^0(X,A)$, is $\X_{g,k}(X,A)$.\\

\noindent
If $J$ is an almost complex structure on $X$, let 
$$\Ga_{g,k}^{0,1}(X,A;J)\!\lra\!\X_{g,k}(X,A)$$
be the bundle of $(TX,J)$-valued $(0,1)$-forms. 
In other words, the fiber of $\Ga_{g,k}^{0,1}(X,A;J)$ over a point
$[b]\!=\![\Si,j;u]$ in $\X_{g,k}(X,A)$ is the space
$$\Ga_{g,k}^{0,1}(X,A;J)\big|_{[b]}=\Ga^{0,1}(b;J)\big/\hbox{Aut}(b),
\quad\hbox{where}\quad
\Ga^{0,1}(b;J)=\Ga\big(\Si;\La_{J,j}^{0,1}T^*\Si\!\otimes\!u^*TX\big).$$
Here $j$ is the complex structure on $\Si$, the domain of the smooth map~$u$.
The bundle \hbox{$\La_{J,j}^{0,1}T^*\Si\!\otimes\!u^*TX$} over $\Si$
consists of $(J,j)$-antilinear homomorphisms:
$$\La_{J,j}^{0,1}T^*\Si\!\otimes\!u^*TX=\big\{
\eta\!\in\!\hbox{Hom}(T\Si,u^*TX)\!:J\!\circ\!\eta\!=\!-\eta\!\circ\!j\big\}.$$
The total space of the bundle $\Ga_{g,k}^{0,1}(X,A;J)\!\lra\!\X_{g,k}(X,A)$ 
is topologized using $L^p$-convergence on compact subsets of smooth points of the domain
and certain convergence requirements near the nodes.
The restriction of $\Ga_{g,k}^{0,1}(X,A;J)$ to each stratum
$\X_{\T}(X)$ is a smooth vector orbibundle of infinite rank.\\

\noindent
We define a continuous section of the bundle 
$\Ga_{g,k}^{0,1}(X,A;J)\!\lra\!\X_{g,k}(X,A)$  by
$$\bar{\partial}_J\big([\Si,j;u]\big) = \bar{\partial}_{J,j}u
= \frac{1}{2}\big(du+J\!\circ\!du\!\circ\!j\big).$$
By definition, the zero set of this section is the moduli space $\ov\M_{g,k}(X,A;J)$.
The restriction of~$\bpar_J$ to each stratum of $\X_{g,k}(X,A)$ is smooth.
The section $\bpar_J$ 
is Fredholm, i.e.~the linearization of its restriction to every stratum $\X_{\T}(X)$
has finite-dimensional kernel and cokernel
at every point of $\bar{\partial}_J^{-1}(0)\!\cap\!\X_{\T}(X)$.
The index of the linearization of~$\bar{\partial}_J$ at an element of 
$\M_{g,k}^0(X,A;J)$ is the expected dimension of the moduli space~$\ov\M_{g,k}(X,A;J)$,
$$\dim_{g,k}(X,A) \equiv 2\big(\lr{c_1(TX),A}+(1\!-\!g)(n\!-\!3)+k\big),
\qquad\hbox{where}\qquad 2n=\dim_{\R}\!\!X.$$ 
This is the dimension of the cycle
$$\ov\M_{g,k}(X,A;J,\nu)\equiv 
\big\{\bar{\partial}_J\!+\!\nu\big\}^{-1}(0)$$
for a small generic multivalued perturbation 
$$\nu\in \G_{g,k}^{0,1}(X,A;J)\equiv
\Ga\big(\X_{g,k}(X,A),\Ga_{g,k}^{0,1}(X,A;J)\big)$$ 
of $\bar{\partial}_J$, where $\G_{g,k}^{0,1}(X,A;J)$
is the space of all continuous multisections $\nu$ of $\Ga_{g,k}^{0,1}(X,A;J)$
such that the restriction of $\nu$ to each stratum $\X_{\T}(X)$ is smooth.
We use the term {\it multisection}, or {\it multivalued section},  
of a vector orbi-bundle as defined in Section~3 of~\cite{FuO}.
Since the moduli space $\ov\M_{g,k}(X,A;J)$ is compact, so
is $\ov\M_{g,k}(X,A;J,\nu)$ if $\nu$ is sufficiently small.\\

\noindent
An element $[\Si;u]$ of $\X_{1,k}(X,A)$ is an equivalence class
of pairs consisting of a prestable genus-one Riemann surface $\Si$
and a smooth map $u\!:\Si\!\lra\!X$.
The prestable surface~$\Si$ is a union of the principal component(s)~$\Si_P$,
which is either a smooth torus or a circle of spheres, and
trees of rational bubble components, which together will be denoted by~$\Si_B$.
Let
$$\X_{1,k}^{\{0\}}(X,A)=
\big\{[\Si;u]\!\in\!\X_{1,k}(X,A)\!: u_*[\Si_P]\neq0\in H_2(X;\Z)\big\}.$$
Suppose 
\begin{equation}\label{compsubs_e1}
[\Si;u]\in\X_{1,k}(X,A)-\X_{1,k}^{\{0\}}(X,A),
\end{equation}
i.e.~the degree of $u|_{\Si_P}$ is zero.
Let $\chi^0(\Si;u)$ be the set of components $\Si_i$ of $\Si$ such that 
for every bubble component $\Si_h$ that lies between $\Si_i$ and $\Si_P$,
including $\Si_i$ itself, the degree of $u|_{\Si_h}$ is zero.
The set $\chi^0(\Si;u)$ includes the principal component(s) of~$\Si$.
We give an example of the set  $\chi^0(\Si;u)$ in Figure~\ref{chi_fig}.
In this figure, we show the domain~$\Si$ of the stable map $(\Si;u)$
and shade the components of the domain on which the degree of the map~$u$ is not zero.
Let
$$\Si_u^0=\bigcup_{i\in\chi^0(\Si;u)}\!\!\!\!\!\!\Si_i.$$
Every bubble component $\Si_i\!\subset\!\Si_B$ is a sphere and has 
a distinguished singular point, which will be called the {\it attaching node of~$\Si_i$}.
This is the node of $\Si_i$ that lies either on~$\Si_P$
or on a bubble $\Si_h$ that lies between $\Si_i$ and~$\Si_P$.
We denote by $\chi(\Si;u)$ the set of bubble components $\Si_i$ such that
the attaching node of $\Si_i$ lies on $\Si_u^0$ and the degree of $u|_{\Si_i}$ is not zero.

\begin{dfn}
\label{pert_dfn}
Suppose $(X,\om)$ is a compact symplectic manifold and $\under{J}\!\equiv\!(J_t)_{t\in[0,1]}$ 
is a $C^1$-continuous family of $\om$-tamed almost structures on~$X$.
A continuous family of multisections $\under{\nu}\!\equiv\!(\nu_t)_{t\in[0,1]}$,
with $\nu_t\!\in\!\G_{1,k}^{0,1}(X,A;J_t)$ for all $t\!\in\![0,1]$, 
is  {\tt effectively supported} if for every element 
$$b\!\equiv\![\Si;u] \in \X_{1,k}(X,A)\!-\!\X_{1,k}^{\{0\}}(X,A)$$
there exists a neighborhood $\W_b$ of $\Si_u^0$ in a semi-universal family of
deformations for~$b$ such that
$$\nu_t(\Si';u')\big|_{\Si'\cap\W_b}=0 \qquad
\forall~~  [\Si';u']\in\X_{1,k}(X,A),~t\!\in\![0,1].$$\\
\end{dfn}

\noindent
If $b\!=\![\Si;u]$ is an element of $\X_{1,k}(X,A)$,
a {\tt semi-universal universal family of deformations} 
for~$b$ is  a fibration
$$\si_b\!: \ti\U_b\lra \De_b$$
such that $\De_b/\Aut(b)$ is a neighborhood of $b$ in $\X_{1,k}(X,A)$
and the fiber of $\si_b$ over a point $[\Si';u']$ is~$\Si'$.
If $\under{J}\!\equiv\!(J_t)_{t\in[0,1]}$ is a continuous family of $\om$-tamed almost 
structures on~$X$, we denote the space of effectively supported families
$\under{\nu}$ as in Definition~\ref{pert_dfn} by $\G^{\es}_{1,k}(X,A;\under{J})$.
Similarly, if $J$ is an almost complex structure on $X$, let 
$\G^{\es}_{1,k}(X,A;J)$ be the subspace of elements $\nu$ of $\G^{0,1}_{1,k}(X,A;J)$ 
such that the family $\nu_t\!=\!\nu$ is effectively supported.\\

\begin{figure}
\begin{pspicture}(-1.1,-1.8)(10,1.25)
\psset{unit=.4cm}
\psellipse(8,-1.5)(1.5,2.5)
\psarc[linewidth=.05](6.2,-1.5){2}{-30}{30}\psarc[linewidth=.05](9.8,-1.5){2}{150}{210}
\pscircle[fillstyle=solid,fillcolor=gray](5.5,-1.5){1}\pscircle*(6.5,-1.5){.2}
\pscircle[fillstyle=solid,fillcolor=gray](3.5,-1.5){1}\pscircle*(4.5,-1.5){.2}
\pscircle(10.5,-1.5){1}\pscircle*(9.5,-1.5){.2}
\pscircle[fillstyle=solid,fillcolor=gray](11.91,-.09){1}\pscircle*(11.21,-.79){.2}
\pscircle[fillstyle=solid,fillcolor=gray](11.91,-2.91){1}\pscircle*(11.21,-2.21){.2}
\rput(8.5,1.5){$h_0$}\rput(5.5,0){$h_1$}\rput(3.5,0){$h_2$}
\rput(10.3,0){$h_3$}\rput(13.5,0.1){$h_4$}\rput(13.5,-2.9){$h_5$}
\rput(8,-5){\small ``tacnode"}
\pnode(8,-5){A1}\pnode(6.5,-1.5){B1}
\ncarc[nodesep=.35,arcangleA=-25,arcangleB=-15,ncurv=1]{->}{A1}{B1}
\pnode(8,-4.65){A2}\pnode(10.3,-1.5){B2}
\ncarc[nodesep=0,arcangleA=40,arcangleB=30,ncurv=1]{-}{A2}{B2}
\pnode(11,-.95){B2a}\pnode(11.02,-2.02){B2b}
\ncarc[nodesep=0,arcangleA=0,arcangleB=10,ncurv=1]{->}{B2}{B2a}
\ncarc[nodesep=0,arcangleA=0,arcangleB=10,ncurv=1]{->}{B2}{B2b}
\rput(25,-1.5){\begin{tabular}{l}$\chi^0(\Si;u)\!=\!\{h_0,h_3\}$\\
${}$\\ $\chi(\Si;u)\!=\!\{h_1,h_4,h_5\}$\end{tabular}}
\end{pspicture}
\caption{An Illustration of Definition~\ref{degen_dfn}} 
\label{chi_fig}
\end{figure}

\noindent
Suppose $\nu\!\in\!\G^{\es}_{1,k}(X,A;J)$ and $[\Si;u]$ is an element
of $\ov\M_{1,k}(X,A;J,\nu)$ as in Definition~\ref{pert_dfn}.
Since $\Si_i\!\subset\!\Si_B$ is a sphere, we can represent this element 
by a pair $(\Si;u)$ such that
the attaching node of every bubble component $\Si_i\!\subset\!\Si_B$
is the south pole, or the point $\i\!=\!(0,0,-1)$, of $S^2\!\subset\!\R^3$.
Let $e_{\i}\!=\!(1,0,0)$ be a nonzero tangent vector to $S^2$ at the south pole.
If $i\!\in\!\chi(\Si;u)$, we put
$$\cD_i(\Si;u)= d\big\{u|_{\Si_i}\big\}\big|_{\i}e_{\i}
\in T_{u|_{\Si_i}(\i)}X.$$
We note that $u|_{\Si_u^0}$ is a degree-zero holomorphic map and thus constant.
Thus, $u$  maps the attaching nodes of all elements of $\chi(\Si;u)$
to the same point in~$X$.

\begin{dfn}
\label{degen_dfn}
Suppose $(X,\om,J)$ is a compact almost Kahler manifold, $A\!\in\!H_2(X;\Z)^*$, and 
\hbox{$k\!\in\!\bar{\Z}^+$}.
If $\nu\!\in\!\G^{\es}_{1,k}(X,A;J)$ is an effectively supported perturbation
of the $\bar{\partial}_J$-operator,
the {\tt main component} of the space 
$\ov\M_{1,k}(X,A;J,\nu)$ is the subset $\ov\M_{1,k}^0(X,A;J,\nu)$ 
consisting of the elements $[\Si;u]$ of $\ov\M_{1,k}(X,A;J,\nu)$
such~that\\
${}\quad$ (a) the degree of $u|_{\Si_P}$ is not zero, or\\
${}\quad$ (b) the degree of $u|_{\Si_P}$ is zero and 
$\dim_{\C}\text{Span}_{(\C,J)}\{\cD_i(\Si;u)\!:
i\!\in\!\chi(\Si;u)\}<|\chi(\Si;u)|$.
\end{dfn}

\noindent
We note that $\ov\M_{1,k}^0(X,A;J,\nu)$ is a compact  space by
Theorem~\ref{g1comp2-comp_thm} of~\cite{g1comp2} if $\nu$ is effectively supported
and sufficiently small.
For a generic effectively supported~$\nu$, $\ov\M_{1,k}^0(X,A;J,\nu)$
determines a homology class of the expected dimension
in a small neighborhood of $\ov\M_{1,k}^0(X,A;J)$
in $\X_{1,k}(X,A)$ which is independent of $\nu$ and~$J$;
see Theorem~\ref{g1comp2-reg_thm} and Corollary~\ref{g1comp2-fundclass_crl} in~\cite{g1comp2}.\\

\noindent
If $X$, $A$, $g$, and $k$ are as above and $(\L,\na)\!\lra\!X$ 
is a vector bundle with connection, we denote~by
$$\Ga_{g,k}(\L,A)\!\lra\!\X_{g,k}(X,A)$$ 
the cone such that the fiber of $\Ga_{g,k}(\L,A)$ over $[b]\!=\![\Si;u]$ in 
$\X_{g,k}(X,A)$ is the Banach space
$$\Ga_{g,k}(\L,A)\big|_{[b]}=\Ga(b;\L)\big/\hbox{Aut}(b),
\qquad\hbox{where}\qquad \Ga(b;\L)=L^p_1(\Si;u^*\L).$$
The topology on the total space of $\Ga_{g,k}(\L,A)$ is defined analogously
to the topology on 
$\Ga_{g,k}(TX,A)$ of Section~3 in~\cite{LT1}.
Let
$$\V_{g,k}^A=\big\{[b,\xi]\!\in\!\Ga_{g,k}(\L,A)\!: [b]\!\in\!\X_{g,k}(X,A);
~ \xi\!\in\!\ker\bpar_{\na,b}\!\subset\Ga_{g,k}(b;\L)\big\}
\subset\Ga_{g,k}(\L,A).$$
The cone $\V_{g,k}^A\!\lra\!\X_{g,k}(X,A)$ inherits 
its topology from~$\Ga_{g,k}(\L,A)$.

\subsection{Statement and Proof of Hyperplane Property}
\label{hyperppf_subs}

\noindent
We will call a cohomology class $\psi$ on $\X_{1,k}(X,A)$ {\tt tautological}
if there exists a vector bundle 
$$\W\lra\X_{1,k}(X,A)$$ 
such that $\W|_{\X_{\T}(X)}$ is smooth for every stratum $\X_{\T}(X)$ of $\X_{1,k}(X,A)$
and $\psi\!=\!e(\W)$.\\

\noindent
If $(X,J)$ is an almost complex manifold and
$(\L,\na)\!\lra\!X$ is a complex vector bundle with connection,
we will call a section $s$ of $\L$ {\tt $\na$-holomorphic} if
$$\bpar_{\na}s \equiv \frac{1}{2}\big(\na s+\I\na s\circ J\big)=0.$$

\begin{thm}
\label{hyperpgen_thm}
Suppose $(X,\om,J)$ is a compact almost Kahler manifold, $A\in H_2(X;\Z)^*$, $k\!\in\!\Z^+$, 
$(\L,\na)\!\lra\!X$ is a complex vector bundle with connection,
and $s$ is a $\na$-holomorphic section of $\L$ such that
$J$ is genus-one $A$-regular in the sense of Definition~\ref{g1comp-g1reg_dfn} in~\cite{g1comp},
$s$ is transverse to the zero set in~$\L$, and $(\L,\na)$ splits into $(\om,A)$-positive line
bundles.
If $\V_{1,k}^A\!\lra\!\X_{1,k}(X,A)$ is the cone corresponding to $(\L,\na)$ 
and $Y\!=\!s^{-1}(0)$, 
\begin{equation}\label{hyperpgen_thm_e}
\GW_{1,k}^{0;Y}(A;\psi) =\blr{\psi,\PD_{\ov\M_{1,k}^0(X,A;J)}e(\V_{1,k}^A)}
\end{equation}
for every tautological class $\psi$ on $\X_{1,k}(X,A)$.
\end{thm}

\noindent
Since $J$ is genus-one $A$-regular, $\ov\M_{1,k}^0(X,A;J)$ has the expected 
structure of a topological orbivariety.
By the generalization of the proof of the regularity statement of 
Theorem~\ref{g1comp-reg_thm} in~\cite{g1comp} analogous to
Subsection~\ref{g1comp2-comp1prp_subs} in~\cite{g1comp2}, for all 
$\nu\!\in\!\G_{1,k}^{\es}(X,A;J)$ sufficiently small 
$\ov\M_{1,k}^0(X,A;J,\nu)$ also has the expected structure of a topological orbivariety.
In particular, it is stratified by smooth orbifolds of even dimensions as described
in Subsection~\ref{g1str_subs} below.
We will call $\nu\!\in\!\G_{1,k}^{\es}(X,A;J)$ {\tt $(\na,s)$-compatible} if
$$\na s|_u\circ\nu(\Si;u)=0 \qquad\forall\, [\Si;u]\!\in\!\X_{1,k}(X,A).$$
We note that if $\nu$ is $(\na,s)$-compatible, then the map 
$$(\Si;u)\lra s_{1,k}^A(\Si;u)\!\equiv\!s\!\circ\!u
\in\Ga(\Si;u^*\L)$$
defines a continuous section of the cone $\V_{1,k}^A$ over 
$\ov\M_{1,k}^0(X,A;J)$.\\

\noindent
Since the $\na$-holomorphic section~$s$ is transverse to the zero set in~$\L$,
the $(\I,J)$-linear map
$$\na s\!\!:TX\lra\L$$
is surjective along $Y\!=\!s^{-1}(0)$.
Let $U_s$ be a small neighborhood of~$Y$ in~$X$ such that 
$\na s$ is surjective over~$U_s$.
The kernel of $\na s$ over~$U_s$ is then a complex subbundle
of $(TX,J)|_{U_s}$, which restricts to $TY$ along~$Y$.
We denote this subbundle by $\ti{T}Y$.
If $\nu\!\in\!\G_{1,k}^{\es}(X,A;J)$ is such~that for all $[\Si;u]\!\in\!\X_{1,k}(X,A)$
$$\nu(\Si;u)~
\begin{cases}
\in\Ga(\Si;\La^{0,1}_{J,j}T^*\Si\!\otimes\!\ti{T}Y), 
&\hbox{if}~u(\Si)\!\subset\!U_s;\\
=0,&\hbox{otherwise},
\end{cases}$$
then $\nu$ is $(\na,s)$-compatible.
Thus, every element $\nu_Y\!\in\!\G_{1,k}^{\es}(Y,A;J)$ can extended to 
a $(\na,s)$-compatible element $\nu$ of $\G_{1,k}^{\es}(X,A;J)$.
Furthermore, if $\nu_Y$ is a small, then $\nu$ can also be chosen to be small.\\

\noindent
For a small generic $\nu_Y\!\in\!\G_{1,k}^{\es}(Y,A;J)$, 
$\ov\M_{1,k}^0(Y,A;J,\nu_Y)$ is stratified by smooth orbifolds of even dimensions
so that the largest-dimensional stratum is $\M_{1,k}^0(Y,A;J,\nu_Y)$ and
$$\dim\M_{1,k}^0(Y,A;J,\nu_Y)=\dim_{1,k}(Y,A).$$
Let $\nu$ be an extension of $\nu_Y$ to a small
$(\na,s)$-compatible element of $\G_{1,k}^{\es}(X,A;J)$.
Suppose 
$$\W\lra\X_{1,k}(X,A)$$ 
is a complex vector bundle of rank $\dim_{1,k}(Y,A)/2$
as in the first paragraph of this subsection.
Choose a section $f$ of $\W$ over $\ov\M_{1,k}^0(X,A;J,\nu)$ such that 
$f|_{\U_{\al}}$ is transverse to the zero set in $\W|_{\U_{\al}}$ for
every stratum $\U_{\al}$ of~$\ov\M_{1,k}^0(X,A;J,\nu)$ and of~$\ov\M_{1,k}^0(Y,A;J,\nu_Y)$.
Then,
\begin{gather}
f^{-1}(0) \!\cap\!  \ov\M_{1,k}^0(Y,A;J,\nu_Y) \subset \M_{1,k}^0(Y,A;J,\nu_Y)
\qquad\hbox{and}\notag\\
\label{hyperpgen_e3b}
\GW_{1,k}^{0;Y}(A;\psi) = ~^{\pm}\big|f^{-1}(0)\!\cap\!\ov\M_{1,k}^0(Y,A;J,\nu_Y)\big|.\\
\end{gather}

\noindent
On the other hand, since $\nu$ is $(\na,s)$-compatible,
$s_{1,k}^A$ is a section of
$$\V_{1,k}^A \lra \ov\M_{1,k}^0(X,A;J,\nu).$$
Furthermore, 
\begin{gather}
\big\{s_{1,k}^A\big\}^{-1}(0)=\ov\M_{1,k}^0(X,A;J,\nu)\!\cap\!\X_{1,k}(Y,A)
\equiv \ov\M_{1,k}^0(Y,A;J,\nu_Y);\notag\\
\label{hyperpgen_e5b}
\Lra\qquad f^{-1}(0)\!\cap\!\big\{s_{1,k}^A\big\}^{-1}(0)
= f^{-1}(0)\!\cap\!\ov\M_{1,k}^0(Y,A;J,\nu_Y) \subset \M_{1,k}^0(X,A;J,\nu).
\end{gather}
Note that if $[b]\!=\![\Si;u]\!\in\!\M_{1,k}^0(Y,A;J,\nu_Y)$,
\begin{equation}\label{hyperpgen_e6a}\begin{split}
\ker ds_{1,k}^A|_b 
& =\big\{\xi\!\in\!\ker D_{J,\nu;b}\!: \na s|_u\!\circ\!\xi\!=\!0\big\}\\
& =\ker D_{J,\nu;b}\!\cap\!\Ga(\Si;u^*TY) =\ker D_{J|_Y,\nu_Y;b},
\end{split}\end{equation}
where $D_{J,\nu;b}$ and $D_{J|_Y,\nu_Y;b}$ are the linearizations of 
$\bpar_J\!+\!\nu$ and $\bpar_{J|_Y}\!+\!\nu_Y$  at $b$.
The second equality above is immediate from the transversality of~$s$. 
By~\e_ref{hyperpgen_e6a}, 
\begin{equation}\label{hyperpgen_e6b}\begin{split}
\dim_{\R}\Im  ds_{1,k}^A|_{(\Si;u)}
& =\dim\ker D_{J,\nu;b}-\dim\ker D_{J|_Y,\nu_Y;b}\\
& =\dim_{1,k}(X,A)-\dim_{1,k}(Y,A) 
=2\lr{c_1(\L),A} =\dim_{\R}\V_{1,k}^A|_{[b]}.
\end{split}\end{equation}
The second equality above follows from our assumption that 
the operators $D_{J,\nu;b}$ and $D_{J|_Y,\nu_Y;b}$ are surjective;
The last equality is a consequence of the $(\om,\L)$-positivity assumption.
By~\e_ref{hyperpgen_e6b},  $s_{1,k}^A$ is transverse to the zero set in $\V_{1,k}^A$
along $\M_{1,k}^0(X,A;J,\nu)$.
Since $f$ is transverse to the zero set in $\W$ along $\M_{1,k}^0(X,A;J,\nu_Y)$,
it then follows from~\e_ref{hyperpgen_e5b} that 
\begin{equation}\label{hyperpgen_e9}
\blr{\psi,\PD_{\ov\M_{1,k}^0(X,A;J)}e(\V_{1,k}^A)}
\equiv~^{\pm}\big|f^{-1}(0)\!\cap\!\big\{s_{1,k}^A\big\}^{-1}(0)\big|
=~^{\pm}\big|f^{-1}(0)\!\cap\!\ov\M_{1,k}^0(Y,A;J,\nu_Y)\big|.
\end{equation}
Theorem~\ref{hyperpgen_thm} follows from \e_ref{hyperpgen_e3b} and \e_ref{hyperpgen_e9}.

\section{Ingredients in Proof of Theorem~\ref{main_thm}}
\label{main_sec}

\subsection{Notation: Genus-Zero Maps}
\label{notation0_subs}

\noindent
In this subsection we describe our detailed notation for bubble maps 
from genus-zero Riemann surfaces and for related objects.
In general, moduli spaces of stable maps can stratified by the dual graph.
However, in the present situation, it is more convenient to make use
of {\it linearly ordered sets}:

\begin{dfn}
\label{index_set_dfn1}
(1) A finite nonempty partially ordered set $I$ is a 
{\tt linearly ordered set} if 
for all \hbox{$i_1,i_2,h\!\in\!I$} such that $i_1,i_2\!<\!h$, 
either $i_1\!\le\!i_2$ or $i_2\!\le\!i_1$.\\
(2) A linearly ordered set $I$ is a {\tt rooted tree} if
$I$ has a unique minimal element, 
i.e.~there exists \hbox{$\hat{0}\!\in\!I$} such that $\hat{0}\!\le\!i$ 
for {all $i\!\in\!I$}.
\end{dfn}

\noindent
If $I$ is a linearly ordered set, let $\hat{I}$ be 
the subset of the non-minimal elements of~$I$.
For every $h\!\in\!\hat{I}$,  denote by $\io_h\!\in\!I$
the largest element of $I$ which is smaller than~$h$, i.e.
$\io_h\!=\!\max\big\{i\!\in\!I:i\!<\!h\big\}$.\\

\noindent
We identify $\C$ with $S^2\!-\!\{\i\}$ via 
the stereographic projection mapping the origin in $\C$ 
to the north pole, or the point $(0,0,1)$, in~$S^2$.
Let $M$ be a finite set.
A {\tt genus-zero $X$-valued bubble map with $M$-marked points} is a~tuple
$$b=\big(M,I;x,(j,y),u\big),$$
where $I$ is a rooted tree, and
\begin{equation}\label{stablemap_e1}
x\!:\hat{I}\!\lra\!\C\!=\!S^2\!-\!\{\i\},\quad  j\!:M\!\lra\!I,\quad
y\!:M\!\lra\!\C,        \hbox{~~~and~~~} 
u\!:I\!\lra\!C^{\i}(S^2;X)
\end{equation}
are maps such that $u_h(\i)\!=\!u_{\io_h}(x_h)$ for all $h\!\in\!\hat{I}$.
We associate such a tuple with Riemann surface
\begin{equation}\label{stablemap_e2}
\Si_b=
\Big(\bigsqcup_{i\in I}\Si_{b,i}\Big)\Big/\!\sim,
\hbox{~~where}\qquad \Si_{b,i}=\{i\}\!\times\!S^2
\quad\hbox{and}\quad
(h,\i)\sim (\io_h,x_h)
~~\forall h\!\in\!\hat{I},
\end{equation}
with marked points 
$$y_l(b)\!\equiv\!(j_l,y_l)\in\Si_{b,j_l} \qquad\hbox{and}\qquad  
y_0(b)\!\equiv\!(\hat{0},\i)\in\Si_{b,\hat{0}},$$
and continuous map $u_b\!:\Si_b\!\lra\!X$,
given by $u_b|_{\Si_{b,i}}\!=\!u_i$ for \hbox{all $i\!\in\!I$}.
The general structure of bubble maps is described
by tuples ${\cal T}\!=\!(M,I;j,\under{A})$, where 
$$A_i=u_{i*}[S^2]\in H_2(X;\Z) \qquad\forall i\!\in\!I.$$ 
We call such tuples {\tt bubble types}.
We denote by $\ti\X_{\T}(X)$ be the space of all bubble maps of type~$\T$.
For $l\!\in\!\{0\}\!\sqcup\!M$, let 
$$\ev_l\!:\ti\X_{\T}(X)\lra X$$ 
be the evaluation map corresponding to the marked point $y_l$.\\

\noindent
With notation as above, suppose
$$b\!\equiv\!\big(M,I;x,(j,y),u\big) \in \ti\X_{\T}(X).$$
In particular, $I$ is a linearly ordered set with minimal element $\hat{0}$
and the special marked point is the~point 
$$y_0(b)=(\hat{0},\i)\in\Si_{b,\hat{0}}.$$
Let $\chi^0(b)$ be the set of components $\Si_{b,i}$ of $\Si_b$ such that 
for every component $\Si_{b,h}$ that lies between $\Si_i$ and  $\Si_{b,\hat{0}}$,
including $\Si_{b,i}$ and $\Si_{b,\hat{0}}$, the degree of $u|_{\Si_{b,h}}$ is zero.
The set $\chi^0(b)$ is empty if and only if the degree of the restriction of $u_b$ 
to the component containing the special marked point is not zero.
Let
$$\Si_b^0=\big\{(\hat{0},\i)\big\} \cup 
\bigcup_{i\in\chi^0(b)}\!\!\!\Si_{b,i}.$$
We denote by 
$$\chi(b)\!=\!\chi(\T) \subset I$$ 
the set of components $\Si_{b,i}$ of $\Si_b$ such that $\Si_{b,i}$ has a point 
in common with $\Si_b^0$ and the degree of $u_b|_{\Si_{b,i}}$ is not zero,
i.e.~$\Si_{b,i}$ is not an element of $\chi^0(b)$.\\

\noindent
If $b\!=\!(\Si_b;u_b)$ is a bubble map with a special marked point as above
and $i\!\in\!\chi(b)$, we~put
$$\cD_ib= du_{b,i}\big|_{\i}e_{\i} \in T_{u_{b,i}(\i)}X,$$
where $u_{b,i}\!=\!u_b|_{\Si_{b,i}}$.
Similarly, if $(\L,\na)$ is a complex line bundle with connection over $X$ 
and $\xi\!=\!(\xi_h)_{h\in I}$ is an element of $\Ga(b;\L)$, we~put
$$\D_{b,i}\xi=\na_{e_{\i}}^{u_b} \xi_i\big|_{\i} \in \L_{u_{b,i}(\i)}.$$
Note that if $\xi\!\in\!\ker\bpar_{\na,b}$,
\begin{equation}\label{secder_e}
\na_{c\cdot e_{\i}}^{u_b} \xi_i\big|_{\i}=c\cdot\D_{b,i}\xi
\qquad\forall\,c\!\in\!\C.
\end{equation}
If in addition $\L$ is $(\om,u_{b*}[\Si_b])$-positive, then
the linear operator $\bpar_{\na,b}$ and the linear map
$$\bpar_{\na,b}^{y_0(b)}\!: \ker\bpar_{\na,b}\lra\L_{\ev_0(b)}, \qquad
\xi\lra\xi\big(y_0(b)\big),$$
are surjective. This can be seen by an argument similar to Subsection~6.2 in~\cite{g2n2and3}.

\subsection{Notation: Genus-One Maps}
\label{notation1_subs}

\noindent
We next set up notation for maps from genus-one Riemann surfaces.
In this case, in contrast to the genus-zero case,
we also need to specify the structure of the principal component.
We describe it by {\it enhanced linearly ordered sets}:

\begin{dfn}
\label{index_set_dfn2}
An {\tt enhanced linearly ordered set} is a pair $(I,\aleph)$,
where $I$ is a linearly ordered set, $\aleph$ is a subset of $I_0\!\times\!I_0$,
and $I_0$ is the subset of minimal elements of~$I$,
such that if $|I_0|\!>\!1$, 
$$\aleph=\big\{(i_1,i_2),(i_2,i_3),\ldots,(i_{n-1},i_n),(i_n,i_1)\big\}$$
for some bijection $i\!:\{1,\ldots,n\}\!\lra\!I_0$.
\end{dfn}

\begin{figure}
\begin{pspicture}(-1.1,-2)(10,1)
\psset{unit=.4cm}
\pscircle*(6,-3){.2}
\pscircle*(4,-1){.2}\psline[linewidth=.06]{->}(4.14,-1.14)(5.86,-2.86)
\pscircle*(8,-1){.2}\psline[linewidth=.06]{->}(7.86,-1.14)(6.14,-2.86)
\pscircle*(2,1){.2}\psline[linewidth=.06]{->}(2.14,.86)(3.86,-.86)
\pscircle*(6,1){.2}\psline[linewidth=.06]{->}(5.86,.86)(4.14,-.86)
\pscircle*(18,-3){.2}\psline[linewidth=.06](17.86,-3.14)(17.5,-3.5)
\psarc(18,-4){.71}{135}{45}\psline[linewidth=.06]{->}(18.5,-3.5)(18.14,-3.14)
\pscircle*(16,-1){.2}\psline[linewidth=.06]{->}(16.14,-1.14)(17.86,-2.86)
\pscircle*(20,-1){.2}\psline[linewidth=.06]{->}(19.86,-1.14)(18.14,-2.86)
\pscircle*(14,1){.2}\psline[linewidth=.06]{->}(14.14,.86)(15.86,-.86)
\pscircle*(18,1){.2}\psline[linewidth=.06]{->}(17.86,.86)(16.14,-.86)
\pscircle*(30,-2){.2}\pscircle*(30,-4){.2}\pscircle*(29,-3){.2}\pscircle*(31,-3){.2}
\psline[linewidth=.06]{->}(29.86,-2.14)(29.14,-2.86)
\psline[linewidth=.06]{->}(29.14,-3.14)(29.86,-3.86)
\psline[linewidth=.06]{->}(30.14,-3.86)(30.86,-3.14)
\psline[linewidth=.06]{->}(30.86,-2.86)(30.14,-2.14)
\pscircle*(27,-1){.2}\psline[linewidth=.06]{->}(27.14,-1.14)(28.86,-2.86)
\pscircle*(33,-1){.2}\psline[linewidth=.06]{->}(32.86,-1.14)(31.14,-2.86)
\pscircle*(25,1){.2}\psline[linewidth=.06]{->}(25.14,.86)(26.86,-.86)
\pscircle*(29,1){.2}\psline[linewidth=.06]{->}(28.86,.86)(27.14,-.86)
\end{pspicture}
\caption{Some Enhanced Linearly Ordered Sets}
\label{index_set_fig}
\end{figure}
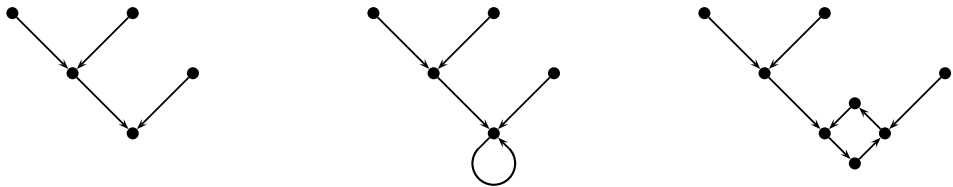

\noindent
An enhanced linearly ordered set can be represented by an oriented connected graph.
In Figure~\ref{index_set_fig}, the dots denote the elements of~$I$.
The arrows outside the loop, if there are any, 
specify the partial ordering of the linearly ordered set~$I$.
In fact, every directed edge outside of the loop
connects a non-minimal element $h$ of $I$ with~$\io_h$.
Inside of the loop, there is a directed edge from $i_1$ to $i_2$
if and only if $(i_1,i_2)\!\in\!\aleph$.\\

\noindent
The subset $\aleph$ of $I_0\!\times\!I_0$ will be used to describe
the structure of the principal curve of the domain of stable maps in 
a stratum of~$\X_{1,M}(X,A)$.
If $\aleph\!=\!\eset$, and thus $|I_0|\!=\!1$,
the corresponding principal curve $\Si_P$ 
is a smooth torus, with some complex structure.
If $\aleph\!\neq\!\eset$, the principal components form a circle of spheres:
$$\Si_P=\Big(\bigsqcup_{i\in I_0}\{i\}\!\times\!S^2\Big)\Big/\sim,
\qquad\hbox{where}\qquad
(i_1,\i)\sim(i_2,0)~~\hbox{if}~~(i_1,i_2)\!\in\!\aleph.$$
A {\tt genus-one $X$-valued bubble map with $M$-marked points} is a tuple
$$b=\big(M,I,\aleph;S,x,(j,y),u\big),$$
where $S$ is a smooth Riemann surface of genus one if $\aleph\!=\!\eset$
and the circle of spheres $\Si_P$ otherwise.
The objects $x$, $j$, $y$, $u$, and $(\Si_b,u_b)$ are as in 
\e_ref{stablemap_e1} and \e_ref{stablemap_e2}, except 
the sphere $\Si_{b,\hat{0}}$ is replaced by the genus-one curve $\Si_{b;P}\!\equiv\!S$.
Furthermore, if $\aleph\!=\!\eset$, and thus $I_0\!=\!\{\hat{0}\}$ is a single-element set,
$u_{\hat{0}}\!\in\!C^{\i}(S;X)$ and $y_l\!\in\!S$ if $j_l\!=\!\hat{0}$.
In the genus-one case, the general structure of bubble maps is encoded by
the tuples of the form ${\cal T}\!=\!(M,I,\aleph;j,\under{A})$.
Similarly to the genus-zero case,
we denote by $\ti\X_{\T}(X)$ be the space of all bubble maps of type~$\T$.
Let 
$$\X_{\T}(X)=\big\{[b]\!\in\!\X_{1,M}(X,A)\!: b\!\in\!\ti\X_{\T}(X)\big\}.$$\\

\noindent
If $\nu$ is an element of $\G_{1,M}^{\es}(X,A)$, we put
$$\U_{\T,\nu}(X;J)=\big\{[b]\!\in\!\X_{\T}(X)\!:
\{\bar{\partial}_J\!+\!\nu\}(b)=0\big\}.$$
If $\T\!=\!(M,I,\aleph;j,\under{A})$ is a bubble type such that $A_i\!=\!0$
for all minimal elements $i$ of~$I$ and $[\Si;u]$ is an element of $\U_{\T,\nu}(X;J)$,
the map $u|_{\Si_P}$ is constant.
Let
$$\ev_P\!:\U_{\T,\nu}(X;J) \lra X$$
be the map sending each element $[\Si;u]$ of $\U_{\T,\nu}(X;J)$
to the image of the principal component $\Si_P$ of $\Si$, 
i.e.~the point $u(\Si_P)$ in~$X$. We note that the map
$$\ti\X_{\T}(X)\lra 2^I, \qquad b\lra \chi(b),$$
is constant. We denote its value by $\chi(\T)$.\\

\noindent
Suppose $b\!=\!(\Si_b;u_b)$ is an element of $\X_{1,k}(X,A)$ as above
and $(\L,\na)$ is a complex vector bundle with connection over $X$.
If $\xi\!=\!(\xi_h)_{h\in I}$ is an element of $\Ga(b;\L)$ and $i\!\in\!I\!-\!I_0$,
similarly to the genus-zero case, we~put
$$\D_{b,i}\xi=\na_{e_{\i}}^{u_b} \xi_i\big|_{\i} \in \L_{u_{b,i}(\i)},$$
where $u_{b,i}\!=\!u_b|_{\Si_{b,i}}$.\\

\noindent
Finally, all vector orbi-bundles we encounter will be assumed to be normed.
Some will come with natural norms; for others, we implicitly choose a norm
once and for~all.
If \hbox{$\pi_{\F}\!:\F\!\lra\!\X$} is a normed vector bundle
and $\de\!:\X\!\lra\!\R$ is any function, possibly constant,
let
$$\F_{\de}=\big\{\ups\!\in\!\F\!: |\ups|\!<\!\de(\pi_{\F}(\ups))\big\}.$$
If $\Om$ is any subset of $\F$, we take  $\Om_{\de}\!=\!\Om\cap\F_{\de}$.

\subsection{Topology}
\label{top_subs}

\noindent
In this subsection we prove a general topological result, Proposition~\ref{conetop_prp}.
For the sake of clarity, we state and prove it only in the manifold/section category,
but Proposition~\ref{conetop_prp} and its proof carry over easily to 
the orbifold/multisection category.\\

\noindent
Similarly to Subsection~\ref{eulercone_subs}, by a cone $\pi\!:\W\!\lra\!\ov\M$ we mean
a continuous map between two topological spaces such that 
$\W_x\!\equiv\!\pi^{-1}(x)$ is vector space for each $x\!\in\!\ov\M$
and the vector space operations induce continuous functions on $\W\!\times\!\W$
and $\C\!\times\!\W$.

\begin{dfn}
\label{coneprop_dfn}
(1) A cone $\pi\!:\W\!\lra\!\ov\M$ is {\tt regular} if for every $b\!\in\!\ov\M$,
there exist a neighborhood $U_b$ of $b$ in $\ov\M$, $n_b\!\in\!\Z^+$, and a bundle map
$$\vph_b\!:  \W|_{U_b}\lra U_b\!\times\!\C^{n_b}$$
over $U_b$ such that $\vph_b$ is a homeomorphism onto its image and the restriction
of $\vph_b$ to each fiber is linear.\\
(2) A cone $\pi\!:\W\!\lra\!\ov\M$ is {\tt obstruction-free} if for every $b\!\in\!\ov\M$,
$\xi\!\in\!\W_b$, and a sequence $b_r\!\in\!\ov\M$ converging to~$b$,
there exists a sequence $\xi_r\!\in\!\W_{b_r}$ converging to $\xi$ in~$\W$.
\end{dfn}

\noindent
If $\W\!\lra\!\ov\M$ is a cone, for each $r\!\in\!\Z$, let
$$\ov\M_r(\W)=\big\{b\!\in\!\ov\M\!: \rk\W_b\!=\!r\big\}.$$
Note that if $\W$ is obstruction-free, then the set 
$\bigcup_{r\le q}\!\ov\M_r(\W)$ is closed in~$\ov\M$.

\begin{lmm}
\label{conetop_lmm}
Suppose $\ov\M$ is a compact Hausdorff space that has a countable basis at each point, 
$A$ is a closed subset of $\ov\M$, $\W\!\lra\!\ov\M$ is a regular obstruction-free cone, 
and $s$ is a section of $\W$ over~$A$.
If
$$r_A=\min\big\{\rk\W_b\!: b\!\in\!\ov\M\!-\!A\big\},$$
$s$ extends to a continuous section $\ti{s}$ of $\W$ over $A\!\cup\!\ov\M_{r_A}(\W)$.
\end{lmm}

\noindent
{\it Remark 1:} It is enough to assume that $\ov\M$ is a paracompact (Hausdorff)
space that has a countable basis at each point.\\

\noindent
{\it Remark 2:} An immediate corollary of this lemma is that $s$ extends to 
a continuous section of $\W$ over~$\ov\M$.\\

\noindent
Let $\{U_b\}_{b\in\A}$ be a finite open cover of the compact set $A\!\cup\!\ov\M_{r_A}(\W)$
by open subspaces of $\ov\M$ as in~(1) of Definition~\ref{coneprop_dfn}.
Since $\ov\M$ is normal, we can choose an open cover $\{U_b'\}_{b\in\A}$ 
of $A\!\cup\!\ov\M_{r_A}(\W)$ such that $\bar{U}_b'\!\subset\!U_b$ for all $b\!\in\!\A$.
Since $\bar{U}_b'$ is normal for each $b\!\in\!\A$,
by the Tietze Extension Theorem the continuous section $f\!\equiv\!\vph_b\!\circ\!s$ 
of $U_b\!\times\!\C^{n_b}$ over $A\!\cap\!\bar{U}_b'$ extends to a continuous section
$\ti{f}$ over~$\bar{U}_b'$.
Let 
$$\pi_b^-\!: U_b\!\times\!\C^{n_b} \lra \Im\vph_b$$
be the orthogonal projection map.
We will show in the next paragraph that the section $\pi_b^-\!\circ\!\ti{f}$ 
is continuous over $(A\!\cup\!\ov\M_{r_A}(\W))\!\cap\!\bar{U}_b'$.
Since $A\!\cup\!\ov\M_{r_A}(\W)$ is normal, we can choose a partition of unity 
$\{\eta_b\}_{b\in\A}$ subordinate to~$\{U_b'\}_{b\in\A}$.
The section 
$$\ti{s}=\sum_{b\in\A}\eta_b\cdot \big(\vph_b^{-1}\!\circ\!\pi_b^-\!\circ\!\ti{f}\big)$$
is continuous over $A\!\cup\!\ov\M_{r_A}(\W)$.
Since $f(x)\!\in\!\Im\vph_b$ for all $x\!\in\!A\!\cap\!\bar{U}_b'$, 
$$\pi_b^-\ti{f}(x)=\pi_b^-f(x)=f(x) \qquad\forall\,x\!\in\!A\!\cap\!\bar{U}_b'.$$
Thus, $\ti{s}|_A\!=\!s$ as required.\\

\noindent
It remains to show that the section $\pi_b^-\!\circ\!\ti{f}$
is continuous over $(A\!\cup\!\ov\M_{r_A}(\W))\!\cap\!\bar{U}_b'$.
Since 
$$\pi_b^-\ti{f}|_{A\cap \bar{U}_b'}=f|_{A\cap \bar{U}_b'},$$
$\pi_b^-\ti{f}$ is continuous along the closed subset  $A\!\cap\!\bar{U}_b'$
of $\bar{U}_b'$.
Thus, we need to show that if 
$$x_r \in \ov\M_{r_A}(\W)\!\cap\!\bar{U}_b'$$ 
is a sequence converging to $x\!\in\!\bar{U}_b'$, 
then $\pi_b^-\ti{f}(x_r)$ converges to $\pi_b^-\ti{f}(x)$.
Suppose first
$$x\in\ov\M_{r_A}(\W)\!\cap\!\bar{U}_b'.$$ 
We will show that $\pi_b^-|_{\ov\M_{r_A}(\W)\cap\bar{U}_b'}$ is continuous at $x$.
Let $\{\xi_i\}_{i\in[r_A]}$ be an orthonormal basis for $\Im\vph_b|_x$.
By~(2) of Definition~\ref{coneprop_dfn}, for each $i\!\in\![r_A]$ there exists 
a sequence $\xi_{i;r}\!\in\!\Im\vph_b|_{x_r}$ converging to~$\xi_i$.
Since $b,b_r\!\in\!\ov\M_{r_A}(\W)$, 
$\{\xi_{i;r}\}_{i\in[r_A]}$ is basis for $\Im\vph_b|_{x_r}$ for all $r$ sufficiently large.
Since $\xi_{i;r}\!\lra\!\xi_i$ for all $i\!\in\![r_A]$ and 
$\{\xi_i\}_{i\in[r_A]}$ is an orthonormal basis,
$$\lim_{r\lra\i}\lr{\xi_{i;r},\xi_{j;r}}=\begin{cases}
1,& \hbox{if}~i\!=\!j;\\ 0,& \hbox{if}~i\!\neq\!j.\end{cases}$$
Thus, applying the Gramm-Schmidt normalization procedure, we can find an orthonormal
basis $\{\ti\xi_{i;r}\}_{i\in[r_A]}$ for $\Im\vph_b|_{x_r}$ such that
$\ti\xi_{i;r}\!\lra\!\xi_i$ for all $i\!\in\![r_A]$.
It follows that  $\pi_b^-|_{x_r}\!\lra\!\pi_b^-|_x$ as needed.
On the other hand, suppose that $x\!\in\!A\!\cap\!\bar{U}_b'$.
We will view $f$, $\ti{f}$, and $\pi_b^-\ti{f}\!=\!\pi_x\ti{f}$
as $\C^{n_b}$-valued functions.
By~(2) of Definition~\ref{coneprop_dfn}, 
there exists a sequence $\xi_r\!\in\!\Im\vph_b|_{x_r}$ converging to~$f(x)$.
Since $\pi_{x_r}\xi_r\!=\!\xi_r$,
\begin{equation*}\begin{split}
\big|\pi_x\ti{f}(x)-\pi_{x_r}\ti{f}(x_r)\big|
& =\big|\ti{f}(x)-\pi_{x_r}\ti{f}(x_r)\big|\\
&\le \big|\ti{f}(x)-\xi_r\big|+
\big|\pi_{x_r}\xi_r-\pi_{x_r}\ti{f}(x)\big|
+\big|\pi_{x_r}\ti{f}(x)-\pi_{x_r}\ti{f}(x_r)\big|\\
&\le 2\big|f(x)-\xi_r\big|+\big|\ti{f}(x)-\ti{f}(x_r)\big|.
\end{split}\end{equation*}
The last two terms above approach $0$ by the assumption on $\xi_r$ and 
the continuity of~$\ti{f}_r$.\\

\noindent
{\it Remark:} The projection $\pi_b^-$ is {\it not} continuous over $U_b$
unless the rank of $\W$ is constant over~$U_b$.
Similarly, the section $\pi_b^-\!\circ\!\ti{f}$ may not be continuous over $\bar{U}_b'$
unless the rank of $\W$ is constant over~$\bar{U}_b'\!-\!A$.

\begin{dfn}
\label{strat_dfn}
If $\ov\M$ is a topological space and $(\A,\prec)$ is a finite partially ordered set,
a collection $\{\U_{\al}\}_{\al\in\A}$ of subspaces of $\ov\M$ is
a {\tt stratification of $\ov\M$} if $\U_{\al}$ is 
a smooth manifold for all $\al\!\in\!\A$,
$$\partial\bar\U_{\al}\!\equiv\!\bar\U_{\al}\!-\!\U_{\al}\subset
\bigcup_{\be\prec\al}\!\U_{\be}  \qquad\forall\,\al\!\in\!\A,
\qquad\hbox{and}\qquad \ov\M=\bigsqcup_{\al\in\A}\U_{\al}.$$
\end{dfn}

\begin{prp}
\label{conetop_prp}
Suppose $\ov\M$ is a compact Hausdorff space that has a countable basis at each point,
$\{\U_{\al}\}_{\al\in\A}$ is a stratification of $\ov\M$,
and   $\W\!\lra\!\ov\M$ is a regular obstruction-free cone.
If $\W|_{\U_\al}\!\lra\!\U_{\al}$ is a smooth vector bundle for all $\al\!\in\!\A$,
the cone $\W$ admits a continuous section $s$ over $\ov\M$ such that 
$s|_{\U_\al}$ is smooth and transverse to the zero set in $\W|_{\U_\al}$ 
for all $\al\!\in\!\A$.
\end{prp}

\noindent
Choose an ordering $<$ on the partially ordered set $(\A,\prec)$ 
such that for all $\al,\be\!\in\!\A$, 
\begin{alignat}{1}
\rk\W|_{\U_\be}\!<\!\rk\W|_{\U_\al} \qquad&\Lra\qquad \be\!<\!\al;
\label{conetopprp_e1} \\
\rk\W|_{\U_\be}\!=\!\rk\W|_{\U_\al},~~ \be\!\prec\!\al  \qquad&\Lra\qquad \be\!<\!\al.
\notag
\end{alignat}
Since $\W$ is obstruction-free, $\bigcup_{r\le q}\!\ov\M_r(\W)$ is closed in~$\ov\M$, and
\begin{equation}\label{conetopprp_e2}
\partial\bar\U_{\al}\subset\bigcup_{\be<\al}\!\U_{\be} 
\qquad\forall\, \al\!\in\!\A,
\end{equation} 
by the closure condition of Definition~\ref{strat_dfn}.
Suppose $\al\!\in\!\A$ and we have defined a continuous section $s$ of $\W$ over
the closed set
$$A\equiv\bigcup_{\be<\al}\U_{\be}$$
such that $s|_{\U_\be}$ is smooth and transverse to the zero set in $\W|_{\U_\be}$
for all $\be\!<\!\al$.
By~\e_ref{conetopprp_e1} and~\e_ref{conetopprp_e2},
$$r_A\!\equiv\!\min\big\{\rk\W_b\!: b\!\in\!\ov\M\!-\!A\big\}=\rk\W|_{\U_\al}.$$
Thus, by Lemma~\ref{conetop_lmm}, $s$ extends to a continuous section~$\ti{s}$
over
$$A\!\cup\!\U_{\al}\subset A\!\cup\!\ov\M_{r_A}(\W).$$
Perturbing $\ti{s}$ over $\U_{\al}$, without changing it over $A$,
we obtain a continuous section~$s$ over $A\!\cup\!\U_{\al}$ 
such that $s|_{\U_\be}$ is smooth and transverse to the zero set in $\W|_{\U_\be}$
for all $\be\!\le\!\al$.
This construction implies Proposition~\ref{conetop_prp}.\\

\noindent
{\it Remark 1:} In the orbifold/multisection category as needed for the purposes 
of this paper, a stratum $\U_{\al}$ locally is a union of finitely many smooth
suborbifolds of a smooth orbifold~$\X_{\al}$. 
We will still call such unions smooth orbifolds.
The bundle $\W|_{\U_{\al}}$ is the restriction of a smooth orbibundle over~$\X_{\al}$.\\

\noindent
{\it Remark 2:} The cone $\V_{1,k}^A$ is not obstruction-free, but is regular.
In Subsection~\ref{g1conestr_subs}, we describe a subcone $\W_{1,k}^A\!\subset\!\V_{1,k}^A$
which is obstruction-free and sufficiently large for the purposes of Theorem~\ref{main_thm}.

\subsection{The Structure of the Moduli Space $\ov\M_{1,k}^0(X,A;J,\nu)$}
\label{g1str_subs}

\noindent
In this subsection, we describe the strata of the moduli space $\ov\M_{1,k}^0(X,A;J,\nu)$
for a small generic element $\nu$ of $\G^{\es}_{1,k}(X,A;J)$.
If $k\!\in\!\Z$, we denote by $[k]$ the set of positive integers
that do not exceed~$k$.
Let $2n\!=\!\dim_{\R}\!X$.

\begin{lmm}
\label{g1mainstr_lmm}
Suppose $(X,\om,J)$, $A$, $k$, and  $\nu$ are as in Theorem~\ref{main_thm}.
If
$$\T=([k],I,\aleph;j,\under{A})$$
 is a bubble type such that
$\sum_{i\in I}\!A_i\!=\!A$ and $A_i\!\neq\!0$ for some minimal element $i$ of $I$,
then $\U_{\T,\nu}(X;J)$ is a smooth orbifold and
$$\dim\U_{\T,\nu}(X;J)=\dim_{1,k}(X,A)-2\big(|\aleph|\!+\!|\hat{I}|\big).$$\\
\end{lmm}

\noindent
The statement that $\U_{\T,\nu}(X;J)$ is smooth should be interpreted 
as in Remark~1 at the end of the previous subsection.
The branches of $\U_{\T,\nu}(X;J)$ correspond to the branches of~$\nu$.
For a generic~$\nu$, the linearization $D_{J,\nu;b}$ of the bundle section
$\bpar_J\!+\!\nu$
at $[b]$ is surjective for every element $b$ in $\ov\M_{1,k}^0(X,A;J,\nu)$
such that $u_b|_{\Si_{b;P}}$ is not constant.
Thus, Lemma~\ref{g1mainstr_lmm} is obtained by a standard Contraction Principle argument,
such as in Chapter~3 of~\cite{McSa}.

\begin{lmm}
\label{g1bdstr_lmm}
Suppose $(X,\om,J)$, $A$, $k$, and $\nu$ are as in Theorem~\ref{main_thm}.
If
$$\T=([k],I,\aleph;j,\under{A})$$
is a bubble type such that
$\sum_{i\in I}\!A_i\!=\!A$ and $A_i\!=\!0$ for all minimal elements $i$ of~$I$, then
for each $m\!\in\!\Z^+$,
$$\U_{\T,\nu;1}^m(X;J) \equiv \big\{
[b]\!\in\!\U_{\T,\nu}(X;J)\!:\dim_{\C}\hbox{Span}_{(\C,J)}
\{\cD_ib\!:i\!\in\!\chi(\T)\}=|\chi(\T)|\!-\!m\big\}$$
is a smooth orbifold and 
$$\dim\U_{\T,\nu;1}^m(X;J) =  \dim_{1,k}(X,A)
-2\big(|\aleph|\!+\!|\hat{I}|-n+(m\!+\!n\!-\!|\chi(\T)|)m\big).$$\\
\end{lmm}

\noindent
If $\T$ is as in Lemma~\ref{g1bdstr_lmm} and $b\!\equiv\!(\Si_b;u_b)\!\in\!\U_{\T,\nu}(X;J)$,
$u_b|_{\Si_{b;P}}$ is constant. Let
\begin{alignat*}{1}
\Ga_B(b) &=\big\{\ze\!\in\!\Ga(\Si_b;u_b^*TX)\!: \ze|_{\Si_{b;P}}\!=\!0\big\};\\
\Ga_B^{0,1}(b;J) &=\big\{\eta\!\in\!\Ga(\Si_b;\La^{0,1}_{J,j}T^*\Si_b\!\otimes\!u_b^*TX)\!:
\eta|_{\Si_{b;P}}\!=\!0\big\}.
\end{alignat*}
If $\nu$ is genetic, the operator 
$$D_{J,\nu;b}^B\!: \Ga_B(b) \lra \Ga_B^{0,1}(b;J)$$
induced by $D_{J,\nu;b}$ is surjective. 
Thus, the space $\U_{\T,\nu}(X;J)$ is a smooth orbifold of dimension
$$\dim\U_{\T,\nu;1}^m(X;J) =  \dim_{1,k}(X,A)
-2\big(|\aleph|\!+\!|\hat{I}|\big)+2n.$$
We note that 
\begin{equation}\label{mbound_e}
\dim\U_{\T,\nu;1}^m(X;J) \neq \eset \qquad\Lra\qquad
 \max(1,|\chi(\T)|\!-\!n)\le m\le |\chi(\T)|.
\end{equation}
As at the end of Subsection~\ref{g1cone-g1str_subs} in~\cite{g1cone}, 
we can construct a vector bundle $F$ over $\U_{\T,\nu}(X;J)$ of rank $|\chi(\T)|$,
a vector bundle $V$ over 
$$\pi_m\!:\Gr_mF\lra \U_{\T,\nu;1}^m(X;J)$$ 
of rank $mn$, and a transverse section $\cD_m$ of~$V$ such that 
$$\pi_m\!: \cD_m^{-1}(0) \lra \bigcup_{m'\ge m}\!\!\U_{\T,\nu;1}^{m'}(X;J)$$
is surjective, and the restriction of $\pi_m|_{\cD_m^{-1}(0)}$ 
to the preimage of $\U_{\T,\nu;1}^m(X;J)$ is an embedding.
This observation implies the dimension claim of Lemma~\ref{g1bdstr_lmm}.\\

\noindent
The spaces $\U_{\T,\nu}(X;J)$ and $\U_{\T,\nu;1}^m(X;J)$ of Lemmas~\ref{g1mainstr_lmm}
and~\ref{g1bdstr_lmm} are disjoint.
By Definition~\ref{degen_dfn}, their union is $\ov\M_{1,k}^0(X,A;J,\nu)$.
Let $\A^*$ be the set of equivalence classes of bubble types~$\T$ as in Lemma~\ref{g1mainstr_lmm}
and $\A^0$ the set of equivalence classes of pairs~$(\T,m)$ consisting of 
a bubble type~$\T$ as in Lemma~\ref{g1bdstr_lmm} and an integer~$m$ as in~\e_ref{mbound_e}.
We define a partial ordering on the set $\A\!\equiv\!\A^*\!\sqcup\!\A^0$ as follows.
Suppose
$$\T=([k],I,\aleph;j,\under{A}) \qquad\hbox{and}\qquad
\T'=([k],I',\aleph';j',\under{A}')$$
are two bubble types as in Lemma~\ref{g1mainstr_lmm} and/or in Lemma~\ref{g1bdstr_lmm}.
We write 
$$\T'\!\prec\!\T   \qquad\Llra\qquad
|\hat{I}'|\!\ge\!|\hat{I}|,~~|\ale'|\!\ge\!|\ale|,~~
|\hat{I}'|\!+\!|\ale'|\!>\!|\hat{I}|\!+\!|\ale|.$$
If $\T\!\in\!\A^*$ and $(\T',m')\!\in\!\A^0$, we define
$$\T\!\prec\!(\T',m') ~~~\Llra~~~ \T\!\prec\!\T'; \qquad
(\T',m')\!\prec\!\T ~~~\Llra~~~ \T'\!\prec\!\T.$$
Finally, if $(\T,m),(\T',m')\!\in\!\A^0$, we define
$$(\T',m')\!\prec\!(\T,m) \qquad\Llra\qquad
\T'\!\prec\!\T \quad\hbox{OR}\quad \T'\!=\!\T,~m'\!>\!m.$$
By definition of the stable-map topology, 
$$\partial\bar\U_{\T,\nu}(X;J)\cap \U_{\T',\nu}(X;J)\neq\eset
\qquad\Lra\qquad \T'\!\prec\!\T.$$
Thus, the closure requirement of Definition~\ref{strat_dfn} follows
from the continuity of the maps~$\cD_i$ on $\ti\X_{\T}(X)$ with $\T$ 
as in Lemma~\ref{g1bdstr_lmm}.

\subsection{The Structure of the Cone $\V_{1,k}^A$}
\label{g1conestr_subs}

\noindent
In this subsection, we describe an obstruction-free subcone $\W_{1,k}^A$ of the cone 
$$\V_{1,k}^A\lra\ov\M_{1,k}^0(X,A;J,\nu).$$ 
The cone $\V_{1,k}^A$ over $\X_{1,k}^A$ can be shown to be regular by standard arguments;
see Remark~2 at the end of Subsection~\ref{gluing1_subs}.
Thus, $\W_{1,k}^A$ is regular as~well.\\

\noindent
If $\T$ is a bubble type as in Lemma~\ref{g1bdstr_lmm} and 
$[b]\!\in\![\Si_b,u_b]\!\in\!\U_{\T,\nu}(X;J)$, let
\begin{alignat*}{1}
F_b^1 &=\big\{(w_i)_{i\in\chi(\T)}\!\in\!\C^{\chi(b)}\!:
\sum_{i\in\chi(\T)}\!\!w_i\!\cdot_J\cD_ib\!=\!0\big\};\\
\ti{\Ga}_-(b;\L)&=\big\{\xi\!\in\!\ker\bpar_{\na,b}\!:
\sum_{i\in\chi(\T)}\!\!w_i\,\D_{b,i}\xi\!=\!0~\forall\,
(w_i)_{i\in\chi(\T)}\!\in\!F_b^1\big\};\\
\W_{1,k}^A|_{[b]}&=\big\{[(b,\xi)]\!\in\!\V_{1,k}^A|_{[b]}\!: \xi\!\in\!\ti{\Ga}_-(b;\L)\big\}.
\end{alignat*}
By~\e_ref{secder_e},
the subspace $\W_{1,k}^A|_{[b]}$ of $\V_{1,k}^A|_{[b]}$ is well-defined.
If $\T$ is a bubble type as in Lemma~\ref{g1mainstr_lmm} and $[b]\in\!\U_{\T,\nu}(X;J)$,
let
$$\W_{1,k}^A|_{[b]}=\V_{1,k}^A|_{[b]}.$$
We take
$$\W_{1,k}^A=\bigcup_{[b]\in\ov\M_{1,k}^0(X,A;J,\nu)}\!\!\!\!\!\!\!\!\!\!\!\!\!\!
\W_{1,k}^A|_{[b]} ~\subset \V_{1,k}^A.$$

\begin{prp}
\label{conestr_prp}
If $(X,\om,J)$, $(\L,\na)$, $A$, $k$, and $\nu$ are as in Theorem~\ref{main_thm}, the cone
$$\W_{1,k}^A \lra \ov\M_{1,k}^0(X,A;J,\nu)$$
is regular and obstruction-free.
If $\T$ is a bubble type as in Lemma~\ref{g1mainstr_lmm}, then 
$\W_{1,k}^A|_{\U_{\T,\nu}(X;J)}$ is a smooth vector orbibundle~and
$$\rk \W_{1,k}^A|_{\U_{\T,\nu}(X;J)}=\lr{c_1(\L),A}.$$
If $\T$ is a bubble type as in Lemma~\ref{g1bdstr_lmm} and $m$ is an integer 
as in~\e_ref{mbound_e},  then 
$\W_{1,k}^A|_{\U_{\T,\nu;1}^m(X;J)}$ is a smooth vector orbibundle~and
$$\rk\W_{1,k}^A|_{\U_{\T,\nu;1}^m(X;J)}
=\lr{c_1(\L),A}-(m\!-\!1) >\frac{1}{2}\dim\U_{\T,\nu;1}^m(X;J)
-\frac{1}{2}\dim_{1,k}(X,A;\L).$$\\
\end{prp}

\noindent
If $\T$ is a bubble type as in Lemma~\ref{g1mainstr_lmm} and $[b]\in\!\U_{\T,\nu}(X;J)$,
the operator~$\bpar_{\na,b}$ is surjective,
by the positivity assumption on the bundle $\L$ and the same argument as in Subsection~6.2
of~\cite{g2n2and3}.
In particular, 
$$\dim\V_{1,k}^A|_{[b]}=\ind\bpar_{\na,b}=\lr{c_1(\L),A}
\qquad\forall\, [b]\in\U_{\T,\nu}(X;J).$$
By standard arguments, the surjectivity of $\bpar_{\na,b}$ 
for every $[b]\!\in\!\U_{\T,\nu}(X;J)$ implies that
$$\W_{1,k}^A|_{\U_{\T,\nu}(X;J)}\equiv\V_{1,k}^A|_{\U_{\T,\nu}(X;J)}$$
is a smooth vector bundle and that the restriction of $\V_{1,k}^A$
to a neighborhood of $\U_{\T,\nu}(X;J)$ in~$\X_{1,k}(X,A)$ is a vector bundle.
Local trivializations can be constructed using the homomorphisms $R_{\ups,\ze}$
as in Subsection~\ref{gluing0_subs}.
In particular, the cone $\W_{1,k}^A$ satisfies the requirements of~(1)
and~(2) of Definition~\ref{coneprop_dfn} for every $[b]\in\!\U_{\T,\nu}(X;J)$.\\

\noindent
If $\T$ is a bubble type as in Lemma~\ref{g1bdstr_lmm} and
$[b]\!\equiv\![\Si_b;u_b]\!\in\!\U_{\T,\nu}(X;J)$, $u_b|_{\Si_{b;P}}$ is constant. 
Let
\begin{alignat*}{1}
\Ga_B(b;\L) &=\big\{\xi\!\in\!\Ga(\Si_b;u_b^*\L)\!: \xi|_{\Si_{b;P}}\!=\!\const\big\};\\
\Ga_B^{0,1}(b;\L) &=\big\{\eta\!\in\!\Ga(\Si_b;\La^{0,1}_{\I,j}T^*\Si_b\!\otimes\!u_b^*\L)\!:
\eta|_{\Si_{b;P}}\!=\!0\big\}.
\end{alignat*}
By the positivity assumption on the bundle $\L$ and the same argument as 
in Subsection~6.2 of~\cite{g2n2and3}, the operator 
$$\bpar_{\na,b}^B\!: \Ga_B(b;\L) \lra \Ga_B^{0,1}(b;\L)$$
induced by $\bpar_{\na,b}$ is surjective. 
In particular, 
$$\dim\V_{1,k}^A|_{[b]}=\dim\ker\bpar_{\na,b}=\dim\ker\bpar_{\na,b}^B
=\ind\bpar_{\na,b}^B=\lr{c_1(\L),A}+1
\qquad\forall~ [b]\!\in\!\U_{\T,\nu}(X;J).$$
Thus, $\V_{1,k}^A|_{\U_{\T,\nu}(X;J)}$ is a smooth vector bundle.
Similarly to Subsection~\ref{g1cone-g1conelocalstr_subs2} in~\cite{g1cone},
for every $m$ as in~\e_ref{mbound_e} we can construct a vector bundle 
$F^1\!\lra\!\U_{\T,\nu;1}^m(X;J)$ of rank~$m$ and a surjective bundle homomorphism
$$\D\!: \V_{1,k}^A\lra\Hom(F^1,\ev_P^*\L)$$
over $\U_{\T,\nu;1}^m(X;J)$ such that the kernel of $\D$ is 
$\W_{1,k}^A|_{\U_{\T,\nu;1}^m(X;J)}$.
Thus, $\W_{1,k}^A|_{\U_{\T,\nu;1}^m(X;J)}$ is a smooth vector bundle of the claimed rank.

\section{Proof of Proposition~\ref{conestr_prp}}
\label{gluing_sec}

\subsection{Outline}
\label{gluingsumm_subs}

\noindent
In this section we prove a generalization of Proposition~\ref{conestr_prp}.
It implies that the Poincare dual of the euler class of $\V_{1,k}^A$
defined as the zero set of a generic section of~$\V_{1,k}^A$ over
$\ov\M_{1,k}^0(X,A;J,\nu)$ is independent of $J$ and~$\nu$.\\

\noindent
Suppose $(X,\om)$ is a compact symplectic manifold, 
$(\L,\na)$ is an $(\om,A)$-positive line bundle with connection over~$X$,
$A\!\in\!H_2(X;\Z)^*$, $k\!\in\!\bar\Z^+$, $\under{J}\!\equiv\!(J_t)_{t\in[0,1]}$
is a continuous family of $\om$-tamed almost complex structures on $X$, and 
$$\under{\nu}\!\equiv\!(\nu_t)_{t\in[0,1]} \in 
\G^{\es}_{1,k}(X,A;\under{J})$$
is a family of sufficiently small perturbations of the $\bpar_{J_t}$-operators
on~$\X_{1,k}(X,A)$.
Let $t_r$ and $[b_r]$ be sequences of elements in $[0,1]$ and 
in $\ov\M_{1,k}^0(X,A;J_{t_r},\nu_{t_r})$ such that 
$$\lim_{r\lra\i}t_r=0  \quad\hbox{and}\quad
\lim_{r\lra\i}[b_r]=[b]\in\ov\M_{1,k}^0(X,A;J_0,\nu_0).$$
We need to show that for every $\xi\!\in\!\W_{1,k}^A|_{[b]}$ 
there exists a sequence $\xi_r\!\in\!\W_{1,k}^A|_{[b_r]}$ converging to~$\xi$.
By the paragraph following Proposition~\ref{conestr_prp}, it is sufficient to assume 
that $[b]$ is an element of $\U_{\T,\nu_0}(X;J_0)$ for a bubble type
$$\T=\big([k],I,\aleph;j,\under{A})$$
such that $A_i\!=\!0$ for all minimal elements $i\!\in\!I$.\\

\noindent
We can also assume that for some bubble type
$$\T'=\big([k],I',\aleph';j',\under{A}')$$
$[b_r]\!\in\!\U_{\T',\nu_{t_r}}(X;J_{t_r})$ for all $r$.
We note that by Definition~\ref{pos_dfn}, for every map 
$u\!:\bP^1\!\lra\!X$ such that 
$$\lr{\om,u_*[\bP^1]}<\lr{\om,A}$$
the linear operators
$$\bpar_{\na,u}\!:\Ga(\bP^1;u^*\L)\lra\Ga\big(\bP^1;
 \La^{0,1}_{\I,j}T^*\bP^1\!\otimes\!u^*\L\big) \qquad\hbox{and}\qquad
 \ker\bpar_{\na,u}\lra\L_{u(\i)},~~\xi\lra\xi(\i),$$
are surjective.
Thus, it is sufficient to consider two possibilities for~$\T'$:\\
${}\quad$ (1) $A_i'\!=\!0$ for all $i\!\in\!I_0'$ and 
$\{i\!\in\!I'\!:A_i'\!\neq\!0\}=\chi(\T')$;\\
${}\quad$ (2) $A_i'\!\neq\!0$ for some $i\!\in\!I_0'$ and 
$\hat{I}'\!=\!\eset$,\\
where $I_0'$ is the subset of minimal elements of $I'$.
In the first case, for every $[b_r]\!\in\!\U_{\T',\nu_{t_r}}(X;J_{t_r})$,
the map $u_{b_r}$ is constant on the principal component $\Si_{b_r;P}$ of~$\Si_{b_r}$,
and thus so is every element $\xi\!\in\!\ker\bpar_{\na,b_r}$.
In this case, the question of existence of a sequence $\xi_r$
as above is an issue concerning the behavior of holomorphic bundle sections
for genus-zero $(J,\nu)$-holomorphic maps, 
for a certain class of perturbations~$\nu$ of the $\bpar_J$-operator.
This class is induced from the class of effectively supported perturbations
of Definition~\ref{pert_dfn} and is described in Definition~\ref{pert0_dfn1}
at the end of this subsection.
The existence of a desired sequence in case~(1) follows
from Lemma~\ref{cuspmap_lmm} below.
In the second case, $\Si_{b_r;P}\!=\!\Si_{b_r}$ is either a smooth torus or
a circle of spheres, depending on whether $\aleph'$ is empty or not.
There are no bubble components.
In this case, the desired result follows from Lemma~\ref{torus_lmm}.

\begin{lmm}
\label{cuspmap_lmm}
Suppose $(X,\om)$ is a compact symplectic manifold, 
$(\L,\na)$ is an $(\om,A)$-positive line bundle with connection over~$X$,
$A\!\in\!H_2(X;\Z)^*$, $M$ is a finite set,
$\under{J}\!\equiv\!(J_t)_{t\in[0,1]}$
is a continuous family of $\om$-tamed almost complex structures on $X$, and 
$$\under{\nu}\!\equiv\!(\nu_t)_{t\in[0,1]} \in 
\G^{\es}_{0,\{0\}\sqcup M}(X,A;\under{J})$$
is a family of sufficiently small perturbations of the $\bar{\partial}_{J_t}$-operators
on $\X_{0,\{0\}\sqcup M}(X,A)$.
Let $t_r$ and $[b_r]$ be sequences of elements in $[0,1]$ and 
in $\M_{0,\{0\}\sqcup M}^0(X,A;J_{t_r},\nu_{t_r})$  such that
$$\lim_{r\lra\i}t_r=0 \quad\hbox{and}\quad
\lim_{r\lra\i}[b_r]=[b]\in \ov\M_{0,\{0\}\sqcup M}(X,A;J_0,\nu_0).$$
Then there exist $(c_{r,i})_{i\in\chi(b)}\!\in\!(\C^*)^{\chi(b)}$, $\ep_r\!\in\!\R^+$,
and isomorphisms
$$R_{b_r,b}\!:  \ker\bpar_{\na,b} \lra \ker\bpar_{\na,b_r}$$
such that
\begin{gather}
\label{cuspmaplmm_e1}
\Big|\cD_{\hat{0}}b-\sum_{i\in\chi(b)}\!\!\!c_{r,i}\!\cdot_{J_0}\!\cD_ib\Big|
\le \ep_r\sum_{i\in\chi(b)}\!\!\!|c_{r,i}| \qquad\forall\,r\!\in\!\Z^+;\\
\label{cuspmaplmm_e2}
\Big|\D_{b,\hat{0}}R_{b_r,b}\xi-\sum_{i\in\chi(b)}\!\!\!c_{r,i}\D_{b,i}\xi\Big|
\le \ep_r\sum_{i\in\chi(b)}\!\!\!|c_{r,i}| \cdot\|\xi\|
\qquad\forall\,\xi\!\in\!\ker\bpar_{\na,b},\, r\!\in\!\Z^+;\\
\label{cuspmaplmm_e3}
\lim_{r\lra\i}\![R_{b_r,b}\xi]=[\xi]\in\V_{0,M}^A|_{[b]}
~~\forall\,\xi\!\in\!\ker\bpar_{\na,b}, \quad
\lim_{r\lra\i}\!\ep_r\!=0,  \quad
|c_{r,i}|\le1~~\forall\, i\!\in\!\chi(b).
\end{gather}
\end{lmm}

\begin{lmm}
\label{torus_lmm}
Suppose $(X,\om)$, $(\L,\na)$, $\under{J}$, $A$, and $M$  are as in Lemma~\ref{cuspmap_lmm} and 
$$\under{\nu}\!\equiv\!(\nu_t)_{t\in[0,1]} \in 
\G^{\es}_{1,M}(X,A;\under{J})$$
is a family of sufficiently small perturbations of the $\bar{\partial}_{J_t}$-operators
on $\X_{1,M}(X,A)$.
Let  $t_r$ and $[b_r]$ be sequences of elements in $[0,1]$ and 
in $\ov\M_{1,M}^0(X,A;J_{t_r},\nu_{t_r})$  such that
$$\Si_{b_r;P}=\Si_{b_r}~~\forall\,r, \qquad
\lim_{r\lra\i}t_r=0, \quad\hbox{and}\quad
\lim_{r\lra\i}[b_r]=[b]\in \ov\M_{1,M}^0(X,A;J_0,\nu_0).$$
If $b\!=\!(\Si_b;u_b)$ is such that the degree of $u_b|_{\Si_{b;P}}$ is zero, 
there exist $(w_i)_{i\in\chi(b)}\in\!(\C^*)^{\chi(b)}$
such that 
$$\sum_{i\in\chi(b)}\!\!w_i\!\cdot_{J_0}\!\cD_ib=0$$
and a subsequence of $\{b_r\}$, which we still denote by~$\{b_r\}$, such that
$$\xi_r\!\in\!\ker\bpar_{\na,b_r}, \quad
\lim_{r\lra\i}\![\xi_r]=[\xi]\in\V_{1,M}^A|_{[b]}
\qquad\Lra\qquad
\sum_{i\in\chi(b)}\!\!\!w_i\,\D_{b,i}\xi=0.$$\\
\end{lmm}

\noindent
Suppose $\{b_r\}$ and $b$ are as in Lemma~\ref{torus_lmm}, 
$\xi\!\in\!\ti\Ga_-(b;\L)$, and there exists no sequence 
\begin{equation}\label{toruslmm_e1}
\xi_r\in\ker\bpar_{\na,b_r} \qquad\hbox{s.t.}\qquad
\lim_{r\lra\i}\![\xi_r]=[\xi].
\end{equation}
After passing to a subsequence of $\{b_r\}$ if necessary, we can assume that
\begin{equation}\label{toruslmm_e2}
[\xi]\not\in \ov{\bigcup_{r=1}^{\i}\V_{1,M}^A|_{[b_r]}} \subset \V_{1,M}^A.
\end{equation}
Choose $\under{w}\!\equiv\!(w_i)_{i\in\chi(b)}\in\!(\C^*)^{\chi(b)}$ and 
a subsequence of $\{b_r\}$, which we still denote by $\{b_r\}$, 
as at the end of Lemma~\ref{torus_lmm}.
Let $\{\xi^j_r\}$ be a basis for $\ker\bpar_{\na,b_r}$ 
which is orthonormal, with respect to a regularization at~$b$ 
as in~(1) of Definition~\ref{coneprop_dfn} for example.
After passing to a subsequence if necessary, for some
linearly independent $\xi^j\!\in\!\ker\bpar_{\na,b}$,
$$\lim_{r\lra\i}\![\xi_r^j]=[\xi^j]\in\V_{1,M}^A|_{[b]} \qquad\forall\,j.$$
By Lemma~\ref{torus_lmm},
$$\xi^j\in \Ga_-(b;\L;[\under{w}]) \!\equiv\!\big\{\xi\!\in\!\ker\bpar_{\na,b}\!:
\sum_{i\in\chi(b)}\!\!\!w_i\,\D_{b,i}\xi\!=\!0\big\}.$$
By the positivity assumption of Definition~\ref{pos_dfn},
\begin{gather*}
\dim\ker\bpar_{\na,b_r}=\ind\bpar_{\na,b_r}=\lr{c_1(\L),A} \qquad\hbox{and}\\
\dim\Ga_-(b;\L;[\under{w}])=\dim\ker\bpar_{\na,b}-1=
\ind\bpar_{\na,b}^B-1=\lr{c_1(\L),A}.
\end{gather*}
Thus, $\{\xi^j\}$ is a basis for $\Ga_-(b;\L;[\under{w}])$.
Since $\ti\Ga_-(b;\L)\!\subset\!\Ga_-(b;\L;[\under{w}])$,
there exists a sequence $\xi_r\!\in\!\ker\bpar_{\na,b_r}$ as in~\e_ref{toruslmm_e1},
with $\{b_r\}$ replaced by a subsequence.
However, this contradicts~\e_ref{toruslmm_e2}.\\

\noindent
An element $[b_r]\!\in\!\U_{\T',\nu_{t_r}}(X;J_{t_r})$, with $\T'$ as in the first case above,
corresponds to the genus-one curve $\Si_{b_r,0}\!=\!\Si_{b_r}^0$ and genus-zero maps 
$\{b_{r,h}\}_{h\in\chi(\T')}$ such that $\Si_{b_{r,h}}\!=\!\bP^1$.
If $[b]$ is the limit of the sequence $\{b_r\}$, 
$b$ corresponds to a genus-one curve $\Si_{b,0}\!\subset\!\Si_b^0$ and genus-zero maps 
$\{b_h\}_{h\in\chi(\T')}$ such that 
$$\lim_{r\lra\i}\!\Si_{b_r,0}=\Si_{b,0}
\qquad\hbox{and}\qquad
\lim_{r\lra\i}\![b_{r,h}]=[b_h] ~~~\forall\,h\!\in\!\chi(\T').$$
Similarly, $\under\xi_r\!\in\!\ker\bpar_{\na,b_r}$ and $\xi\!\in\!\ker\bpar_{\na,b}$
correspond to $\xi_{r,h}\!\in\!\ker\bpar_{\na,b_{r,h}}$ and 
$\under\xi\!\in\!\ker\bpar_{\na,b_h}$, with $h\!\in\!\chi(\T')$, such that
$$\xi_{r,h}(y_0(b_{r,h}))=\xi_{r,h'}(y_0(b_{r,h'})) \quad\hbox{and}\quad
\xi_h(y_0(b_h))=\xi_{h'}(y_0(b_{h'})) \qquad\forall\,h,h'\!\in\!\chi(\T').$$
Furthermore,
$$\lim_{r\lra\i}\![\under\xi_r]=[\under\xi]\in\V_{1,k}^A|_{[b]} \qquad\Llra\qquad
\lim_{r\lra\i}\![\xi_{r,h}]=[\xi_h]\in\V_{1,M_h}^A|_{[b_h]} ~~~\forall\,h\!\in\!\chi(\T'),$$
where $M_h$ is the index set for the marked points of $b_{r,h}$ and~$b_h$.
We will assume that $[b_r]\!\in\!\U_{\T',\nu_{t_r}}^{m'}(X;J_{t_r})$ for some
$m'\!\in\!\Z^+$ and for all~$r$.\\

\noindent
With $(c_{r,i})_{i\in\chi(b_h)}\in\!(\C^*)^{\chi(b_h)}$ for each $h\!\in\!\chi(\T')$
as in Lemma~\ref{cuspmap_lmm}, let
$$\al_r\!: \C^{\chi(\T')}\lra \C^{\chi(b)}\!\equiv\!
\prod_{h\in\chi(\T')}\!\!\!\!\!\C^{\chi(b_h)}$$
be the injective homomorphism defined by
$$\al_r\big((w_h)_{h\in\chi(\T')}\big)
=(c_{r,i}w_h)_{i\in\chi(b_h),h\in\chi(\T')}.$$
We denote by $\ti{F}_r^1\!\subset\!\C^{\chi(b)}$ the image of
$$F_{b_r}^1\equiv\big\{(w_h)_{h\in\chi(\T')}\!\in\!\C^{\chi(\T')}\!:
\sum_{h\in\chi(\T')}\!\!\!\!w_h\!\cdot_{J_0}\!\cD_hb_r\!\equiv\!
\sum_{h\in\chi(\T')}\!\!\!\!w_h\!\cdot_{J_0}\!\cD_{\hat{0}}b_{r,h}=0\big\}$$
under~$\al_r$.
By our assumption on $b_r$, $\dim\ti{F}_r^1\!=\!m'$ for all $r$.
Let 
$$\{(v_{r,i}^l)_{i\in\chi(b)}\}_{l\in[m']} \qquad\hbox{and}\qquad
\big\{\under\xi_r^j\!\equiv\!(\xi_{r,h}^j)_{h\in\chi(\T')}\big\}$$
be orthonormal bases for~$\ti{F}_r^1$ and for
\begin{equation*}\begin{split}
\ti{\Ga}_-(b_r;\L) &\equiv\big\{\under\xi\!\equiv\!(\xi_h)_{h\in\chi(\T')}
\!\in\!\ker\bpar_{\na,b_r}\!:\\
&\qquad\qquad \sum_{h\in\chi(\T')}\!\!\!\!w_h\D_{b_r,h}\under\xi\!\equiv\!
\sum_{h\in\chi(\T')}\!\!\!\!w_h\D_{b_{r,h},\hat{0}}\xi_h=0
~\forall(w_h)_{h\in\chi(\T')}\!\in\!F_{b_r}^1\big\},
\end{split}\end{equation*}
respectively.
After passing to a subsequence if necessary, we can find 
$$\under{w}^l\!\equiv\!(w_i^l)_{i\in\chi(b)}\in\C^{\chi(b)} \qquad\hbox{and}\qquad
\under{\xi}^j\!\equiv\!(\xi^j_h)_{h\in\chi(\T')}\in\ker\bpar_{\na,b}$$
such that
\begin{alignat*}{1}
\lim_{r\lra\i}\!(v_{r,i}^l)_{i\in\chi(b)}&=\under{w}^l\in\C^{\chi(b)}
 \quad\forall\,l  \qquad\hbox{and}\\
\lim_{r\lra\i}\![\under{\xi}_r^j]&=[\under{\xi}^j] \in\V_{1,k}^A|_{[b]}
\quad\forall\,j.
\end{alignat*}
Each of the sets $\{\under{w}^l\}$ and $\{\under{\xi}^j\}$ is orthonormal and thus
linearly independent.
By Lemma~\ref{cuspmap_lmm},
\begin{alignat*}{1}
\bigg|\sum_{i\in\chi(b)}\!\!v_{r,i}^l\!\cdot\!_{J_0}\cD_ib\bigg|
\equiv\bigg|\sum_{h\in\chi(\T')}\sum_{i\in\chi(b_h)}\!\!
v_{r,i}^l\!\cdot\!_{J_0}\cD_{\hat{0}}b_h\bigg| &\le \ti\ep_r \qquad\forall\,l
\qquad\hbox{and}\\
\bigg|\sum_{i\in\chi(b)}\!\!v_{r,i}^l\D_{b,i}\under\xi^j\bigg|
\equiv\bigg|\sum_{h\in\chi(\T')}\sum_{i\in\chi(b_h)}\!\!
v_{r,i}^l\D_{b_h,\hat{0}}\xi^j_h\bigg| &\le \ti\ep_r \qquad\forall\,j,
\end{alignat*}
for some sequence $\ti\ep_r$ converging to~$0$. Thus,
\begin{alignat*}{1}
\under{w}^l&\in F_b^1\equiv\big\{(w_i)_{i\in\chi(b)}\!\in\!\C^{\chi(b)}\!:
\sum_{i\in\chi(b)}\!\!w_i\!\cdot_{J_0}\!\cD_ib\!=\!0\big\} \quad\forall\,l \qquad\hbox{and}\\
\under\xi^j&\in\Ga_-\big(b;\L;\{\under{w}^l\}_{l\in[m']}\big)
\equiv\!\big\{\xi\!\in\!\ker\bpar_{\na,b}\!:
\sum_{i\in\chi(b)}\!\!w^l_i\,\D_{b,i}\xi\!=\!0
~\forall\,l\!\in\![m']\big\} \quad\forall\,j.
\end{alignat*}
By Definition~\ref{pos_dfn},
\begin{equation*}\begin{split}
\dim\Ga_-\big(b;\L;\{\under{w}^l\}_{l\in[m']}\big) &=\dim\ker\bpar_{\na,b}-m'\\
&=\dim\ker\bpar_{\na,b_r}-m'=\dim\ti\Ga_-(b_r;\L).
\end{split}\end{equation*}
Thus, $\{\under\xi^j\}$ is a basis for $\Ga_-(b;\L;\{\under{w}^l\}_{l\in[m']})$.
Since $\under{w}^l\!\in\!F_b^1$ for all~$l$,
$$\ti\Ga_-(b;\L) \subset \Ga_-\big(b;\L;\{\under{w}^l\}_{l\in[m']}\big)
\qquad\hbox{and}\qquad
\W_{1,k}^A|_{[b]}\subset\ov{\bigcup_{r=1}^{\i}\V_{1,M}^A|_{[b_r]}} \subset \V_{1,M}^A.$$
As in the first paragraph after Lemma~\ref{torus_lmm}, this implies that for every 
$\xi\!\in\!\ti\Ga_-(b;\L)$ there exists a sequence $\xi_r\!\in\!\ti\Ga_-(b_r;\L)$
such that
$$\lim_{r\lra\i}\![\xi_r]=[\xi]\in\V_{1,k}^A|_{[b]}.$$\\

\noindent
{\it Remark:} In~\e_ref{cuspmaplmm_e1} and~\e_ref{cuspmaplmm_e2},
the differences are taken via a parallel transport along the shortest geodesic,
with respect to a metric on~$X$, between $\ev_0(b_r)$ and~$\ev_0(b)$.\\

\noindent
Lemmas~\ref{cuspmap_lmm}-\ref{torus_lmm} are proved in the next two subsections by extending
the gluing constructions of Subsections~\ref{g1comp2-g0prp_subs} 
and~\ref{g1comp2-comp1prp_subs} in~\cite{g1comp2} from $J$-holomorphic maps to
holomorphic bundle sections.
These extensions parallel constructions in Subsections~\ref{g1cone-gluing_subs2}
and~\ref{g1cone-gluing_subs3} in~\cite{g1cone}. 
In the rest of this subsection we recall the definition of
the type of perturbations~$\nu$ of the $\bpar_J$-operator on space of 
genus-zero stable maps that appears in Lemma~\ref{cuspmap_lmm}; 
see Subsection~\ref{g1comp2-notation0_subs} in~\cite{g1comp2} for details.

\begin{dfn}
\label{pert0_dfn1}
Suppose $(X,\om)$ is a compact symplectic manifold, $\under{J}\!\equiv\!(J_t)_{t\in[0,1]}$ 
is a continuous family of $\om$-tamed almost structures on~$X$,
$A\!\in\!H_2(X;\Z)^*$, and $M$ is a finite set.
A continuous family of multisections $\under{\nu}\!\equiv\!(\nu_t)_{t\in[0,1]}$,
with $\nu_t\!\in\!\G_{0,\{0\}\sqcup M}^{0,1}(X,A;J_t)$ for all $t\!\in\![0,1]$, 
is  {\tt effectively supported} if for every element 
$b$ of $\X_{0,\{0\}\sqcup M}(X,A)$ 
there exists a neighborhood $\W_b$ of $\Si_b^0$ in a semi-universal family of
deformations for~$b$ such that
$$\nu_t(b')\big|_{\Si_{b'}\cap\W_b}=0 \qquad
\forall~~  [b']\in\X_{0,\{0\}\sqcup M}(X,A),~t\!\in\![0,1].$$\\
\end{dfn}

\noindent
We denote the space of effectively supported families $\under\nu$ as 
in Definition~\ref{pert0_dfn1} by $\G^{\es}_{0,\{0\}\sqcup M}(X,A;\under{J})$.
If  $\nu\!\in\!\G^{\es}_{0,\{0\}\sqcup M}(X,A;\under{J})$, $t\!\in\![0,1]$, 
$[b]$ is an element of 
$$\ov\M_{0,\{0\}\sqcup M}(X,A;J_t,\nu_t) \equiv
\big\{\bar\partial_J\!+\!\nu_t\big\}^{-1}(0),$$
and $i\!\in\!\chi(b)$, then $u_b|_{\Si_{b,i}}$ is $J_t$-holomorphic 
on a neighborhood of $\i$ in~$\Si_{b,i}$ and
$\C\!\cdot_{J_t}\!\cD_ib$ is determined by~$b$, just as in Subsection~\ref{review_subs}.
Furthermore, in this case $u_b|_{\Si_b^0}$ is a degree-zero holomorphic map and thus is constant.
Thus, $u_b$  maps the attaching nodes of all elements of $\chi(b)$
to the same point in~$X$, as in the genus-one case of Subsection~\ref{review_subs}.

\subsection{Proof of Lemma~\ref{cuspmap_lmm}}
\label{gluing0_subs}

\noindent
In this subsection we review the genus-zero gluing construction  of 
Subsection~\ref{g1comp2-g0prp_subs} in~\cite{g1comp2} and extend it 
to holomorphic bundle sections in a manner similar to 
Subsection~\ref{g1cone-gluing_subs2} in~\cite{g1cone}.
This construction essentially constitutes the first step of the two-step gluing
construction described in Subsection~\ref{gluing1_subs}.
Throughout this subsection we assume that $M$ is a finite set, $A\!\in\!H_2(X;\Z)$, and 
$\T\!=\!(M,I;j,\under{A})$ is a bubble such that $\hat{0}$ is the minimal element of~$I$,
$$\sum_{i\in I}A_i=A \qquad\hbox{and}\qquad \lr{\om,A_i}\ge 0~~~\forall\,i\!\in\!I.$$
Let $(\L,\na)\!\lra\!X$ be an $(\om,A)$-positive line bundle with connection.\\

\noindent
We put
$$\ti\X_{\T;B}(X)=\big\{(\Si_b;u_b)\!\in\!\ti\X_{\T}(X)\!: 
u_b|_{\Si_b^0}\!=\!\const\big\}.$$
We denote by
$$\wt\cF\equiv \ti\X_{\T;B}(X)\times\C^{\hat{I}}$$
the bundle of smoothing parameters and by $\wt\cF^{\eset}$ the subset 
of $\wt\cF$ consisting of the elements with all components nonzero.
For each 
$$b\in\ti\X_{\T;B}(X), \qquad i\in\chi(\T)\!\equiv\!\chi(b), \qquad
\ups\!\equiv\!(b,v)\!=\!\big(b,(v_h)_{h\in\hat{I}}\big)\in\wt\cF,$$
we~put
$$\rho_i(\ups)=\prod_{\hat{0}<h\le i}\!\!\!v_h\in\C
\qquad\hbox{and}\qquad 
x_i(\ups)=\!\!\sum_{\hat{0}<i'\le i}\!\Big(x_{i'}(b)
\!\!\!\prod_{\hat{0}<h<i'}\!\!\!\!v_h\Big)~\in~\C,$$
where $x_i(b)$ is the point of $\Si_{b,\io_i}$ to which the bubble $\Si_{b,i}$ is
attached; see \e_ref{stablemap_e2} and Figure~\ref{derivest_fig}.\\

\noindent
For each sufficiently small element $\ups\!=\!(b,v)$ of $\wt\cF^{\eset}$, let 
$$q_{\ups}\!:\Si_{\ups}\lra\Si_b$$
be the basic gluing map constructed in Subsection~2.2 of~\cite{gluing}.
In this case, $\Si_{\ups}$ is the projective line $\bP^1$ with $|M|\!+\!1$ marked points.
The map~$q_{\ups}$ collapses $|\hat{I}|$ circles on~$\Si_{\ups}$.
It induces a metric $g_{\ups}$ on~$\Si_{\ups}$ such that
$(\Si_{\ups},g_{\ups})$ is obtained from $\Si_b$ by replacing 
the $|\hat{I}|$ nodes of $\Si_b$ by thin necks.\\

\noindent
We put
$$u_{\ups}=u_b\!\circ\!q_{\ups}, \qquad b(\ups)=(\Si_{\ups};u_{\ups}), 
\quad\hbox{and}\quad \bpar_{\na,\ups}=\bpar_{\na,b(\ups)}.$$
Fix a metric $g$ on $X$ and denote 
the corresponding Levi-Civita connection by~$\na^X$.
By the same construction as in Subsection~3.3 of~\cite{gluing},
the map~$q_{\ups}$ induces weighted $L^p_1$-norms $\|\cdot\|_{\ups,p,1}$ on the spaces
$$\big\{\ze\!\in\!\Ga(\Si_{\ups};u_{\ups}^*TX)\!: \ze(\i)\!=\!0\big\}  
\qquad\hbox{and}\qquad  \Ga(\Si_{\ups};u_{\ups}^*\L)$$
and a weighted $L^p$-norm $\|\cdot\|_{\ups,p}$ on the space
$\Ga(\Si_{\ups};\La^{0,1}_{\I,j}T^*\Si_{\ups}\!\otimes\!u_{\ups}^*\L)$.
We denote the corresponding completions by $\Ga(\ups)$, $\Ga(\ups;\L)$,
and~$\Ga^{0,1}(\ups;\L)$.
The norms $\|\cdot\|_{\ups,p,1}$ and $\|\cdot\|_{\ups,p}$
are analogous to the ones used in Section~3 of~\cite{LT1} for the bundle~$TX$.
We~put
\begin{alignat}{1}
\Ga_-(\ups;\L)&=
\big\{\xi\!\circ\!q_{\ups}\!: \xi\!\in\!\ker\bpar_{\na,b}\big\}
\subset \Ga(\ups;\L);\notag\\
\label{g0reg_lmm2e0}
\Ga_+(\ups;\L)&=\big\{\xi\!\in\!\Ga(\ups;\L)\!: \xi(\i)\!=\!0;~
\llrr{\xi,\xi'}_{\ups,2}\!=\!0~\forall~\xi'\!\in\!\Ga_-(\ups;\L)
~\st~\xi'(\i)\!=\!0\big\}.
\end{alignat}
By the construction of the map~$q_{\ups}$ in 
Subsection~2.2 of~\cite{gluing},
\begin{equation}\label{g0reg_lmm2e1}
\big\|\bpar_{\na,\ups}\xi\big\|_{\ups,p} \le C(b) |\ups|^{1/p}\|\xi\|_{\ups,p,1}
 \qquad \forall~\xi\!\in\!\Ga_-(\ups;\L).
\end{equation}
On the other hand, for the same reasons as in Section~3 of~\cite{LT1}, for some
$\de,C\!\in\!C(\ti\X_{\T;B}(X);\R^+)$ and
for all $\ups\!=\!(b,v)\!\in\!\ti\cF_{\de}^{\eset}$,
\begin{equation}\label{g0reg_lmm1e1}
C(b)^{-1}\big\|\xi\big\|_{\ups,p,1} \le  \big\|\bpar_{\na,\ups}\xi\big\|_{\ups,p} 
\le C(b)\big\|\xi\big\|_{\ups,p,1}
\qquad\forall~\xi\in\Ga_+(\ups;\L);
\end{equation}
see Subsections~5.3 and~5.4 in~\cite{gluing}.
In particular, the operators $\bpar_{\na,\ups}$ are surjective,
since $\bpar_{\na,b}$ and  $\bpar_{\na,b}^{y_0(b)}$ are.\\

\noindent
If $\ze\!\in\!\Ga(\ups)$, we~set
$$b(\ups,\ze)=\big(\Si_{\ups};\exp_{u_{\ups}}\!\ze)
\qquad\hbox{and}\qquad u_{\ups,\ze}=\exp_{u_{\ups}}\!\ze,$$
where $\exp$ is the exponential map with respect to the connection $\na^X$.
Let
$$\Pi_{\ups,\ze}\!: \Ga(\Si_{\ups};u_{\ups}^*\L)
\lra \Ga(\Si_{\ups};u_{\ups,\ze}^*\L)$$
be the isomorphism induced by the $\na$-parallel transport along the $\na^X$-geodesics
$\tau\!\lra\!\exp_{\ups}\!\tau\ze$ with $\tau\!\in\![0,1]$.
By a direct computation,
\begin{equation}\label{g0reg_lmm1e2}
\big\|\Pi_{\ups,\ze}^{-1}\bpar_{\na,b(\ups,\ze)}\Pi_{\ups,\ze}\xi
-\bpar_{\na,\ups}\xi\big\|_{\ups,p}
\le C(b)\|\ze\|_{\ups,p,1}^2\|\xi\|_{\ups,p,1}
\qquad\forall~\xi\!\in\!\Ga(\ups;\L), ~\ze\!\in\!\Ga(\ups);
\end{equation}
see the proof of Corollary~2.3 in~\cite{analyt}.
By \e_ref{g0reg_lmm2e1}-\e_ref{g0reg_lmm1e2}, 
for every $\ze\!\in\!\Ga(\ups)$ sufficiently small and $\xi\!\in\!\ker\bpar_{\na,b}$, 
there exists a unique
$$\xi_{\ups,\ze}(\xi)\in\Ga_+(\ups;\L) \qquad\st\qquad
\Pi_{\ups,\ze}\big(\xi\!\circ\!q_{\ups}\!+\!\xi_{\ups,\ze}(\xi)\big)
\in \ker\bpar_{\na,b(\ups,\ze)}.$$
Furthermore, 
\begin{equation}\label{g0reg_lmm2e2}
\big\|\xi_{\ups,\ze}(\xi)\big\|_{\ups,p,1} 
\le C(b)\big(|\ups|^{1/p}\!+\!\|\ze\|_{\ups,p,1}^2\big)\|\xi\|_{\ups,p,1}
\qquad \forall~\xi\!\in\!\ker\bpar_{\na,b}.
\end{equation}
We define the isomorphism 
$$R_{\ups,\ze}\!: \ker\bpar_{\na,b} \lra \ker\bpar_{\na,b(\ups,\ze)}
\qquad\hbox{by}\qquad
R_{\ups,\ze}\xi=\Pi_{\ups,\ze}
\big(\xi\!\circ\!q_{\ups}\!+\!\xi_{\ups,\ze}(\xi)\big).$$\\

\noindent
We will use a convenient family of connections in the vector bundles 
$u^*\L\!\lra\!\Si$, which is provided by Lemma~\ref{conn_lmm} below. 
First, if $b\!=\!(\Si_b;u_b)$ is a stable bubble map,
$g_b$ is a Hermitian metric in the vector bundle $u_b^*\L\!\lra\!\Si_b$,
and  $\na^b$ is a connection in~$u_b^*\L$, we will call the pair
{\tt $(g,\na)$-admissible}~if\\
${}\quad$ ($g\na1$) $\na^b$ is $g_b$-compatible and $\bpar_{\na,b}$-compatible;\\
${}\quad$ ($g\na2$) $g_b\!=\!g_{u_b}$ and $\na^b\!=\!\na^{u_b}$ on $\Si_b^0$,\\
where $g_{u_b}$ is the Hermitian metric in $u_b^*\L$ induced
from the standard metric in~$\L$.
The second condition in~($g\na1$) means~that
$$\bpar_{\na,b}
\equiv\frac{1}{2}\big(\na^{u_b}+\I\na^{u_b}\circ j\big)
=\frac{1}{2}\big(\na^b+\I\na^b\circ j\big),$$ 
where $\na^{u_b}$ is as in Subsection~\ref{eulercone_subs}.\\

\noindent
If $b$ is any genus-zero bubble map and $\de\!\in\!\R^+$, we put
\begin{equation}\label{mainpart_e1}
\Si_b^0(\de)=\Si_b^0 \cup \!\bigcup_{i\in\chi(b)} \!\!\! A_{b,i}(\de),
\quad\hbox{where}\quad
A_{b,i}(\de)=\big\{(i,z)\!:|z|\!\ge\!\de^{-1/2}/2\big\}\subset \Si_{b',i}\!\approx\!S^2.
\end{equation}
If $\ups\!=\!(b,v)$ is a gluing parameter such that the map $q_{\ups}$ is defined, let
\begin{equation}\label{mainpart_e2}
\Si_{\ups}^0(\de)=q_{\ups}^{-1}\big(\Si_b^0(\de)\big).
\end{equation}

\begin{lmm}
\label{conn_lmm}
If $(X,\om)$, $(\L,\na)$, $A$, $M$, $\T$, $\ti\X_{\T;B}(X)$, $\wt\cF$ are as above,
there exist 
$$\de,C\in C(\ti\X_{\T;B}(X);\R^+)$$
with the following property.
For every 
$$b\in\ti\X_{\T;B}(X), \qquad \ups\!=\!(b,v)\in\wt\cF^{\eset}_{\de}, 
\quad\hbox{and}\quad \ze\in\Ga(\ups)~~\st~~\|\ze\|_{\ups,p,1}\le\de(b)$$
there exist metrics $g_b$ and $g_{(\ups,\ze)}$ and 
connections $\na^b$ and $\na^{(\ups,\ze)}$ in the vector bundles
$$u_b^*\L \lra \Si_b \qquad\hbox{and}\qquad u_{\ups,\ze}^*\L\lra\Si_{\ups}$$
such~that\\
${}\quad$ (1) all pairs $(g_b,\na^b)$ and 
$(g_{(\ups,\ze)},\na^{(\ups,\ze)})$ are admissible;\\
${}\quad$ (2) the curvatures of $\na^b$ and $\na^{(\ups,\ze)}$ vanish on 
$\Si_b^0(2\de)$ and $\Si_{\ups}^0(2\de)$, respectively;\\ 
${}\quad$ (3) $\|\Pi_{\ups,\ze}^{-1}\na^{(\ups,\ze)}(\Pi_{\ups,\ze}\xi\!\circ\!q_{\ups})
-\na^b\xi\circ dq_{\ups}\|_{\ups,p}
\le C\big(|\ups|^{1/p}\!+\!\|\ze\|_{\ups,p,1}\big)\|\xi\|_{b,p,1}$ 
for all $\xi\!\in\!\Ga(b;\L)$;\\
${}\quad$ (4) the map $b\!\lra\!(g_b,\na^b)$ is continuous.
\end{lmm}

\noindent
This lemma is proved by exactly the same argument as 
Lemma~\ref{g1cone-conn_lmm} in~\cite{g1cone}.\\

\noindent
If $b\!\in\!\ti\X_{\T;B}(X)$ and  $i\!\in\!\chi(b)$, let $w_i$ be
the standard holomorphic coordinate centered at the point $\i$ in $\Si_{b,i}\!=\!S^2$.
If $m\!\in\!\Z^+$ and $\xi\!=\!(\xi_i)_{i\in I}\!\in\!\Ga(b;\L)$, we~put
$$\D_{b,i}^{(m)}\xi = \frac{1}{m!}
\frac{d^m}{dw_i^m} \xi_i(w_i)\Big|_{w_i=0} \in \L_{\ev_0(b)},$$
where the derivatives are taken with respect to the connection $\na^b$.
Similarly, for all 
$$\ups\!\!=\!(b,v)\in\wt\cF^{\eset}_{\de}, \qquad
\ze\in\Ga(\ups) ~~\hbox{s.t.}~~\|\ze\|_{\ups,p,1}\!\le\!\de(b), 
\quad\hbox{and}\quad \xi\in\Ga(\Si_{\ups};u_{\ups,\ze}^*\L),$$ 
let
$$\D_{(\ups,\ze),\hat{0}}^{(m)}\xi = \frac{1}{m!}
\frac{d^m}{dw^m} \xi(w)\Big|_{w=0} \in \L_{\ev_0(b)},$$
where $w$ is the standard holomorphic coordinate centered 
at the point $\i$ in $\Si_{\ups}\!\approx\!S^2$
and the derivatives are taken with to the connection $\na^{(\ups,\ze)}$.
We note that
\begin{equation}\label{g0reg_lmm4e0}
\D_{b,i}^{(1)}\xi=\D_{b,i}\xi \quad\forall~ i\!\in\!\chi(b),~\xi\!\in\!\Ga(b;\L),
\quad\hbox{and}\quad    \D_{(\ups,\ze),\hat{0}}^{(1)}\xi=\D_{b(\ups,\ze),\hat{0}}\xi 
\quad\forall~ \xi\!\in\!\Ga(\Si_{\ups};u_{\ups,\ze}^*\L),
\end{equation}
by the second condition  in~($g\na2$) above.\\

\noindent
A key step in understanding the obstruction to extending holomorphic bundle sections
from singular to smooth domains is the following power series expansion.
For every 
$$i\!\in\!\chi(b), \qquad k\!\in\!\Z^+, \qquad 
\ups\!=\!(b,v)\in\wt\cF^{\eset}_{\de}, \quad\hbox{and}\quad
\ze\!\in\!\Ga(\ups),$$
there exist 
$$\ve_i^{(k)}(\ups,\ze)\in\Hom\big(\ker\bpar_{\na,b},\L_{\ev_0(b)}\big)$$
such that for all $\xi\!\in\!\ker\bpar_{\na,b}$
\begin{equation}\label{derivexp_e1}
\D_{b(\ups,\ze),\hat{0}}^{(m)}R_{\ups,\ze}\xi =
\sum_{k=1}^{k=m} \binom{m\!-\!1}{k\!-\!1}\sum_{i\in\chi(b)} \!\!\!
x_i^{m-k}(\ups)\rho_i^k(\ups) \big\{ \D_{b,i}^{(k)}\xi\!+\!\ve_i^{(k)}(\ups,\ze)\xi\big\}
\end{equation}
and 
\begin{equation}\label{derivexp_e2}
\big|\ve_i^{(k)}(\ups,\ze)\xi\big| \le C\de^{-k/2}
\big(|\ups|^{1/p}\!+\!\|\ze\|_{\ups,p,1}\big)\|\xi\|_{b,p,1}. 
\end{equation}
The expansion~\e_ref{derivexp_e1} is obtained by exactly the same integration-by-parts
argument as the expansion in Theorem~2.8 of~\cite{gluing}; see also the paragraph 
following Lemma~\ref{g1cone-derivest_lmm} in~\cite{g1cone}.
We point out that $\ve_i^{(k)}$ is independent of~$m$.
The $m\!=\!1$ case of the estimates~\e_ref{derivexp_e1} 
and~\e_ref{derivexp_e2} is illustrated in Figure~\ref{derivest_fig}.\\

\begin{figure}
\begin{pspicture}(-1.1,-1.8)(10,1.25)
\psset{unit=.4cm}
\pscircle(5,-1.5){1.5}\pscircle*(5,-3){.24}\rput(5,-3.8){$(\hat{0},\i)$}
\pscircle[fillstyle=solid,fillcolor=gray](3.23,.27){1}\pscircle*(3.94,-.44){.2}
\pscircle[fillstyle=solid,fillcolor=gray](1.23,.27){1}\pscircle*(2.23,.27){.2}
\pscircle(6.77,.27){1}\pscircle*(6.06,-.44){.2}
\pscircle[fillstyle=solid,fillcolor=gray](6.77,2.27){1}\pscircle*(6.77,1.27){.2}
\pscircle[fillstyle=solid,fillcolor=gray](8.77,.27){1}\pscircle*(7.77,.27){.2}
\rput(3.3,1.8){$h_1$}\rput(1.3,1.8){$h_2$}\rput(5.3,1){$h_3$}
\rput(8.4,2.4){$h_4$}\rput(10.4,.4){$h_5$}
\pnode(2,-3){A1}\rput(.75,-3){$x_{h_1}(b)$}\pnode(4,-.46){B1}
\ncarc[nodesep=.35,arcangleA=-30,arcangleB=-50,ncurv=.8]{->}{A1}{B1}
\pnode(7.3,-3){A2}\rput(8.5,-2.9){$x_{h_3}(b)$}\pnode(6,-.45){B2}
\ncarc[nodesep=.35,arcangleA=30,arcangleB=40,ncurv=.8]{->}{A2}{B2}
\pnode(8.8,-1.55){A3}\rput(10,-1.4){$x_{h_5}(b)$}\pnode(7.85,.25){B3}
\ncarc[nodesep=.35,arcangleA=45,arcangleB=110,ncurv=1]{->}{A3}{B3}
\rput(24.9,0){\begin{small}\begin{tabular}{l}
$\chi(b)\!=\!\{h_1,h_4,h_5\}$\\
$\rho(\ups)=(v_{h_1},v_{h_3}v_{h_4},v_{h_3}v_{h_5})$\\
$x_{h_5}(\ups)=x_{h_3}(b)+v_{h_3}x_{h_5}(b)$\\
\\
$\D_{b(\ups,\ze),\hat{0}}^{(1)}R_{\ups,\ze}\xi
\cong v_{h_1}\big(\D_{b,h_1}^{(1)}\xi\big)+
v_{h_3}v_{h_4}\big(\D_{b,h_4}^{(1)}\xi\big)
 + \!v_{h_3}v_{h_5}\big(\D_{b,h_5}^{(1)}\xi\big)$
\end{tabular}\end{small}}
\end{pspicture}
\caption{An Example of the Estimates \e_ref{derivexp_e1} and~\e_ref{derivexp_e2}}
\label{derivest_fig}
\end{figure}

\noindent
Let $t_r$ and $b_r$ be as in Lemma~\ref{cuspmap_lmm}.
Since the sequence $[b_r]$ converges to $[b]$, for all $r$ sufficiently large
there exist\\
$$b_r'\in\ti\X_{\T;B}(X), \qquad
\ups_r=(b_r',v_r)\in\wt{\cal F}^{\eset}, \qquad\hbox{and}\qquad
\ze_r\in\Ga(\ups_r)$$
such that
\begin{gather}\label{cuspmap_e0}
\lim_{r\lra\i}b_r'=b, \qquad \lim_{r\lra\i}|v_r|=0, \qquad
\lim_{r\lra\i}\|\ze_r\|_{\ups_r,p,1}=0,\\
\hbox{and}\qquad b_r\!\equiv\!\big(\Si_{b_r};u_{b_r}\big)=
\big(\Si_{\ups_r};\exp_{u_{\ups_r}}\!\ze_r\big).\notag
\end{gather}
The last equality holds for a representative $b_r$ for $[b_r]$.
By the estimate~\e_ref{g1comp2-cuspmap_e4} in~\cite{g1comp2} and~\e_ref{cuspmap_e0},
\begin{equation}\label{cuspmap_e4a}
\Big|\cD_{\hat{0}}b_r-
\sum_{i\in\chi(b)}\!\!\rho_i(\ups_r)\big(\cD_ib_r'\big)\Big|
\le \ep_r\sum_{i\in\chi(b)}\!\!\big|\rho_i(\ups_r)\big|
\end{equation}
for a subsequence $\ep_r$ converging to~$0$.
Furthermore, by the $m\!=\!1$ case of \e_ref{derivexp_e1} and~\e_ref{derivexp_e2},
\begin{equation}\label{cuspmap_e4b}
\Big|\D_{b,\hat{0}}R_{\ups_r,\ze_r}\xi-
\sum_{i\in\chi(b)}\!\!\!\rho_i(\ups_r)\D_{b_r',i}\xi\Big|
\le \ep_r\sum_{i\in\chi(b)}\!\!\!|\rho_i(\ups_r)| \cdot\|\xi\|_{b,p,1}
\qquad\forall\,\xi\!\in\!\ker\bpar_{\na,b_r'}.
\end{equation}
For the purposes of Lemma~\ref{cuspmap_lmm}, we take 
$$R_{b_r,b}\!:\ker\bpar_{\na,b}\lra \ker\bpar_{\na,b_r}$$
to be the composition of 
$$R_{\ups_r,\ze_r}\!:\ker\bpar_{\na,b_r'}\lra \ker\bpar_{\na,b_r}$$
with an isomorphism 
$$R_{b_r',b}\!:\ker\bpar_{\na,b}\lra \ker\bpar_{\na,b_r'}
\qquad\hbox{s.t.}\qquad
\lim_{r\lra\i}R_{b_r',b}\xi=\xi ~~\forall\,\xi\!\in\!\ker\bpar_{\na,b}.$$
We take $c_{r,i}\!=\!\rho_i(\ups_r)$.
It is immediate that the requirements~\e_ref{cuspmaplmm_e3} of Lemma~\ref{cuspmap_lmm}
are satisfied.
Since 
$$\lim_{r\lra\i}\!\cD_ib_r'=\cD_ib  \quad\hbox{and}\quad
\lim_{r\lra\i}\!\D_{b_r',i}R_{b_r',b}\xi=\D_{b,i}\xi
\qquad\forall\, \xi\!\in\!\ker\bpar_{\na,b},\, i\!\in\!\chi(b),$$
the requirements~\e_ref{cuspmaplmm_e1} and~\e_ref{cuspmaplmm_e2} 
are satisfied as well. This concludes the proof of Lemma~\ref{cuspmap_lmm}.\\

\noindent
{\it Remark:} 
A regularization $\vph_b$ of the cone ${\cal V}_{0,M}^A$ near $b$, 
as in Definition~\ref{coneprop_dfn}, can be constructed using the isomorphisms~$R_{\ups,\ze}$
as above and a description of open subsets in $\X_{0,M}(X,A)$ as in Section~3 of~\cite{LT1}.
In this case, $\vph_b$ is a homeomorphism.

\subsection{Proof of Lemma~\ref{torus_lmm}}
\label{gluing1_subs}

\noindent
In this subsection, we review the two-step gluing construction
used in Subsection~\ref{g1comp2-comp1prp_subs} of~\cite{g1comp2} and extend it 
to holomorphic bundle sections in a way similar  to 
Subsection~\ref{g1cone-gluing_subs2} in~\cite{g1cone}.
We assume that 
$$\T=(M,I,\aleph;j,\under{A})$$ 
is a bubble type such that 
$$\sum_{i\in I}A_i=A,\qquad 
\lr{\om,A_i}\ge 0~~\forall\,i\!\in\!I, \qquad
A_i\!=\!0~~\forall~i\!\in\!I_0,  \quad\hbox{and}\quad I_0\!\subsetneq\!I,$$
where $I_0$ is the subset of minimal elements of~$I$  as before.
Let $(\L,\na)\!\lra\!X$ be an $(\om,A)$-positive vector bundle with connection.
Throughout this subsection we focus on the case 
$$[b_r]\in\M_{1,M}^0(X,A;J,\nu),$$
i.e.~$\Si_{b_r;P}\!=\!\Si_{b_r}$ is a smooth torus for all~$r$.\\

\noindent
Similarly to Subsection~\ref{gluing0_subs}, we put
$$\ti\X_{\T;B}(X)=\big\{(\Si_b;u_b)\!\in\!\ti\X_{\T}(X)\!: 
u_b|_{\Si_b^0}\!=\!\const\big\}.$$
Let
$$I_1=\big\{h\!\in\!\hat{I}\!:\io_h\!\in\!I_0\big\}  \qquad\hbox{and}\qquad
k_0=|I_1|+\big|\{l\!\in\!M\!:j_l\!\in\!I_0\}\big|.$$
We denote by $\ov\cM_{1,k_0}$ the moduli space of genus-one curves with $k_0$ marked 
points and by
\begin{alignat*}{2}
&\pi_P\!:\ti\X_{\T;B}(X)\lra\ov\cM_{1,k_0}, &\qquad &b\lra[\Si_{b;P}],\\
&\ev_P\!:\ti\X_{\T;B}(X)\lra X, &\qquad &b\lra u_b(\Si_{b;P}),
\end{alignat*}
the maps sending each element $b$ of $\ti\X_{\T;B}(X)$ to the equivalence class
of the principal component(s) $\Si_{b;P}$ of its domain and to the image of~$\Si_{b;P}$ 
in~$X$.
Let 
$$\E\lra\ov\cM_{1,k_0}$$
be the Hodge line bundle, i.e.~the line bundle of holomorphic differentials.\\

\noindent
Let
$$\wt\cF\lra \ti\X_{\T;B}(X)$$
be the bundle of gluing parameters.
It has three distinguished components:
\begin{gather*}
\wt\cF=\wt\cF_{\aleph} \oplus \wt\cF_0 \oplus \wt\cF_1, 
\qquad\hbox{where}\qquad
\wt\cF_{\aleph}=\ti\X_{\T;B}(X)\times\C^{\aleph}, \quad
\wt\cF_1=\ti\X_{\T;B}(X)\times\C^{\hat{I}-I_1},\\
\wt\cF_0=\bigoplus_{h\in I_1}\wt\cF_h, \quad\hbox{and}\quad
\wt\cF_h\big|_b=T_{x_h(b)}\Si_{b;P}~~\forall~b\!\in\!\ti\X_{\T;B}(X).
\end{gather*}
The total space of $\wt\cF_h$ has a natural topology;
see Subsection~\ref{g1comp-notation1_subs} in~\cite{g1comp}.
We denote by $\wt\cF^{\eset}$ the subset of  $\ti\cF$ consisting
of the elements with all components nonzero.
If $i\!\in\!\hat{I}$, let $h(i)\!\in\!I_1$ be the unique element such that
$h(i)\!\le\!i$.
For each $\ups\!=\!(b,v)$, where $b\!\in\!\ti\X_{\T;B}(X)$
and $v\!=\!(v_i)_{i\in\aleph\sqcup\hat{I}}$, and $i\!\in\!\chi(b)$, we~put
\begin{gather*}
\ups_0=\big(b,(v_i)_{i\in\aleph\sqcup I_1}\big), \qquad
\ups_1=\big(b,(v_i)_{i\in\hat{I}- I_1}\big),\\
\ti\rho_i(\ups)=\!\prod_{h(i)< h\le i}\!\!\!\!\!\! v_h \in\C, \quad
\rho_i(\ups)=\ti\rho_i(\ups)\cdot v_{h(i)}\in T_{x_{h(i)}(b)}\Si_{b;P},\\
\hbox{and}\qquad 
\rho(\ups)=\big(b,\big(\rho_i(\ups)\big)_{i\in\chi(b)}\big)
\in \wt\F \equiv \bigoplus_{i\in\chi(b)}\cF_{h(i)}.
\end{gather*}
The component $\ups_1$ of $\ups$ consists of the smoothings of the nodes 
of $\Si_b$ that lie away from the principal component.
In the case of Figure~\ref{str_fig}, 
these are the attaching nodes of the bubbles $h_2$, $h_4$, and~$h_5$.
For each element $\ti{\ups}\!=\!(b,(\ti{v}_i)_{i\in\chi(b)})$ of $\wt\F$,
we define the linear map
$$\D_{\ti\ups}\!: \ker\bpar_{\na,b} \lra \E_{\pi_P(b)}^*\!\otimes\!\L_{\ev_P(b)}
\qquad\hbox{by}\qquad
\big\{\D_{\ti\ups}\xi\big\}(\psi)=\sum_{i\in\chi(b)}\!\!\!
\psi_{x_{h(i)}(b)}\big(\ti{v}_i\big)\cdot\D_{b,i}\xi,$$
if $\psi\!\in\!\E_{\pi_P(b)}$.\\

\begin{figure}
\begin{pspicture}(-1.1,-1.8)(10,1.25)
\psset{unit=.4cm}
\psellipse(5,-1.5)(1.5,2.5)
\psarc[linewidth=.05](3.2,-1.5){2}{-30}{30}\psarc[linewidth=.05](6.8,-1.5){2}{150}{210}
\pscircle[fillstyle=solid,fillcolor=gray](2.5,-1.5){1}\pscircle*(3.5,-1.5){.2}
\pscircle[fillstyle=solid,fillcolor=gray](.5,-1.5){1}\pscircle*(1.5,-1.5){.2}
\pscircle(7.5,-1.5){1}\pscircle*(6.5,-1.5){.2}
\pscircle[fillstyle=solid,fillcolor=gray](8.91,-.09){1}\pscircle*(8.21,-.79){.2}
\pscircle[fillstyle=solid,fillcolor=gray](8.91,-2.91){1}\pscircle*(8.21,-2.21){.2}
\rput(2.5,0){$h_1$}\rput(.5,0){$h_2$}\rput(7.3,0){$h_3$}
\rput(10.5,0.1){$h_4$}\rput(10.5,-2.9){$h_5$}
\rput(25,-1.5){\begin{small}\begin{tabular}{l}
$\chi(b)\!=\!\{h_1,h_4,h_5\}$,\\
\\
$\rho(\ups)\!=\!(b;v_{h_1},v_{h_3}v_{h_4},v_{h_3}v_{h_5})$\\
\end{tabular}\end{small}}
\end{pspicture}
\caption{An Example of $\rho(\ups)$} 
\label{str_fig}
\end{figure}

\noindent
For each sufficiently small element $\ups\!=\!(b,v)$ of $\wt\cF^{\eset}$, let 
$$q_{\ups_1}\!:\Si_{\ups_1}\lra\Si_b$$
be the basic gluing map constructed in Subsection~2.2 of~\cite{gluing}.
In this case, the  principal component $\Si_{\ups_1;P}$ of $\Si_{\ups_1}$
is the same as principal component $\Si_{b;P}$ of~$\Si_b$,
and $\Si_{\ups_1}$ has $|I_1|$ bubble components $\Si_{\ups_1,h}$,
with $h\!\in\!I_1$, attached directly to~$\Si_{\ups_1;P}$.
The map~$q_{\ups_1}$ collapses $|\hat{I}\!-\!I_1|$ circles on the bubble components
of~$\Si_{\ups_1}$.
It induces a metric $g_{\ups_1}$ on~$\Si_{\ups_1}$ such that
$(\Si_{\ups_1},g_{\ups_1})$ is obtained from $\Si_b$ by replacing 
$|\hat{I}\!-\!I_1|$ nodes by thin necks.
Let
$$u_{\ups_1}=u_b\circ q_{\ups_1}, \qquad
b(\ups_1)=(\Si_{\ups_1};u_{\ups_1}), \quad\hbox{and}\quad
\bpar_{\na,\ups_1}=\bpar_{\na,b(\ups_1)}.$$
The map~$q_{\ups_1}$ induces weighted $L^p_1$-norms $\|\cdot\|_{\ups_1,p,1}$ 
on the spaces 
$$\big\{\ze\!\in\!\Ga(\Si_{\ups_1};u_{\ups_1}^*TX)\!:\ze|_{\Si_{\ups_1;P}}\!=\!0\big\}
\qquad\hbox{and}\qquad
\big\{\xi\!\in\!\Ga(\Si_{\ups_1};u_{\ups_1}^*\L)\!:\xi|_{\Si_{\ups_1;P}}\!=\!0\big\}$$
and a weighted $L^p$-norm $\|\cdot\|_{\ups_1,p}$ on 
$$\big\{\eta\!\in\!\Ga(\Si_{\ups_1};\La^{0,1}_{\I,j}T^*\Si_{\ups_1}\!\otimes\!u_{\ups_1}^*\L)\!:
\eta|_{\Si_{\ups_1;P}}\!=\!0\big\};$$
see Subsection~3.3 of~\cite{gluing} and the first remark in 
Subsection~\ref{g1comp2-comp1prp_subs} in~\cite{g1comp2}.
We denote the corresponding completions by $\Ga_B(\ups_1)$, $\Ga_B(\ups_1;\L)$, 
and~$\Ga_B^{0,1}(\ups_1;\L)$.\\

\noindent 
For each $\ze\!\in\Ga_B(\ups_1)$, let 
$$u_{\ups_1,\ze}=\exp_{u_{\ups_1}}\!\ze \qquad\hbox{and}\qquad 
b(\ups_1,\ze)=(\Si_{\ups_1};u_{\ups_1,\ze}).$$
For $\de\in\!C(\ti\X_{\T;B}(X);\R^+)$ sufficiently small,
$$\ups\!\equiv\!(b,v)\in\ti\cF_{\de}^{\eset},  \qquad\hbox{and}\qquad
\ze\!\in\!\Ga_B(\ups_1)~~\hbox{s.t.}~~\|\ze\|_{\ups_1,p,1}\le\de(b),$$
the isomorphisms $R_{\ups,\ze}$ of Subsection~\ref{gluing0_subs} corresponding
to the restriction of $b(\ups_1,\ze)$ to $\Si_{\ups_1,h}$, with $h\!\in\!I_1$,
induce an isomorphism 
$$R_{\ups_1,\ze}\!: \ker\bpar_{\na,b}\lra\ker\bpar_{\na,b(\ups_1,\ze)}$$
such that 
$$\big\|R_{\ups_1,\ze}\xi\big\|_{\ups_1,p,1}\le 2\|\xi\|_{b,p,1}.$$
Furthermore, by the $m\!=\!1$ case of~\e_ref{derivexp_e1} and~\e_ref{derivexp_e2},
\begin{equation}\label{g1reg_lmm1e4}
\Big|\D_{b(\ups_1,\ze),h}R_{\ups_1,\ze}\xi-\!\!
\sum_{i\in\chi(b),h(i)=h}\!\!\!\!\!\!\!\!\!\ti\rho_i(\ups)\D_{b,i}\xi\Big|
\le C(b)\big(|\ups|^{1/p}\!+\!\|\ze\|_{\ups_1,p,1}\big)\!
\sum_{i\in\chi(b),h(i)=h}\!\!\!\!\!\!\!\!\!|\ti\rho_i(\ups)|\cdot\|\xi\|_{b,p,1},
\end{equation}
for all $h\!\in\!I_1$ and $\xi\!\in\!\ker\bpar_{\na,b}$.
Let $\na^{\ups_1,\ze}$ and $g_{\ups_1,\ze}$ be the connection and metric in the line
bundle $u_{\ups_1,\ze}^*\L$ induced by the connections and metrics of Lemma~\ref{conn_lmm}.
For each $h\!\in\!I_1$ and $\de\!\in\!\R^+$, let
$$A_{\ups_1,h}(\de)= q_{\ups}^{-1}\big(A_{b,h}(\de)\big).$$
From the estimates \e_ref{derivexp_e1} and~\e_ref{derivexp_e2}, we find that
\begin{equation}\label{derivest_e3}
\big\|\na^{\ups_1,\ze}R_{\ups_1,\ze}\xi
 \big\|_{C^0( A_{\ups_1,h}(\de(b))),g_{\ups_1}} \le C(b)
 \big(|\ups|^{1/p}\!+\!\|\ze\|_{\ups_1,p,1}\big)
 \sum_{i\in\chi(b),h(i)=h}\!\!\!\!\!\!\!\!\!|\ti\rho_i(\ups)|\cdot\|\xi\|_{b,p,1},
\end{equation}
for all $h\!\in\!I_1$ and $\xi\!\in\!\ker\bpar_{\na,b}$.\\

\noindent
Fix a smooth function $\ep\!:\ti\X_{\T;B}(X)\!\lra\!\R^+$ such that for every $b\!\in\!\ti\X_{\T;B}(X)$ and $h\!\in\!I_1$ the disk of radius of $8\ep_b$
in $\Si_{b;P}$ around the node $x_h(b)$ contains no other special, i.e.~singular or
marked, point of~$\Si_b$.
For each 
$$\ups\!\equiv\!(b,v)\!\equiv\!\big(b,(v_h)_{h\in\aleph\sqcup\hat{I}} 
\in\wt\cF^{\eset}$$ 
sufficiently small, let
$$q_{\ups_0;2}\!:\Si_{\ups}\lra\Si_{\ups_1}
\qquad\hbox{and}\qquad
\ti{q}_{\ups_0;2}\!:\Si_{\ups}\lra\Si_{\ups_1}$$
be the basic gluing map of Subsection~2.2 in~\cite{gluing} corresponding 
to the gluing parameter~$v_0$ and the modified basic gluing map defined
in the middle of Subsection~\ref{g1comp-reg1_subs2} in~\cite{g1comp}
with the collapsing radius~$\ep_b$.
In this case, $\Si_{\ups}$ is a smooth genus-one curve.
For each $h\!\in\!I_1$, the maps $q_{\ups_0;2}$ and $\ti{q}_{\ups_0;2}$
collapse the circles of radii $|v_h|^{1/2}$ and $\ep_b$, respectively, 
around the point $x_h(b)\!\in\!\Si_{\ups_1;P}$. 
As before, the~map 
$$q_{\ups} \equiv q_{\ups_0;2}\!\circ q_{\ups_1}\!: \Si_{\ups}\lra\Si_b$$ 
induces a metric $g_{\ups}$ on~$\Si_{\ups}$ such that $(\Si_{\ups},g_{\ups})$ is 
obtained from $\Si_b$ by replacing all nodes by thin necks.
Let
$$I_1^*=I_1-\big\{h\!\in\!I_1\!: A_i\!=\!0~\forall\, i\!\ge\!h\big\}.$$
The map $\ti{q}_{\ups_0;2}$ is biholomorphic outside $|\aleph|$
thin necks $A_{\ups,h}$, with $h\!\in\!\aleph$, of $(\Si_{\ups},g_{\ups})$
and the $|I_1|$ annuli 
$$\ti\A_{b,h} \equiv \ti\A_{b,h}^-\cup\ti\A_{b,h}^+,$$
with $h\!\in\!I_1$, where
$$\ti\A_{b,h}^{\pm} \equiv \ti\A_{b,h}^{\pm}\big(\de(b)\big)
\subset \Si_{b;P} \!\approx\! \Si_{\ups}$$
are annuli independent of~$\ups$.
In addition,
\begin{gather}\label{glmap_e1a}\begin{split}
&\qquad\qquad\qquad\qquad 
u_{\ups_1,\ze}|_{\ti{q}_{\ups_0;2}(A_{\ups,h})}=\const 
~~\forall\,h\!\in\!\aleph,\\
&u_{\ups_1,\ze}|_{\ti{q}_{\ups_0;2}(\ti\A_{b,h})}=\const 
~~\forall\,h\!\in\!I_1\!-\!I_1^*, \qquad
u_{\ups_1,\ze}|_{\ti{q}_{\ups_0;2}(\ti\A_{b,h}^+)}=\const 
~~\forall\,h\!\in\!I_1^*;
\end{split}\\
\label{glmap_e1b}
\ti{q}_{\ups_0;2}(\ti\A_{b,h}^-) \subset 
  A_{\ups_1,h}\big(|v_h|^2/\de(b)\big) 
  ~~~\hbox{and}~~~
\big\|d\ti{q}_{\ups_0;2}\big\|_{C^0(\ti\A_{b,h}^-)}
\le C(b)|v_h| \quad\forall\, h\!\in\!I_1^*,
\end{gather}
if the $C^0$-norm of $d\ti{q}_{\ups_0;2}$ is computed
 with respect to the metrics $g_{\ups}$ on $\Si_{\ups}$
and $g_{\ups_1}$ on~$\Si_{\ups_1}$.
Furthermore, 
\begin{equation}\label{glmap_e2}
\big\|d\ti{q}_{\ups_0;2}\big\|_{C^0} \le C(b).
\end{equation}\\

\noindent
If $\ups\!=\!(b,v)\!\in\!\wt\cF^{\eset}$ is sufficiently small and
$\ze\!\in\!\Ga_B(\ups_1)$, we put
$$u_{\ups,\ze}= u_{\ups_1,\ze}\circ \ti{q}_{\ups_0;2} \qquad\hbox{and}\qquad
b(\ups,\ze)=\big(\Si_{\ups};u_{\ups,\ze}\big).$$
The map~$q_{\ups_0;2}$ induces weighted $L^p_1$-norms $\|\cdot\|_{\ups,p,1}$ on the spaces
$$\Ga(\Si_{\ups};u_{\ups,\ze}^*TX) \qquad\hbox{and}\qquad
\Ga(\Si_{\ups};u_{\ups,\ze}^*\L)$$
and a weighted $L^p$-norm $\|\cdot\|_{\ups,p}$ on the space 
$$\Ga(\Si_{\ups};\La^{0,1}_{\I,j}T^*\Si_{\ups}\!\otimes\!u_{\ups,\ze}^*TX).$$
Let $\Ga(\ups,\ze)$, $\Ga(\ups,\ze;\L)$,  and~$\Ga^{0,1}(\ups,\ze;\L)$ 
be the corresponding completions.
We~put
$$\Ga_-(\ups,\ze;\L)
=\big\{R_{\ups,\ze}\xi\!\equiv\!R_{\ups_1,\ze}\xi\!\circ\!\ti{q}_{\ups_0;2}\!: 
\xi\!\in\!\ker\bpar_{\na,b}\big\} \subset\Ga(\ups,\ze;\L).$$
By \e_ref{derivest_e3}-\e_ref{glmap_e2}, for all $\xi\in\ker\bpar_{\na,b}$
\begin{equation}\label{toruslmm_e4}
\big\|\bpar_{\na,b(\ups,\ze)}R_{\ups,\ze}\xi\big\|_{\ups,p}
\le C(b)\sum_{h\in I_1^*}~|v_h|\!\!\!\!\!\!\!\!\!\!
\sum_{i\in\chi(b),h(i)=h}\!\!\!\!\!\!\!\!\!|\ti\rho_i(\ups)|\cdot\|\xi\|_{b,p,1}
=C(b)\big|\rho(\ups)\big|\cdot\|\xi\|_{b,p,1}.
\end{equation}
Let $\Ga_+(\ups,\ze;\L)$ be the $L^2$-orthogonal complement of
$\Ga_-(\ups,\ze;\L)$ in $\Ga(\ups,\ze;\L)$.
We denote~by
$$\pi_{\ups,\ze;-}\!: \Ga(\ups,\ze;\L)\lra \Ga_-(\ups,\ze;\L)$$
the $L^2$-projection map.
For the same reasons as before,
\begin{equation}\label{g1reg_lmm1e5}
C(b)^{-1}\big\|\xi\big\|_{\ups,p,1} \le  \big\|\bpar_{\na,b(\ups,\ze)}\xi\big\|_{\ups,p} 
\le C(b)\big\|\xi\big\|_{\ups,p,1}
\qquad\forall~\xi\in\Ga_+(\ups,\ze;\L),
\end{equation}
if $\ups\!=\!(b,v)\!\in\!\cF^{\eset}$ and $\ze\!\in\!\Ga(\ups_1)$ are sufficiently small.
Let $\Ga_+^{0,1}(\ups,\ze;\L)$ be the image of $\Ga_+(\ups,\ze;\L)$ 
under~$\bpar_{\na,b(\ups,\ze)}$.\\

\noindent
The operator $\bpar_{\na,b_{\ups,\ze}}$ is not surjective.
We next describe its cokernel. 
Since the operator $\bpar_{\na,b}^B$ is surjective, the cokernel of $\bpar_{\na,b}$
can be identified with the vector space
$$\Ga^{0,1}_-(b;\L) \equiv \H_{b;P}\!\otimes\!\L_{\ev_P(b)}
\approx\E_{\pi_P(b)}^*\!\otimes\!\L_{\ev_P(b)},$$
where $\H_{b;P}$ is the space of harmonic antilinear differentials on 
the main component $\Si_{b;P}$ of~$\Si_b$.
As in Subsection~\ref{g1cone-gluing_subs3} in~\cite{g1cone},
there exist isomorphisms
$$R^{0,1}_{\ups;P}\!: \H_{b;P} \lra  \H_{\ups;P}\!\equiv\!\H_{b(\ups,\ze);P},
\qquad \ups\!=\!(b,v)\in \wt\cF_{\de},$$
such that the family of induced homomorphisms
$$ \H_{b;P}\lra\Ga^{0,1}(\ups;\C)^*,
\qquad \{R^{0,1}_{\ups;P}\eta\}(\eta')=\llrr{R^{0,1}_{\ups;P}\eta,\eta'}_2
\quad\forall\,\eta\!\in\!\H_{b;P},\,\eta'\!\in\!\Ga^{0,1}(\ups;\C),$$
is continuous on~$\wt\cF_{\de}$, and
\begin{equation}\label{glmap_e6}
R^{0,1}_{\ups;P}|_b=\id \qquad\forall\, b\!\in\!\ti\X_{\T;B}(X).
\end{equation}\\

\noindent
Let $\be\!:\R^+\!\lra\![0,1]$ be a smooth function such that
$$\be(t)\in
\begin{cases}
0,&\hbox{if}~t\!\le\!1;\\
1,&\hbox{if}~t\!\ge\!2.\\
\end{cases}$$
If $r\!\in\!\R^+$, let $\be_r(t)\!=\!\be(t/\sqrt{r})$.
We define $\be_b\!\in\!C^{\i}(\Si_b;\R)$ by
$$\be_b(z)=
\begin{cases}
1,&\hbox{if}~z\!\in\!\Si_{b,i},~i\!\in\!\chi^0(b);\\
1-\be_{\de(b)}(r(z)),&\hbox{if}~z\!\in\!\Si_{b,i},~i\!\in\!\chi(b);\\
0,&\hbox{otherwise},
\end{cases}$$
where $r(z)\!=\!|q_S^{-1}(z)|$ if $q_S\!:\C\!\lra\!S^2$ is
the stereographic projection mapping the origin to the south pole of~$S^2$.
In other words, $\be_b\!=\!1$ on $\Si_b^0(\de(b)/2)$ and vanishes
outside of $\Si_b^0(2\de(b))\!\subset\!\Si_b$.
Let $\be_{\ups}\!=\!\be_b\!\circ\!q_{\ups}$.
If $z\!\in\!\Si_{\ups}^0(2\de(b))$,
we denote by $\Pi_z^{(\ups,\ze)}$ the parallel transport in the line bundle 
$u_{\ups,\ze}^*\L$ along a path from 
$x\!\in\!\ti{q}_{\ups_0;2}^{-1}(\Si_{\ups_1;P})$ 
to $z$ in $\Si_{\ups}^0(2\de(b))$  with respect to 
the connection $\ti{q}_{\ups_0;2}^*\na^{(\ups_1,\ze)}$.
For each 
\begin{equation}\label{upseta_e}
\ups=(b,v)\!\in\!\wt\cF_{\de}^{\eset} \qquad\hbox{and}\qquad
\eta\!\in\!\Ga^{0,1}_-(b;\L),
\end{equation}
let $R_{\ups,\ze}^{0,1}\eta\!\in\!\Ga^{0,1}(\ups,\ze;\L)$ be given by
$$\{R_{\ups,\ze}^{0,1}\eta\}_zw
= \be_{\ups}(z) \, \Pi_z^{\ups,\ze}\eta_z(w) \in \L_{u_{\ups,\ze}(z)}
\qquad\,z\!\in\!\Si_{\ups},\,w\!\in\!T_z\Si_{\ups}.$$
Since the curvature of $\ti{q}_{\ups_0;2}^*\na^{(\ups_1,\ze)}$ vanishes
over $\Si_{\ups}^0(2\de(b))$, $\{R_{\ups,\ze}^{0,1}\eta\}_zw$ is independent
of the choice of~$x$ and path from $x$ to~$z$ above.\\

\noindent
If $\eta\!\in\!\Ga_-^{0,1}(b;\L)$,  we put
$$\|\eta\|=\sum_{h\in I_1^*}|\eta|_{x_h(b)},$$
where $|\eta|_{x_h(b)}$ is the norm of $\eta|_{x_h(b)}$ with respect to
the metric $g_{\pi_P(b)}$ on $\Si_{b;P}$.
If $\ups$ and $\eta$ are as in~\e_ref{upseta_e} and $\|\eta\|\!=\!1$,
we define~by
$$\pi^{0,1}_{\ups,\ze;-}\!: \Ga^{0,1}(\ups,\ze;\L)\lra\Ga_-^{0,1}(b;\L)
\qquad\hbox{by}\qquad
\pi^{0,1}_{\ups,\ze;-}(\eta')=\llrr{\eta',R_{\ups,\ze}^{0,1}\eta}_2\eta
\quad\forall\,\eta'\!\in\!\Ga^{0,1}(\ups,\ze;\L).$$
Since the space $\Ga_-^{0,1}(b;\L)$ is one-dimensional, 
$\pi^{0,1}_{\ups,\ze;-}$ is independent of the choice of~$\eta$.
We note that since $p\!>\!2$, by Holder's inequality
\begin{equation}\label{toruslmm_e5a}
\big\|\pi^{0,1}_{\ups,\ze;-}\eta'\big\|\le C(b)\|\eta'\|_{\ups,p}
\qquad\forall\,\eta'\!\in\!\Ga^{0,1}(\ups,\ze;\L).
\end{equation}
Furthermore, by the proof of Lemma~2.2 in~\cite{gluing},
\begin{equation}\label{toruslmm_e5b}
\big\|\pi^{0,1}_{\ups,\ze;-}\bpar_{\na,b(\ups,\ze)}\xi\big\|
\le C(b)\big|\rho(\ups)\big|\|\xi\|_{\ups,p,1}
\qquad\forall\,\xi\!\in\!\Ga(\ups,\ze;\L).
\end{equation}
With the same restriction on the homomorphisms $R_{\ups;P}^{0,1}$ and
identification of gluing parameters as described in Subsection~\ref{g1comp-reg1_subs2} 
of~\cite{g1comp}, we also have
\begin{equation}\label{toruslmm_e5c}
\pi^{0,1}_{\ups,\ze;-}\bpar_{\na,b(\ups,\ze)}R_{\ups,\ze}\xi
=-2\pi\I\, \D_{\rho(\ups)}\xi
\qquad\forall\,\xi\!\in\!\ker\bpar_{\na,b},
\end{equation}
by the proof of Proposition~4.4 in~\cite{gluing}.\\

\noindent
For each  $\ups\!=\!(b,v)\!\in\!\cF_{\de}^{\eset}$, $\ze\!\in\!\Ga_B(\ups_1)$,
and $\ze'\!\in\!\Ga(\ups,\ze)$, we~put
$$u_{\ups,\ze,\ze'}=\exp_{u_{\ups,\ze}}\!\ze' \qquad\hbox{and}\qquad
b(\ups,\ze,\ze')=\big(\Si_{\ups};u_{\ups,\ze,\ze'}\big).$$
We denote by $\Pi_{\ze'}$ the isomorphisms 
\begin{gather*}
\Ga(\Si_{\ups};u_{\ups,\ze}^*\L) \lra \Ga(\Si_{\ups};u_{\ups,\ze,\ze'}^*\L)
\qquad\hbox{and}\\
\Ga(\Si_{\ups};\La^{0,1}_{\I,j}T^*\Si_{\ups}\!\otimes\!u_{\ups,\ze}^*\L\big)
\lra \Ga(\Si_{\ups};\La^{0,1}_{\I,j}T^*\Si_{\ups}\!\otimes\!u_{\ups,\ze,\ze'}^*\L\big)
\end{gather*}
induced by the $\na$-parallel transport along the geodesics 
$\tau\!\lra\!\exp_{u_{\ups,\ze}}\!\tau\ze'$ with $\tau\!\in\![0,1]$.
Let 
$$L_{\ups,\ze,\ze'}=\Pi_{\ze'}^{-1}\bpar_{\na,b(\ups,\ze,\ze')}\Pi_{\ze'}
-\bpar_{\na,b(\ups,\ze)}\!:
\Ga(\ups,\ze;\L) \lra \Ga^{0,1}(\ups,\ze;\L\big).$$
Similarly to \e_ref{g0reg_lmm1e2},
\begin{equation}\label{g1reg_lmm2e7}
\big\|L_{\ups,\ze,\ze'}\xi\big\|_{\ups,p} \le C(b)\|\ze'\|_{\ups,p,1}^2\|\xi\|_{\ups,p,1}
\qquad\forall~\xi\!\in\!\Ga(\ups,\ze;\L).
\end{equation}\\

\noindent
Let $J$ be an almost complex structure on $X$.
With notation as at the beginning of this subsection and in Subsection~\ref{notation1_subs},
we define the linear bundle map 
$$\cD_{J;\T}\!: \wt\F \lra \pi_P^*\E^*\!\otimes_{\C}\!\ev_P^*(TX,J)$$
over $\ti\X_{\T;B}(X)$ by
$$\big\{\cD_{J;\T}\big(b,(\ti{v}_i)_{i\in\chi(b)}\big)\big\}(\psi)
= \sum_{i\in\chi(b)}\!\!\!\psi_{x_{h(i)}(b)}(\ti{v}_i)\cdot_J\!\cD_ib \in T_{\ev_P(b)}X.$$
Suppose $t_r$ and $b_r$ are as in Lemma~\ref{torus_lmm} and 
$\Si_{b_r;P}\!=\Si_{b_r}$ is a smooth torus.
Since the sequence $[b_r]$ converges to $[b]$, 
by Subsection~\ref{g1comp2-comp1prp_subs} in~\cite{g1comp2} 
there exist $C\!\in\!\R^+$ and for all $r$ sufficiently large 
\begin{equation}\label{torus_e1}
b_r'\in\ti\X_{\T;B}(X), \quad \ups_r=(b_r',v_r)\in\wt\cF_{\de}^{\eset},  \quad
\ze_r\in\Ga(\ups_r), \quad \ze_r'\in\Ga(\ups_r,\ze_r), \quad\hbox{and}\quad \ep_r\in\R^+
\end{equation}
such that
\begin{gather}\label{torus_e1a}
\lim_{r\lra\i}b_r'=b, \quad \lim_{r\lra\i}|v_r|=0, \quad
\lim_{r\lra\i}\|\ze_r\|_{\ups_{r;1},p,1}=0,  \quad
\|\ze_r'\|_{\ups_r,p,1}\le C\big|\rho(\ups)\big|,\\
\label{torus_e1b}
\lim_{r\lra\i}|\ep_r|=0, \qquad 
\big|\cD_{J_0;\T}\rho(\ups)\big|\le \ep_r\big|\rho(\ups)\big|,\\
\hbox{and}\qquad b_r\!\equiv\!\big(\Si_{b_r};u_{b_r}\big)=
\big(\Si_{\ups_r},\exp_{u_{\ups_r,\ze_r}}\!\ze_r'\big).\notag
\end{gather}
The last equality holds for a representative $b_r$ for $[b_r]$.\\

\noindent
By~\e_ref{g1reg_lmm2e7} and the last inequality in~\e_ref{torus_e1a},
for some $C\!\in\!\R^+$ and for all $r$ sufficiently large
\begin{equation}\label{toruslmm_e21}
\big\|L_{\ups_r,\ze_r,\ze_r'}\xi\big\|_{\ups_r,p} \le 
C\big|\rho(\ups_r)\big|^2\|\xi\|_{\ups_r,p,1}
\qquad\forall~\xi\!\in\!\Ga(\ups_r,\ze_r;\L).
\end{equation}
Thus, by \e_ref{toruslmm_e4} and~\e_ref{g1reg_lmm1e5},
\begin{equation}\label{toruslmm_e22}
\big\|\pi_{\ups_r,\ze_r;-}\Pi_{\ze_r'}^{-1}\xi-\Pi_{\ze_r'}^{-1}\xi
\big\|_{\ups_r,p,1}\le C\big|\rho(\ups_r)\big|\|\xi\|_{\ups_r,p,1}
\qquad\forall~\xi\!\in\!\ker\bpar_{\na,b_r}.
\end{equation}
Since 
$$\big\{\bpar_{\na,b(\ups_r,\ze_r)}\!+\!L_{\ups_r,\ze_r,\ze_r'}\big\}
\Pi_{\ze_r'}^{-1}\xi=0 \qquad\forall~\xi\!\in\!\ker\bpar_{\na,b_r},$$
by~\e_ref{toruslmm_e5a}-\e_ref{toruslmm_e5c}, \e_ref{toruslmm_e21}, 
and~\e_ref{toruslmm_e22},
\begin{equation}\label{toruslmm_e25}
\big|\D_{\rho(\ups)}\xi\big|\le C\big|\rho(\ups_r)\big|^2\cdot\|\xi\|_{b_r',p,1}
\qquad\forall\,R_{\ups_r,\ze_r}\xi\in \pi_{\ups_r,\ze_r;-}\Pi_{\ze_r'}^{-1}\ker\bpar_{\na,b_r}.
\end{equation}\\

\noindent
After passing to a subsequence if necessary, let 
$$\under{w}'\!\equiv\!(w_i')_{i\in\chi(b)} \in \wt\F_b-\{0\}$$
be such that
$$\lim_{r\lra\i}\!\big[\rho(\ups_r)\big]=\big[\under{w}'\big]
\in \bP\wt\F_b.$$
Since $\cD_ib_r'\!\lra\!\cD_ib$ for all $i\!\in\!\chi(b)$,
\begin{equation}\label{toruslmm_e27a}
\big\{\cD_{J_0;\T}\big(b,(w_i')_{i\in\chi(b)}\big)\big\}(\psi)
\equiv \sum_{i\in\chi(b)}\!\!\!\psi_{x_{h(i)}(b)}(w_i')\cdot_{J_0}\!\cD_ib 
=0 \in T_{\ev_P(b)}X \quad\forall\,\psi\!\in\!\E_{\pi_P(b)},
\end{equation}
by~\e_ref{torus_e1b}.
If $\xi_r\!\in\!\ker\bpar_{\na,b_r}$ and $\xi\!\in\!\ker\bpar_{\na,b}$
are such that
$$\lim_{r\lra\i}\![\xi_r]=[\xi] \in \V_{1,M}^A|_{[b]},$$
then by \e_ref{toruslmm_e22} and~\e_ref{toruslmm_e25},
\begin{equation}\label{toruslmm_e27b}
\big|\{\D_{\under{w}'}\xi\}\psi\big| \equiv 
\bigg|\sum_{i\in\chi(b)}\!\!\!\psi_{x_{h(i)}(b)}(w_i')\cdot\cD_{b,i}\xi\bigg|
\le \ti\ep_r\big|\under{w}'\big|\|\psi\| 
\qquad\forall\,\psi\in\E_{\pi_P(b)}
\end{equation}
for a sequence $\ti\ep_r$ converging to zero.
Thus, by~\e_ref{toruslmm_e27a} and~\e_ref{toruslmm_e27b}, 
for the purposes of Lemma~\ref{torus_lmm} we can take
$$w_i=\psi_{x_{h(i)}(b)}(w_i')\in\C^*,$$
where $\psi\!\in\!\E_{\pi_P(b)}$ is any nonzero element.\\

\noindent
{\it Remark~1:} If $\Si_{b_r;P}\!=\!\Si_{b_r}$ is a circle of spheres,
i.e.~$\aleph'\!\neq\!\eset$ in the notation of Subsection~\ref{gluingsumm_subs},
the proof of Lemma~\ref{torus_lmm} is formally the same, but some details change
in a way analogous to Subsection~3.9 in~\cite{gluing}.
In particular, in~\e_ref{torus_e1},
$$\ups_r\in\wt\cF^{\aleph_0}\equiv\big\{
\big(b,(v_h)_{h\in\aleph\sqcup\hat{I}}\big)\in\wt\cF\!: 
v_h\!=\!0\,\Llra\, h\!\in\!\aleph_0\big\},$$
for a nonempty subset $\aleph_0$ of $\aleph$.
If $\ups\!\in\!\wt\cF_{\de}^{\aleph_0}$, $\Si_{\ups}$ is a circle of spheres
with nodes~$\aleph_0$.
If in addition $\ze\!\in\!\Ga(\ups_1)$, $\Ga(\ups,\ze)$ consists of the vector fields
on the $|\aleph_0|$ components of $\Si_{\ups}$ that agree at the nodes of~$\Si_{\ups}$.
Similarly, $\Ga(\ups,\ze;\L)$ consists of the sections of $u_{\ups,\ze}^*\L$ over
the components of $\Si_{\ups}$ that agree at the nodes.
If $\eta\!\in\!\Ga^{0,1}_-(b;\L)$, the $u_{\ups,\ze}^*\L$-valued $(0,1)$-form
$R_{\ups,\ze}^{0,1}\eta$ has poles at the nodes of $\Si_{\ups}$ with residues that
add~up to zero at each node.
In particular, $R_{\ups,\ze}^{0,1}\eta$ is not an element of $\Ga^{0,1}(\ups,\ze;\L)$, 
but the homomorphism $\pi^{0,1}_{\ups,\ze;-}$ is well defined and still satisfies
\e_ref{toruslmm_e5a}-\e_ref{toruslmm_e5c}.
Finally, the argument of Subsection~\ref{g1comp2-comp1prp_subs} in~\cite{g1comp2}
easily generalizes to show that $(\ups_r,\ze_r,\ze_r',\ep_r)$ as in
\e_ref{torus_e1}-\e_ref{torus_e1b} exist in this situation.\\

\noindent
{\it Remark~2:} A regularization $\vph_b$ of the cone ${\cal V}_{1,M}^A$ near $b$, 
as in Definition~\ref{coneprop_dfn}, can be constructed using
the description of open subsets in $\X_{1,M}(X,A)$ of Section~3 in~\cite{LT1}
and the corresponding analogues of the isomorphisms~$R_{\ups,\ze}$ 
and the injective homomorphisms $\pi_{\ups,\ze;-}\Pi_{\ze'}^{-1}$ as above.\\

\vspace{.2in}

{\it 
\begin{tabbing}
${}\qquad$
\= Department of Mathematics, Stanford University, Stanford,
CA 94305-2125\\ 
\> jli@math.stanford.edu, azinger@math.stanford.edu
\end{tabbing}}

\end{document}